\input amstex.tex
\documentstyle{amsppt}
\UseAMSsymbols\TagsAsMath\widestnumber\key{ASSSS}
\magnification=\magstephalf\pagewidth{6.2in}\vsize8.0in
\parindent=6mm\parskip=3pt\baselineskip=16pt\tolerance=10000\hbadness=500
\NoRunningHeads\loadbold\NoBlackBoxes\nologo
\def\today{\ifcase\month\or January\or February\or March\or April\or May\or June\or
     July\or August\or September\or October\or November\or December\fi \space\number\day, \number\year}

\def\ns#1#2{\mskip -#1.#2mu}\def\operator#1{\ns07\operatorname{#1}\ns12}
\def\pa{\partial}\def\na{\nabla\!}
\def\dt{{\ns18{\Cal D}_{\ns23 t\!} }}
\def\tr{\operatorname{tr}}
\def\div{\operator{div}}\def\curl{\operator{curl}}
\def\dist{\operatorname{dist} }
\def\leftalignspace{\mskip -45mu}\def\rightalignspace{\mskip -90mu}

\def\dooot#1{\dddot{{#1}}^{{}^{}}}\def\doooot#1{\ddddot{{#1}}^{{}^{}}}

\topmatter
\title
Well-posedness for the motion of an incompressible liquid with
free surface boundary
\endtitle
\author  Hans Lindblad
\endauthor
\thanks  The author was supported in part by the National
Science Foundation.  \endthanks
\address
University of California at San Diego
\endaddress
\email
lindblad\@math.ucsd.edu
\endemail

\abstract We study the motion of an incompressible perfect liquid
body in vacuum. This can be thought of as a model for the motion
of the ocean or a star. The free surface moves with the velocity
of the liquid and the pressure vanishes on the free surface. This
leads to a free boundary problem for Euler's equations, where the
regularity of the boundary enters to highest order. We prove local
existence in Sobolev spaces assuming a "physical condition",
related to the fact that the pressure of a fluid has to be
positive.

\endabstract
\endtopmatter
\document

\head{1. Introduction}\endhead
We consider Euler's equations
describing the motion of a perfect incompressible fluid in vacuum:
$$\align
\big(\partial_t +V^k\partial_{k} \big) v_j +\partial_{j} p&=0,\quad j=1,...,n
\quad \text{in}\quad {\Cal D},\tag 1.1\\
\div V &=\partial_{k} V^k=0 \quad \text{in}\quad {\Cal D}\tag 1.2
\endalign
$$
where $\partial_i=\pa/\pa x^i$ and ${\Cal D}=\cup_{\, 0\leq t\leq T}\,
\{t\}\times {\Cal D}_t $,
${\Cal D}_t\subset \Bbb R^n$.
Here  $V^k\!\ns15=\delta^{ki} v_i=v_k$,
and we use the convention that repeated upper and lower indices
are summed over. $V\!$ is the velocity vector field of the fluid,
$p$ is the pressure and
${\Cal D}_t$ is the domain the fluid occupies at time $t$.
We also require boundary conditions on the free boundary
$\partial {\Cal D}=\cup_{\, 0\leq t\leq T}\, \{t\}\times\partial {\Cal D}_t$;
$$
\align p=0,\quad &\text{on}\quad \partial {\Cal D},\tag 1.3\\
       (\partial_t+V^k\partial_{k})|_{\partial {\Cal D}}
&\in T(\partial {\Cal D}),\tag 1.4
\endalign
$$
Condition (1.3) says that the pressure $p$ vanishes outside the domain
and condition (1.4) says that the boundary moves with the velocity $V$
of the fluid particles at the boundary.

Given a domain
${\Cal D}_0\subset \Bbb R^n$, that is homeomorphic to the unit ball,
 and initial data $v_0$, satisfying the constraint
(1.2), we want to find a set
${\Cal D}=\cup_{\, 0\leq t\leq T}\, \{t\}\times{\Cal D}_t$,
${\Cal D}_t\subset{\Bbb R}^n$
and a vector field $v$ solving (1.1)-(1.4) with initial conditions
$$
 \{x;\, (0,x)\in {\Cal D} \}={\Cal D}_0,\qquad\text{and}\qquad
v=v_0,\quad\text{on}\quad \{0\}\times {\Cal D}_0 \tag 1.5
$$
Let ${\Cal N}$ be the exterior unit normal to the
free surface $\pa{\Cal D}_t$.
Christodoulou\cite{C2} conjectured
 that the initial value problem (1.1)-(1.5),
is well posed in Sobolev spaces if
$$
\na_{\Cal N}\, p\leq -c_0 <0,\quad\text{on}\quad \partial {\Cal D},
\qquad \text{where}\quad \na_{\Cal N}={\Cal N}^i\partial_{x^i}. \tag 1.6
$$

Condition (1.6) is a natural {\it physical condition} since the pressure
$p$ has to be positive in the interior of the fluid.
It is essential for the well posedness in Sobolev spaces.
A condition related to Rayleigh-Taylor instability
in \cite{BHL,W1} turns out to be equivalent to (1.6), see \cite{W2}.
Taking the divergence of (1.1) gives:
$$
-\triangle p= (\partial_{j} V^k)\partial_{k} V^j,
\qquad \text{in}\quad {\Cal D}_t,\qquad p=0,\quad \text{on}\quad\partial
{\Cal D}_t\tag 1.7
$$
In the irrotational case, when
$\curl {\,\ns18}v_{\ns10 \, ij}\!=\!\partial_i v_j\!-\partial_j v_i\!=\!0$, then
 $\triangle p\!\leq\!0$ so $p\!\geq\!0$ and (1.6) holds by the strong
maximum principle. Furthermore
Ebin \cite{E1} showed that the equations are ill posed when (1.6) is
not satisfied and the pressure is negative and
Ebin \cite{E2} announced an existence result when one adds
surface tension to the boundary condition which has a regularizing effect so
(1.6) is not needed then.

The incompressible perfect fluid is to be thought of as an idealization of a liquid.
For small bodies like water drops surface tension should help holding it together and for
larger denser bodies like stars its own gravity should play a role.
Here we neglect the influence of such forces. Instead it is the incompressibility condition that prevents the
body from expanding and it is the fact that the pressure is positive
that prevents the body from breaking up in the interior.
Let us also point out that, from
a physical point of view one can alternatively think of the pressure as being a small
positive constant on the boundary instead of vanishing.
What makes this problem difficult is that the regularity
of the boundary enters to highest order.
Roughly speaking, the velocity tells the boundary where to move and the boundary
is the zero set of the pressure that determines the acceleration.

In general it is possible to prove local existence for analytic data
for the free interface between two fluids.
However, this type of problem might be subject to instability in Sobolev norms,
in particular Rayleigh-Taylor instability, which occurs when a heavier fluid is
on top of a lighter fluid.
Condition (1.6) prevents Rayleigh-Taylor instability from occurring.
Indeed,  if condition (1.6) is violated  Rayleigh-Taylor instability occurs
in a linearized analysis.

Some existence results in Sobolev spaces were known in the irrotational case,
for the closely related water wave problem which describes
the motion of the surface of the ocean under the influence of earth's gravity.
The gravitational field can be considered as uniform
and it reduces to our problem by going to an
accelerated frame. The domain ${\Cal D}_t$ is unbounded for the water wave problem
coinciding with a half-space in the case of still water.
Nalimov\cite{Na} and Yosihara\cite{Y} proved local existence in Sobolev spaces
in two space dimensions for initial conditions sufficiently close to still water.
Beale, Hou and Lowengrab\cite{BHL} have given an argument to show
that problem is linearly well posed in a weak sense in Sobolev spaces,
assuming a condition, which can be shown to be equivalent to (1.6).
The condition (1.6) prevents the Rayleigh-Taylor instability
from occurring when the water wave turns over.
Finally Wu\cite{W1,2} proved local existence in the general
irrotational case in two and three dimensions
for the water wave problem.
The methods of proofs in these papers uses that the vector field is irrotational
to reduce to equations on the boundary  and
do not generalize to deal with the case of nonvanishing curl.

We consider the general case of nonvanishing curl.
With Christodoulou $\ns40$\cite{CL}$\ns40$ we proved local
{\it $\ns30$a$\ns30$ priori$\ns10$} bounds in Sobolev
spaces in the general case of non vanishing curl, assuming (1.6) hold initially.
Usually if one has {\it a priori} estimates, existence follows from
similar estimates for some regularization or iteration scheme for the equation,
but the sharp estimates in \cite{CL} use all the symmetries of the
equations and so only hold for perturbations of the equations that preserve the
symmetries. In \cite{L1} we proved existence for the linearized equations, but
the estimates for the solution of the linearized equations
looses regularity compared to  the solution we linearize around,
so existence for the nonlinear problem does not follow directly.
Here we use improvements of the estimates in \cite{L1} together with the
Nash-Moser technique to show local existence for the nonlinear problem in
the smooth class:

\proclaim{Theorem {1.}1} Suppose that $v_0$ and $\pa{\Cal D}_0$ in (1.5) are smooth,
${\Cal D}_0$ is diffeomorphic to the unit ball,
and that (1.6) hold initially when $t=0$. Then there is a $T>0$ such that
 (1.1)-(1.5) has a smooth solution for $0\leq t\leq T$,
and (1.6) hold with $c_0$ replaced by $c_0/2$  for $0\leq t\leq T$.
\endproclaim

In \cite{CL} we proved
local energy bounds in Sobolev spaces. It now follows from the bounds there
that the solution remains smooth as
long as it is $C^2$ and the physical condition (1.6) hold.
The existence for smooth data now implies existence in
the Sobolev spaces we considered in \cite{CL}.
Moreover, method here also works for the compressible case \cite{L2,L3}.

Let us now describe the main ideas and difficulties in the proof.
In order to construct an iteration scheme we must first introduce some parametrization
in which the moving domain becomes fixed and express Euler's equations
in this fixed domain. This is achieved by the Lagrangian coordinates
given by following the flow lines of the velocity vector field of the
fluid particles.

In \cite{L1} we studied the linearized equations
of Euler's equations expressed in Lagrangian coordinates.
We proved that the linearized operator is invertible
at a solution of Euler's equations. The linearized equations become an evolution equation
for what we called the normal operator, (2.17).
 The normal operator is unbounded
and not elliptic but it is symmetric and positive on divergence free vector
fields if (1.6) hold. This leads to energy bounds and
existence for the linearized equations follows from a delicate regularization argument.
The solution of the linearized equations however looses
regularity compared to the solution we linearize around so existence for the nonlinear problem
does not follow directly from an inverse function theorem
in a Banach space but we must use the Nash-Moser technique.

We first define a nonlinear functional
whose zero will be a solution of Euler's equations expressed in the Lagrangian coordinates.
Instead of defining our map by the left hand side
of (1.1) and (1.2) expressed in the Lagrangian coordinates
we let our map be given by the left hand side of (1.1)
and we let pressure be implicitly defined by (1.7)
 satisfying the boundary condition (1.3).
This is because one has to make sure that the pressure vanishes on the boundary
at each step of an iteration or else the linearized operator is
ill posed. One can see this by looking at the irrotational case
where one gets an evolution equation on the boundary.
If the pressure vanishes on the boundary then one has an evolution equation for
a positive elliptic operator but if it does not vanish on the boundary there
will also be some tangential derivative, no matter how small coefficients
they come with the equation will have exponentially growing Fourier modes.

In order to use the Nash-Moser technique one
has to be able to invert the linearized operator
in a neighborhood of a solution of Euler's equations
or at least do so up to a quadratic error \cite{Ha}.
In this paper we generalize the existence in \cite{L1} so the linearized
operator is invertible in a neighborhood of a solution of Euler's equations
and outside the class of divergence free vector fields.
This does present a difficulty because the normal operator, introduced in \cite{L1}
is only symmetric on divergence free vector fields and in general it looses regularity.
Overcoming this difficulty requires two new observations.
The first is that, also for the linearized equations there is an identity for the curl that
gives a bound that is better than expected.
The second is that one can bound any first order derivative
of a vector field by the curl,
the divergence and the normal operator times one over the constant $c_0$ in (1.6).
Although the normal operator is not elliptic on general vector fields it is elliptic on
irrotational divergence free vector fields and in general one can invert it
if one also have bounds for the curl and the divergence.

The methods here and in \cite{CL} are on a technical level very different but
there are philosophical similarities.
First we fix the boundary by introducing Lagrangian coordinates.
Secondly, we take the geometry of the boundary into account.
Here in terms of the normal operator and
 Lie derivatives with respect to tangential vector fields
and in \cite{CL} in terms of the second fundamental form of the boundary
and tangential components of the tensor of higher order derivatives.
Thirdly, we use interior estimates to pick up the curl and the divergence.
Lastly, we get rid of a difficult term,
the highest order derivative of the pressure, by projecting.
Here we use the orthogonal
projection onto divergence free vector fields whereas in \cite{CL}
we used the local projection of a tensor onto the tangent space of the boundary.

The paper is organized as follows.
In section 2 we reformulate the problem in the Lagrangian coordinates
and give the nonlinear functional which a solution of Euler's
equations is a zero of and  we derive the linearized equations in this formulation.
In section 2 we also give an outline of the proof and state the main steps
that we will prove.
The main part of the paper, sections 3 to 13 are devoted to proving existence
and tame energy estimates for the inverse of the linearized operator.
Once this is proven, the remaining sections 14 to 18 are devoted to setting
up the Nash-Moser theorem we are using.

\head{2. Lagrangian coordinates and the linearized operator}\endhead
%\subheading{Euler's equations expressed in Lagrangian coordinates}

Let us first introduce the Lagrangian coordinates in which the boundary becomes fixed.
By a scaling we may assume that ${\Cal D}_0$ has the volume of the unit ball
$\Omega$ and since we assumed that
${\Cal D}_0$ is diffeomorphic to the unit ball we can, by a theorem in \cite{DM},
find a volume preserving diffeomorphism $f_0:\Omega\to {\Cal D}_0$,
i.e. $\det{ \,(\pa f_0/\pa y)}=1$.
Assume that $v(t,x)$, $p(t,x)$, $(t,x)\in  {\Cal D}$ are given satisfying the boundary
conditions (1.3)-(1.4).
The Lagrangian coordinates $x=x(t,y)=f_t(y)$ are given by solving
$$
\frac{d x(t,y)}{dt}=V(t,x(t,y)), \qquad x(0,y)=f_0(y),\quad y \in \Omega\tag 2.1
$$
Then $f_t:\Omega\to {\Cal D}_t$ is a volume preserving diffeomorphism,
if $\div V=0$,
and the boundary becomes fixed in the new $y$ coordinates.
Let us introduce the material derivative:
$$
\quad D_t=\frac{\partial }{\partial t}\Big|_{y=constant}=
\frac{\partial }{\partial t}\Big|_{x=constant}
+\,  V^k\frac{\partial}{\partial x^k},\tag 2.2
$$
The partial derivatives
$\pa_i=\pa/\pa x^i$ can then be expressed in terms of partial derivatives
$\pa_a=\pa/\pa y^a$ in the Lagrangian coordinates.
We will use letters $a,b,c,...,f$ to denote partial differentiation in the
Lagrangian coordinates and $i,j,k,...$ to denote partial differentiation
in the Eulerian frame.

In  these coordinates Euler's equation (1.1) become
$$
D_t^2 x_i+\pa_i p=0,\quad \qquad(t,y)\in[0,T]\times\Omega,
\tag 2.3
$$
where now $x_i=x_i(t,y)$ and $p=p(t,y)$ are functions on $[0,T]\times \Omega$,
$D_t$ is just the partial derivative with respect to $t$
and $\pa_i=(\pa y^a/\pa x^i)\pa_a$, where $\pa_a$ is differentiation with respect to $y^a$.
(1.7) become
$$
\triangle p+(\pa_i V^k)\pa_k V^i=0,\qquad
p\Big|_{\pa \Omega}\!\!\!=0,\qquad\text{where}\quad  V^i=D_t x^i. \tag 2.4
$$
Here
$$
\triangle p\!=\!\sum_{i=1}^n \pa_i^2 p\!=\!
\kappa^{-1}\pa_a\big(\kappa g^{ab}\pa_b p\big)\qquad\text{where}\quad
g_{ab}=\delta_{ij}\frac{\pa x^i}{\pa y^a}\frac{\pa x^j}{\pa y^b}.\tag 2.5
$$
$g^{ab}$ is the inverse of the metric $g_{ab}$ and $\kappa=\det{(\pa x/\pa y)}
=\sqrt{\det{g}}$.
The initial conditions (1.5) become
$$
x\big|_{t=0}=f_0,\qquad D_{t\,} x\big|_{t=0}=v_0\tag 2.6
$$

 Christodoulou's physical condition (1.6) become
$$
\na_{\Cal N}\, p\leq -c_0 <0,\quad\text{on}\quad \partial {\Omega},
\qquad \text{where}\quad \na_{\Cal N}={\Cal N}^i\partial_{x^i}. \tag 2.7
$$
This is needed in the proof for the normal operator (2.17) to be positive
which leads to energy bounds.
In addition to (2.7) we also need to assume a coordinate condition having to do with
that we are looking for a solution in the Lagrangian coordinates and we are starting
by composing with a particular diffeomorphism. The coordinate condition is
$$
|\pa x/\pa y|^2+|\pa y/\pa x|^2\leq c_1^2 ,\qquad
\sum_{a,b=1}^n (|g^{ab}|+|g_{ab}|)\leq n c_1^2,\tag 2.8
$$
where $|\pa x/\pa y|^2=\sum_{i,a=1}^n (\pa x^i/\pa y^a)^2$.
This is needed for (2.5) to be invertible.
We note that the second condition in (2.8) follows from the first and the
first follows from the second with a larger constant.
We remark that this condition is fulfilled initially since we are composing with a
diffeomorphism. Furthermore, for solution of Euler's equations, $\div V=0$,
so the volume form $\kappa$ is preserved and hence an upper bound for the metric also implies
a lower bounded for the eigenvalues and an upper bound for the  inverse of
the metric follows. However, in the iteration, we will go outside the divergence free
class and hence we must make sure that both (2.7) and (2.8) hold
at each step of the iteration. We will prove the following theorem:
\proclaim{Theorem {2.}1} Suppose that initial data (2.6)
 are smooth, $v_0$ satisfy the constraint (1.2), and that (2.7) and (2.8)
hold when $t=0$. Then there is $T>0$
such that (2.3)-(2.4) has a solution
$x,p\in C^\infty([0,T]\times\overline{\Omega})$. Furthermore, (2.7)-(2.8) hold,
for $0\leq t\leq T$,
with $c_0$ replaced by $c_0/2$ and $c_1$ replaced by $2c_1$.
\endproclaim

Theorem {1.}1 follows from Theorem {2.}1. In fact, the assumption that ${\Cal
D}_0$ is diffeomorphic to the unit ball, together with that one
then can find a volume preserving diffeomorphism guarantees that
(2.8) hold initially. Once, we obtained a solution to
(2.3)-(2.4), we can hence follow the flow lines of $V$ in (2.1)
and this defines a diffeomorphism of $[0,T]\times \Omega$ to
${\Cal D}$, so we obtain smoothness of $V$ as a function of
$(t,x)$ from the smoothness as a function of $(t,y)$.

In this section we first define a nonlinear functional whose zero is a solution of
Euler's equations, (2.9)-(2.13). Then we derive the linearized
operator in Lemma {2.}2.
The existence will follow from the Nash-Moser inverse function theorem,
once we  proven that the linearized operator is invertible
 and so called tame estimates for the inverse stated in Theorem {2.}3.
Proving that the linearized operator is invertible away from a solution of
Euler's equations and
outside the divergence free class is the main difficulty of the paper.
This is because the normal operator (2.17) is only symmetric and positive
within the divergence free class and in general it looses regularity.
In order to prove that the linearized operator is invertible
 and estimates for its inverse we
introduce a modification (2.31) of the linearized operator that preserves
the divergence free condition, and first prove that the modification is
invertible and estimates for its inverse, stated in  Theorem {2.}4.
The difference between the linearized operator and the modification is
lower order and the estimates for the inverse of the
modified linearized operator leads to
existence and estimates also for the inverse of the linearized operator.

Proving the estimates for the inverse of the modified linearized operator,
stated in Theorem {2.}4, takes up most of the paper, sections 3 to 13.
In this section we also derive certain identities for the curl and
the divergence, see (2.29)-(2.30), need for the proof of Theorem {2.}4.
 Here we also transform the vector field to the Lagrangian frame
and express the operators and identities in the Lagrangian frame,
see Lemma {2}.5. The estimates in Theorem {2.}4 will be derived in
the Lagrangian frame since commutators of the normal
operator with certain differential operators are better behaved in this frame.

In section 3, we introduce the orthogonal projection onto divergence
free vector field and decompose the modified linearized
equation into a divergence free part
and an equation for the divergence.
This is needed to prove Theorem {2.}4
because the normal operator is only
symmetric on divergence free vector fields and in general
it looses regularity. However, we have a better equation for the
divergence which will allow us to obtain the same space regularity
for the divergence as for the vector field itself.

In section 4 we introduce the tangential vector fields and Lie derivatives and calculate
commutators between these and the operators that occur in the modified linearized equation,
in particular the normal operator.
In section 5 we show that any derivative of a vector field can be estimated by
 derivatives of the curl and of the divergence, and tangential derivatives or
tangential derivatives of the normal operator.
In section 6 introduce the $L^\infty$ norms that we will use and
state the interpolation inequalities that we will use.
In section 7 and 8 we give the tame $L^2\infty$ and $L^\infty$ estimates
for the Dirichlet problem.
In section 9 we give the equations and estimates for the curl that we will use.
In section 10 we show existence for the modified linearized equations
in the divergence class.
In section 11 we give the improved estimates for the inverse of the
modified linearized operator within the divergence free class.
These are needed in section 12 to prove existence and estimates for the
inverse of the modified linearized operator.
Finally in section 13 we use this to prove existence and estimates for the
inverse of the linearized operator.

In section 14 we explain what is needed to ensure
that the physical and coordinate conditions (2.7) and (2.8) continue to hold.
In section 15 we summarize the tame estimates for the inverse of the linearized
operator in the formulation that will use with the Nash-Moser theorem.
In section 16 we derive the tame estimates for the second variational derivative.
In section 17 we give the smoothing operators needed for the
proof of the Nash-Moser theorem on a bounded domain.
Finally, in section 18 we state and prove the Nash-Moser theorem in the
form that we will use.

%\subheading{The Euler map}
Let us now define the nonlinear map, that we will use to find a solution
of Euler's equations. Let
$$
\Phi_i=\Phi_i(x)=D_t^2 x_i +\pa_i p,\qquad\text{where}
\quad \pa_i=(\pa y^a/\pa x^i)\pa_a, \tag 2.9
$$
and $p=\Psi(x)$ is given by solving
$$
\triangle p=-(\pa_i V^k)\pa_k V^i,\qquad p\big|_{\pa\Omega}\!=0,
\qquad\text{where}\quad V=D_t x .
\tag 2.10
$$
A solution to Euler's equations is given by
$$
\Phi(x)=0,\quad\text{for}\quad 0\leq t\leq T, \qquad x\big|_{t=0}=f_0,\quad
D_t x\big|_{t=0}=v_0\tag 2.11
$$
We will find $T>0$ and a smooth function  $x$ satisfying (2.11) using
the Nash-Moser iteration scheme.

First we turn (2.11) into a problem with vanishing initial data and a small
inhomogeneous term using a trick from \cite{Ha} as follows.
It is easy, to construct a formal power series solution $x_0$ as $t\to 0$:
$$
D_t^k \Phi(x_0)\big|_{t=0}=0,\qquad k\geq 0,\qquad
x_0\big|_{t=0}=f_0,\quad D_t x_0\big|_{t=0}=v_0\tag 2.12
$$
In fact, the equation (2.10) for the pressure $p$ only depends on one time derivative
of the coordinate $x$ so commuting through time derivatives in (2.10) gives a
Dirichlet problem for $D_t^k p$ depending only on $D_t^m x$, for $m\leq k+1$ and
$D_t^{\ell}\, p$, for $\ell\leq k-1$. Similarly commuting through time derivatives in
Euler's equation, (2.11),
gives $D_t^{2+k}x$
in terms of $D_t^m x$, for $m\leq k$, and $D_t^\ell\, p$, for $\ell\leq k$.
We can hence construct a formal power series solution in $t$ at $t=0$
and by a standard trick we can find a smooth function $x_0$
having this as it power series, see section 10.
We will now solve for $u$ in
$$
\tilde{\Phi}(u)=\Phi(u+x_0)-\Phi(x_0)=F_\delta-F_0=f_\delta,
\qquad u\big|_{t=0}= D_t u\big|_{t=0}=0 \tag 2.13
$$
where $F_\delta$ is constructed as follows. Let $F_0=\Phi(x_0)$ and let
$F_\delta (t,y)=F_0(t-\delta,y)$, when $t\geq \delta$ and $F_\delta(t,y)=0$,
when $t\leq \delta$. Then $F_\delta$ is smooth and $f_\delta =F_\delta-F_0$ tends to $0$
in $C^\infty$ when $\delta\to 0$. Furthermore,
$f_\delta$ vanish
to infinite order as $t\to 0$.
Now, $\tilde{\Phi}(0)=0$ so it will follow
from the Nash-Moser inverse function theorem that $\tilde{\Phi}(u)=f_\delta$
has a smooth solution $u$ if $\delta$ is sufficiently small.
Then $x=u+x_0$ satisfies (2.11) for $0\leq t\leq \delta$.

In order to solve (2.11) or (2.13) we must show that the linearized operator is
invertible. Let us therefore first calculate the linearized equations.
Let $\delta$ be the Lagrangian variation, i.e. derivative w.r.t.
some parameter $r$ when $(t,y)$ are fixed. We have:
\proclaim{Lemma {2.}2} Let
$\overline{x}=\overline{x}(r,t,y)$ be a smooth function of
$(r,t,y)\in K=[-\varepsilon,\varepsilon]\times [0,T]\times\overline{\Omega} $,
$\varepsilon>0$,
such that $\overline{x}\big|_{r=0}=x$.
Then $\Phi(\overline{x})$ is a smooth function of $(r,t,y)\in K$, such that
$\pa \Phi(\overline{x})/\pa r\big|_{r=0}=\Phi^\prime( x) \delta x$,
where $\delta x=\pa \overline{x}/\pa r\big|_{r=0}$
and the linear map $L_0=\Phi^\prime(x)$ is given by
$$
\Phi^\prime(x)\delta x_i=D_t^2\delta x_i +(\pa_k \pa_i p)\delta x^k
+\pa_i \delta p_0
+\pa_i \big(\delta p_1-\delta x^k\pa_k p \big),\tag 2.14
$$
where $p$ satisfies (2.10) and $\delta p_i$, $i=0,1$, are given by solving
$$\align
\triangle & \big( \delta p_1- \delta x^k\pa_k p\big)=0,\qquad\qquad\qquad\qquad\quad
 \delta p_1\big|_{\pa\Omega}\!=0, \tag 2.15\\
\triangle  & \delta p_0
=-2(\pa_k V^i)\pa_i \big( \delta V^k-\delta x^l \pa_l  V^k\big)
,\qquad\quad \delta p_0\big|_{\pa\Omega}\!=0,  \tag 2.16
\endalign
$$
where $\delta v=D_t \delta x$.
Here, the normal operator
$$
A\delta x_i=-\pa_i \big(\pa_k p \,\, \delta x^k-\delta p_1\big)\tag 2.17
$$
restricted to divergence free vector fields is symmetric and positive,
in the inner product $\langle u,w\rangle=\int_{{\Cal D}_t} \delta^{ij} u_i w_j\, dx$,
if the physical condition (2.7) hold.
\endproclaim
\demo{Proof} That $\Phi(\overline{x})$ is a smooth function follows from that
the solution of (2.10) is a smooth function if $\overline{x}$ is, see section 16.
Let us now calculate $\Phi^\prime(x)$.
Since $[\delta, \pa /\pa y^a]=$ it follows that
$$
[\delta,\pa_i] =\Big(\delta \frac{\pa y^a}{\pa x^i}\Big)\frac{\pa}{\pa y^a}
-(\pa_i \delta x^l)\pa_l,\tag 2.18
$$
where we used the formula for the derivative of the inverse of a matrix
$\delta A^{-1}=-A^{-1}(\delta A)A^{-1}$.
It follows that
$[\delta\!-\!\delta x^l\pa_l,\pa_i]=0$ ($\delta\!-\!\delta
x^l\pa_l$ is the Eulerian variation ). Hence
$$\align
\delta \Phi_i-\delta x^k\pa_k \Phi_i&=D_t^2\delta x_i -(\pa_k D_t^2 x_i)\delta x^k
+\pa_i \big(\delta p-\delta x^k\pa_k p \big),\qquad\qquad \text{where}
\rightalignspace \tag 2.19\\
\triangle \big( \delta p- \delta x^k\pa_k p\big)
&=(\delta-\delta x^k\pa_k)\triangle p=
-2(\pa_k V^i)\pa_i \big( \delta V^k-\delta x^l \pa_l  V^k\big),
\qquad \delta p\big|_{\pa\Omega}\!=0.\rightalignspace \tag 2.20
\endalign
$$
The symmetry and positivity of $A$ were proven in
\cite{L1}, see also section 3 here.\qed\enddemo

In order to use the Nash-Moser iteration scheme to obtain a solution of (2.13) we
must  show that linearized operator is invertible and that the inverse
satisfies tame estimates:

\proclaim{Theorem {2.}3} Let
 $$
\||u\||_{a,k}=\sup_{0\leq t\leq T} \|u(t,\cdot)\|_{a,\infty}+...+\|D_t^k u(t,\cdot)\|_{a,\infty}
\tag 2.21
$$
where $\|u(t,\cdot)\|_{a,\infty}$ are the H\"older norms in $\overline{\Omega}$,
see (17.1).

Suppose that (2.7) and (2.8) hold initially, where
$p$ is given by (2.10),
and let $x_0\in C^\infty\big([0,T]\times\overline{\Omega}\big)$
satisfy (2.12).
Then there is a  $T_0=T(x_0)>0$, depending only on upper bounds for
$\||x_0\||_{4,2}$, $c_0^{-1}$ and $c_1$, such
that the following hold. If
$x\in C^\infty\big([0,T]\times\overline{\Omega}\big)$,
$p$ is defined by (2.10),
$$
T\leq T_0,\qquad \quad \||x-x_0\||_{4,2}\leq 1,
\qquad \text{and}\qquad (x-x_0)\big|_{t=0}=D_t\, (x-x_0)\big|_{t=0}=0, \tag 2.22
$$
then (2.7) and (2.8) hold for $0\leq t\leq T$ with $c_0$ replaced by $c_0/2$
and $c_1$ replaced by $2c_1$. Furthermore, linearized equations
$$
\Phi^\prime(x)\delta x=\delta \Phi,\qquad \text{in}\quad
[0,T]\times\overline{\Omega},\qquad\quad
\delta x\big|_{t=0}=D_t\,\delta x\big|_{t=0}=0.\tag 2.23
$$
where $\delta \Phi\in C^\infty\big([0,T]\times\overline{\Omega}\big)$
has a solution
$\delta x\in C^\infty\big([0,T]\times\overline{\Omega}\big)$.
The solution satisfies the estimates
$$
\||\delta x\||_{a,2}
\leq C_a\big( \||\delta \Phi\||_{a+r_0+2,0}
+\||\delta \Phi\||_{1,0}\, \||x-x_0\||_{a+r_0+6,2} \big),\qquad a\geq 0\tag 2.24
$$
where $C_a=C_a(x_0)$ is bounded when $a$ is bounded,
and in fact depends only on upper bounds for $ \||x_0\||_{a+r_0+6,2}$,
$c_0^{-1}$ and $c_1$. Here $r_0=[n/2]+1$, where $n$ is the number of space dimensions.

Furthermore $\Phi$ is twice differentiable and the second derivative
satisfies the estimates
$$
\||{\Phi}^{\prime\prime}(x)(\delta x,\epsilon x)\||_{a,0}
\leq C_a \Big( \||\delta x\||_{a+4,1}
\||\epsilon{x}\||_{2,1}+
 \||\delta x\||_{2,1}
\||\epsilon{x}\||_{a+4,1}
+\||x-x_0\||_{a+4,1} \||\delta x\||_{2,1}
\||\epsilon {x}\||_{2,1}\Big)\tag 2.25
$$
\endproclaim

The proof of Theorem {2.}1 follows from Theorem {2.}3 and Proposition {18.}1.
In Theorem {2.}3 we use norms that only has two time derivatives
and our Nash-Moser theorem, Proposition {18.}1, gives
a solution of (2.13)
$u\in C^2\big([0,T],C^\infty(\overline{\Omega})\big)$.
However, additional regularity in time follows from differentiating the
equations with respect to time. In fact,
if $x\in C^k\big([0,T],C^\infty(\overline{\Omega})\big)$
then $D_t^2 x=-\pa_i p\in C^{k-1}\big([0,T],C^\infty(\overline{\Omega})\big)$,
since (2.10) only depends on one time derivative of $x$, see the proof of Lemma {6.}7,
and it follows that $x\in  C^{k+1}\big([0,T],C^\infty(\overline{\Omega})\big)$.

Theorem {2.}3 follows from Lemma {14.}1, Proposition {15.1} and Proposition {16.}1.
The main point being existence for (2.23) and the tame  estimate (2.24)
given in Proposition {15.}1.
We will now discuss how to prove existence and estimates for the linearized equations.
The terms $(\pa_k \pa_i p)\delta x^k$
and $\pa_i \delta p_0$ in (2.14) are order zero in $\delta x$ and $D_t\delta x$.
The last term is a positive symmetric operator but only on divergence free vector fields
and in general it is an unbounded operator that looses regularity.
In general $\delta x$ is not going to be divergence free but we will derive evolution equations
for the divergence and the curl of $\delta x$, that gain regularity.
These evolution equations comes from that the divergence and the curl of the velocity $v$
are conserved expressed in the Lagrangian coordinates for a solution of Euler's equations,
$\Phi(x)=0$. In fact, since $[D_t,\pa_i]=-(\pa_i V^k)\pa_k$ it follows from (2.9) that
$$
D_t \div V=\div \Phi,\qquad\qquad\qquad
{\Cal L}_{D_t} \curl v=\curl \Phi\tag 2.26
$$
where $\curl v_{ij}\!=\!\pa_i v_j\!-\pa_j v_i$ and
${\Cal L}_{D_t}$ is the space time Lie derivative with respect to $D_t\!=\!(1,\!V)$:
$$
{\Cal L}_{D_t} \sigma_{ij}=D_t \,\sigma_{ij}
+(\pa_i V^l) \sigma_{lj}+(\pa_j V^l)\sigma_{il}\tag 2.27
$$
restricted to the space components. Expressing the two form $\sigma$ in the
Lagrangian frame this is just the time derivative:
$$
D_t \big( a^i_a a^j_b\sigma_{ij}\big)
=a^i_a a^j_b{\Cal L}_{D_t} \sigma_{ij},\qquad \text{where}
\quad a^i_a=\pa x^i/\pa y^a \tag 2.28
$$

We have the following evolution equations for the divergence and the curl of the linearized
operator
$$\align
\leftalignspace
\div\big(\Phi^\prime(x)\delta x\big)&=D_t^2\div \delta x+(\pa_i \delta x^k )\pa_k \Phi^i,
\rightalignspace\tag 2.29\\  \leftalignspace
\curl\,( \Phi^\prime(x)\delta x)&={\Cal L}_{D_t}\curl\,\big( D_t \,\delta x-\delta x^k\pa v_k
\big)+(\pa_i\delta x^k)\pa_j\Phi_k-(\pa_j\delta x^k)\pa_i\Phi_k
\rightalignspace\tag 2.30
\endalign
$$
In fact, since $[\delta,\pa_i]=-(\pa_i\delta x^k)\pa_k$ and
$[D_t ,\pa_i]=-(\pa_i V^k)\pa_k$ it follows that
 $\delta \div D_t x=D_t\div\delta x$ so by (2.26)
$D_t^2 \div \delta x=\delta \div \Phi$ and (2.29) follows.
To prove (2.30) we note that $[\delta,a_a^i a_b^j\pa_i]=
[\delta, a_b^j \pa_a] =(\delta a_b^j) \pa_a=(\pa_b \delta x^j)\pa_a=
a_a^i a_b^k (\pa_k \delta x^j)\pa_i  $ so
$$
\delta \big( a^i_a a^j_b \curl v_{ij}\big)
=a^i_a a^j_b\big(  \curl\delta  v_{ij}
+(\pa_j \delta x^k)\pa_i v_k-(\pa_i \delta x^k)\pa_k v_k\big)
= a^i_a a^j_b  \curl\, (\delta v-\delta x_k\pa V^k)_{ij}\tag 2.31
$$
where $ \curl\, (\delta v-\delta x_k\pa V^k)_{ij}\!
=\pa_i(\delta v_j\!-\delta x^k\pa_j v_k)
-\pa_j(\delta v_i\!-\delta x^k\pa_i v_k)$
and (2.30) follows since by (2.26)-(2.28)
$$
{\Cal L}_{D_t}
 \curl\, (\delta v-\delta x_k\pa V^k)=\curl\, (\delta \Phi-\delta x_k\pa \Phi^k).
\tag 2.32
$$

In \cite{L1} we proved existence and estimates
for the inverse of the linearized operator at a
solution of Euler's equations and within the divergence free class.
We only inverted  $\Phi^\prime(x)\delta x=\delta \Phi$
when $\delta \Phi$ was divergence free and $\Phi(x)=0$, in which case by (2.29)
$\delta x$ is also divergence free. In order to use the Nash-Moser iteration
scheme we will show that the linearized operator
is invertible away from a solution of Euler's equations and outside the divergence
free class.
This does present a  problem since the normal operator is only symmetric on divergence free
vector fields so for general vector fields we loose a derivative. In order to
recover this loss we will use that
one has better evolution equations for the divergence and for the curl that do not
loose regularity.
(2.29)-(2.30), says that we can get bounds for
the divergence and the curl of
$D_t\, \delta x$ if we have bounds for all first order derivatives of
 $ \, \delta x$. In fact (2.29)-(2.30) can be integrated even without
knowing a bound for first order derivatives of $D_t\, \delta x$.

We will now first modify the linearized operator so as to remove the term
$(\pa_i \delta x^k)\pa_k\Phi^i$ in (2.29) without making (2.30) worse.
(2.29) without this term will give
us an evolution equation that allows us to control the divergence.
This together with that the normal operator (2.17) is symmetric and
positive on divergence free vector fields will give us existence for the inverse
of the modified
linearized operator. The modified linearized operator is given by
$$
\align
L_1\delta x^i &=\Phi^\prime(x)\delta x^i
-\delta x^k\pa_k\Phi^i+\delta x^i\div \Phi\tag 2.33\\
&=D_t^2\delta x^i-(\pa_k D_t^2 x^i)\delta x^k
+\pa_i\big(\delta p_1-\delta x^k\pa_k p)+ \delta x^i\div \Phi+\pa_i\delta p_0
\endalign
$$
It follows from (2.29) that
$$
\div\,( L_1\delta x)=D_t^2\div\delta x+ \div\Phi\, \div\delta x
\tag 2.34
$$
The operator $L_1$ reduces to the linearized operator
$L_0=\Phi^\prime(x)$ when $\Phi(x)=0$ and the difference
$L_1-L_0$ is lower order. Furthermore, $L_1$ preserves the divergence free condition. We will first prove existence for the
inverse of the modified linearized operator and the existence of the inverse
of the linearized operator follows since the difference it is a lower order.
The main part of the manuscript is devoted to proving the following
existence and energy estimates:
\proclaim{Theorem {2.}4} Suppose that $x$ is smooth and that
the physical condition (2.7) and the
coordinate condition (2.8) hold for $0\leq t\leq T$. Then
$$
L_1 \delta x=\delta \Phi,\qquad 0\leq t\leq T,
\qquad \quad \delta x\big|_{t=0}=D_t\, \delta x\big|_{t=0}=0\tag 2.35
$$
has a smooth solution $\delta x$ if $\delta \Phi$ is smooth.

Furthermore, there are constants $K_4$ depending only on upper bounds
for $T$, $c_0^{-1}$, $c_1$, $r$ and $\||x\||_{4,2}$ such that the following estimates
hold. If $\div\delta \Phi=0$ then $\div\delta x=0$ and
$$
\|D_t \delta x\|_{r}+\|\delta x\|_r\leq K_4
  \int_0^t\big( \|\delta \Phi\|_r+\||x\||_{r+3,1}\|\delta \Phi\|_0\big) \, d\tau,
\qquad r\geq 0.
\tag 2.36
$$
If $\div\delta \Phi=0$, $\curl\delta\Phi=0$ and $\delta\Phi\big|_{t=0}=0$ then
$$\multline
\|D_t^2 \delta x\|_{r}+\|D_t \delta x\|_{r}+\|\delta x\|_r+c_0 \|\delta x\|_{r+1}\\
\leq K_4\int_0^t\big( \|D_t \delta \Phi\|_r+\|\delta \Phi\|_r
+ \||x\||_{r+3,2}( \|D_t \delta \Phi\|_0+ \|\delta \Phi\|_0)
 \big)\, d\tau,\qquad r\geq 0 \endmultline \tag 2.37
$$
In general
$$
\|D_t \delta x\|_{r-1}+\|\delta x\|_r\leq K_4\int_0^t\big( \|\delta \Phi\|_r
+ \||x\||_{r+3,2} \|\delta \Phi\|_1 \big)\, d\tau,\qquad r\geq 1\tag 2.38
$$
Here $\||x\||_{r,k}$ is as in Theorem {2.}3 and
$$
\|\delta x\|_r=\|\delta x(t,\cdot)\|_r
=\sum_{|\alpha|\leq r} \Big(\int_{{\Omega}}
|\pa_y^\alpha \delta x(t,y)|^2\, dy\Big)^{1/2}.\tag 2.39
$$
\endproclaim

The proof of the existence for (2.23) and the tame estimate (2.24) for the inverse
of the linearized operator in Theorem {2.}3 follows from Theorem {2.}4.
In fact, since the difference
$(L_1-\Phi^\prime(x))\delta x=O(\delta x)$ is lower order, the estimate (2.38)
will then allow us to get existence and the same estimate
 also for the inverse of the linearized operator (2.23), by iteration.
In (2.38) we only have estimates for one time derivative,
but we get estimates for an additional time derivative from also
using the equation. The $L^2$  estimates for (2.23) so obtained  then gives
the $L^\infty$ estimates (2.24) by also using Sobolev's lemma.

The proof of Theorem {2.}4 takes up most of the manuscript.
The proof (2.36) uses the symmetry and positivity of the normal operator (2.17) within
the divergence free class. This leads to energy estimates within the divergence free
class. The proof of (2.37) is obtained by first differentiating the equation with
respect to time and then using that a bound for two time derivatives
also gives a bound for the normal operator (2.17) using the equation.
The normal operator is not elliptic acting on general vector fields.
However, it is elliptic acting on divergence and curl free vector fields
and in general one can invert it and gain a space derivative
if one also has bounds for the curl and the divergence, see Lemma {5.}4.
Here we also need to use the improved estimate for the curl coming from (2.30).
To prove (2.38) we first subtract of a vector field picking up the divergence.
The equation for the divergence from (2.34):
$$
D_t^2 \div \delta x+\div\Phi\,\div\delta x=\div\delta\Phi\tag 2.40
$$
is just an ordinary differential equation that do not loose regularity
and in fact the estimates for (2.40)
gain an extra time derivative compared to the estimate (2.36).
Once we control the divergence we use the orthogonal projection onto divergence free
vector fields to obtain an equation for the divergence free part by projecting the
equation (2.35), see section 3. The equation so obtained is of the form (2.35)
with $\div \delta \Phi=0$ and $\delta \Phi$ depending also on the divergence
$\div \delta x$ that we just calculated.
The interaction term coming from the divergence part looses a space
derivative but it is in the form of a gradient
so we can recover this loss by using
the gain of a space derivative in (2.37).

In order to prove the energy estimates needed to prove Theorem {2.}4
one has to express the vector fields in the
Lagrangian frame, see (2.43). Theorem {2.}4, expressed in the
Lagrangian frame, follows from Theorem {10.}1,
Theorem {11.}1 and Theorem {12.}1.
Below, we will express the equation (2.35) in the Lagrangian frame
and in section 3 we outline the main ideas of how to
decompose the equation into a divergence free part and an equation for
the divergence using the orthogonal projection onto divergence free vector fields
and we show the basic energy estimate within the divergence free class.

%\subheading{The modified time derivative that preserves the divergence free condition}
As described above we now want to invert the modified linearized operator (2.35)
by decomposing it into an operator on the divergence free part and the ordinary differential equation
(2.40) for the divergence.
Hence we first want to be able to invert $L_1$ in the divergence free class.
The normal operator $A$, the third term on the second row in (2.33),
 maps divergence free vector fields onto divergence
free vector fields. We also
 want to modify the time derivative by adding a lower order term so it
 preserves the divergence free condition.
Let the Lie derivative and modified Lie derivative with respected to the time derivative
acting on vector fields be defined by
$$
{\Cal L}_{D_t}\delta x^i=D_t\delta x^i-(\pa_k V^i)\delta x^k,\qquad
\text{and}\qquad
\hat{\Cal L}_{D_t}\delta x^i={\Cal L}_{D_t}\delta x^i+\div V\, \delta x^i
\tag 2.41
$$
As before, ${\Cal L}_{D_t}$ is the space time Lie derivative restricted to the space components.
Then
$$
\div \hat{\Cal L}_{D_t} \delta x=\hat{D}_t \div\delta x,\qquad\quad
\text{where}\qquad \hat{D}_t=D_t+\div V\tag 2.42
$$
i.e. $\hat{D}_t f=D_t f+(\div V)\, f$.

This is  easier to see if we express the vector field in the Lagrangian frame.
Let
$$
 W^a=\frac{\pa y^a}{\pa x^i }\delta x^i\tag 2.43
$$
Then,
$$
D_t\, \delta x^i=D_t \big( W^b\pa x^i/\pa y^b\big)
= (D_t W^b)\pa x^i/\pa y^b+ W^b\pa V^i/\pa y^b
= (D_t W^b)\pa x^i/\pa y^b+ \delta x^k\pa_k V^i\tag 2.44
$$
and multiplying with the inverse $\pa y^a/\pa x^i$ gives
$$
D_t W^a= \frac{\pa y^a}{\pa x^i } {\Cal L}_{D_t} \delta x^i,
\qquad\text{and}\qquad
\hat{D}_t W^a= \frac{\pa y^a}{\pa x^i } \hat{\Cal L}_{D_t} \delta x^i.\tag 2.45
$$
With $\kappa=\det{(\pa x/\pa y)}$, we have
$$
\dot{W}^a=\hat{D}_t W^a=D_t W^a+(\div V)W^a=\kappa^{-1} D_t(\kappa W^a)\tag 2.46
$$
since  $D_t\, \kappa=\kappa \div V$, see \cite{L1}. Since the divergence is invariant
$$
\div \delta x=\div W=\kappa^{-1}\pa_a\big(\kappa W^a\big)\tag 2.47
$$
it therefore follows that
$$
\div \hat{D}_t W=\hat{D}_t \div W\tag 2.48
$$

The idea is now to replace the time derivatives $D_t$ in (2.33) by
$\hat{\Cal L}_{D_t}$ or equivalently express $L_1$ in the Lagrangian
frame and use the modified time derivatives $\hat{D}_t$.
Expressing the operator $L_1$ in the Lagrangian frame we get:
\proclaim{Lemma {2.}5} Let $\dot{W}=\hat{D}_t W$ and $\ddot{W}=\hat{D}_t^2 W$.
Then we can write (2.35) as $L_1 W=F$, where $W$ is given by (2.43), $F^a=\Phi^i \pa y^a/\pa x^i$ and
$$
 L_1 W^a=\ddot{W}^a+AW^a -B(W,\dot{W})^a,\qquad \qquad B(W,\dot{W})^a=
 B_0 W^a +B_1\dot{W}^a .
 \tag 2.49
$$
Here
$$\align
g_{ab}A {W}^b&=-\pa_a\big( (\pa_c p) W^c-q_1\big),\qquad \qquad \qquad  \qquad
\qquad \div A W=0\tag 2.50\\
g_{ab}B_0{W}^b&= \dot{\sigma}\big(D_t {g}_{ac}-\omega_{ac}
-\dot{\sigma}g_{ac}\big)W^c-\pa_a q_3 ,
\qquad\qquad  \div B_0{W}=-\dot{\sigma}^2 \div {W}\tag 2.51\\
g_{ab}B_1\dot{W}^b&= -\big(D_t{{g}}_{ac}-\omega_{ac} -2\dot{\sigma} g_{ac}\big)
\dot{W}^c-\pa_a q_2 ,\qquad \qquad \div B_1\dot{W}=2\dot{\sigma}\div \dot{W},
\tag 2.52
\endalign
$$
where $q_i$, for $i=1,2,3$ are given by solving the Dirichlet problem $q_i\big|_{\pa\Omega}=0$
where $\triangle q_i$ are given by the equations for the divergences above,
$\sigma=\ln\kappa$, $ \dot{\sigma}=D_t\,\sigma=\div V$,
$\ddot{\sigma}=D_t^2\sigma$ and
$$
D_t\, g_{ab}=\frac{\pa x^i}{\pa y^a}\frac{\pa x^j}{\pa y^b}\big(\pa_i v_j\!+\pa_j v_i\big),
\qquad \omega_{ab}
=\frac{\pa x^i}{\pa y^a}\frac{\pa x^j}{\pa y^b}\big(\pa_i v_j\!-\pa_j v_i\big).
\tag 2.53
$$

We have
$$
\div\, (L_1 W)=D_t^2\div W+\ddot{\sigma}\div W\tag 2.54
$$
Let $\underline{L}_1 W_a=g_{ab} L_1 W^b$,  $\dot{w}_a=g_{ab}\dot{W}^b$ and
$\tilde{w}_a=\dot{w}_a-(\omega_{ab}+\dot{\sigma} g_{ab})W^b$. Then
$$\align
\curl \, (\underline{L_1} W)&=D_t \curl \tilde{w} +\curl \underline{B}_4 W
\tag 2.55\\
\curl \, (\underline{L_1} W)&=D_t \curl \dot{w} +\curl \underline{B}_5 \dot{W}
+\curl \underline{B}_6 W\tag 2.56
\endalign
$$
where $\underline{B}_4 W_a=(D_t\, {\omega}_{ab}+\ddot{\sigma} g_{ab})W^b$,
 $\underline{B}_5 \dot{W}_a=-({\omega}_{ab}+\dot{\sigma} g_{ab})\dot{W}^b$
and  $\underline{B}_6 W_a=-\dot{\sigma}(D_t \, g_{ab}-{\omega}_{ab}-\dot{\sigma} g_{ab})W^b$.

Furthermore $L_0=\Phi^\prime(x)$ expressed in the Lagrangian frame is given by
$$
L_0 W^a=L_1 W^a -B_3 W^a ,\qquad\text{where}\qquad
B_3 W^a=-W^c\na_c \Phi^a+W^a\div \Phi\tag 2.57
$$
where $\nabla_c$ is covariant differentiation with respect to the metric $g_{ab}$
and $\Phi^a=\Phi^i\pa y^a/\pa x^i$, i.e.
$\na_c \Phi^a=(\pa x^i/\pa y^c)(\pa y^a/\pa x^j)\pa_i \Phi^j$.
\endproclaim
\demo{Proof}
Differentiating (2.44) once more gives
$$
D_t^2 \delta x^i-(\pa_k D_t V^i)\delta x^k=(D_t^2 W^b)\pa x^i/\pa y^b
+2(D_t W^b)\pa V^i/\pa y^b\tag 2.58
$$
It follows that
$$\aligned
\frac{\pa x^i}{\pa y^a }\big(D_t^2 \delta x^i-(\pa_k D_t V^i)\delta x^k\big)
&=\frac{\pa x^i}{\pa y^a }\frac{\pa x^i}{\pa y^b }D_t^2 W^b
+2(D_t W^b)\frac{\pa x^i}{\pa y^b }\frac{\pa x^j}{\pa y^a }
\pa_i v_j\\
&=g_{ab} D_t^2 W^b+ (D_t\, g_{ab}-\omega_{ab}) D_t W^b
\endaligned\tag 2.59
$$
It follows from (2.33) that
$$
\align
g_{ab} L_1 W^b&=g_{ab} D_t^2 W^b -\pa_a\big( (\pa_c p) W^c-q\big)
+(D_t{g}_{ac}-\omega_{ac}) D_t W^c+\ddot{\sigma}g_{ab} W^b \tag 2.60\\
&=  D_t\big(g_{ab}D_t  W^b-\omega_{ab}W^b\big)
 -\pa_a\big( (\pa_c p) W^c-q\big)+D_t{\omega}_{ab}W^b
 +\ddot{\sigma}g_{ab} W^b
\endalign
$$
where $q=\delta p$ is chosen so that the divergence is
equal to $\div L_1 W=D_t^2 \div W +\div W\, D_t \div V$ in order for it to
be consistent with (2.34).
We have
$\hat{D}_t^2=(D_t+\div V)(D_t+\div V)=
D_t^2+2\dot{\sigma} D_t +\dot{\sigma}^2+\ddot{\sigma}=
D_t^2+2\dot{\sigma} \hat{D}_t +\ddot{\sigma}-\dot{\sigma}^2$
so
$$
D_t^2=\hat{D}_t^2-2\dot{\sigma}\hat{D}_t+\dot{\sigma}^2-\ddot{\sigma},
\qquad D_t=\hat{D}_t-\dot{\sigma}\tag 2.61
$$
Hence, with $\dot{W}=\hat{D}_t W$ and $\ddot{W}=\hat{D}_t^2 W$,
we can write the equation (2.60) as
$$
 L_1 W^a=\ddot{W}^a-g^{ab}\pa_b\big( (\pa_c p) W^c-q_1\big)
 -B^a(W,\dot{W})
 \tag 2.62
$$
where $q_1$ is chosen so the divergence of the second term on the right
vanishes and
$$
g_{ab}B^b(W,\dot{W})=
-\big(D_t{{g}}_{ac}-\omega_{ac} -2\dot{\sigma} g_{ac}\big)
\dot{W}^c
+\big(\dot{\sigma}\big(D_t{g}_{ac}-\omega_{ac}-\dot{\sigma}g_{ac}\big)W^c
-\pa_a q_0\tag 2.63
$$
Here $q_0$ is chosen as follows so that $\div L_1 W=\hat{D}_t^2 \div W-\div B=
D_t^2\div W+\div W\, \ddot{\sigma} $. But $\hat{D}_t^2 \div W=
D_t^2\div W+2\dot{\sigma} \hat{D}_t \div W+(\ddot{\sigma}-\dot{\sigma}^2)\div W$
so we must have $\div B=2\dot{\sigma} \hat{D}_t\div W -\dot{\sigma}^2 \div W$.
Hence $q_0$ is chosen so this is fulfilled and (2.49) follows by writing $q_0=q_2+q_3$.
(2.54) follows from (2.34) or (2.49)
It follows from (2.49) that we write $L_1$ in the two alternative forms:
$$\align \leftalignspace
g_{ab} L_1 W^b&= D_t\big(g_{ab}\dot{W}^b-(\omega_{ab}+\dot{\sigma} g_{ab})W^b\big)
 -\pa_a\big( (\pa_c p) W^c-q_1\big)+(D_t {\omega}_{ab}
 +\ddot{\sigma}g_{ab}) W^b+\pa_a q_0 \rightalignspace \tag 2.64\\ \leftalignspace
g_{ab} L_1 W^b&= D_t\big(g_{ab}\dot{W}^b\big)
 -\pa_a\big( (\pa_c p) W^c-q_1\big)-(\omega_{ab}+\dot{\sigma}g_{ab})
 (\dot{W}^b-\dot{\sigma} W^b)
 -\dot{\sigma}D_t{g}_{ab}W^b
 +\pa_a q_0\rightalignspace\tag 2.65
\endalign
$$
(2.55) and (2.56) follows from these.
 Finally, we also want to express $L_0=\Phi^\prime(x)$ is these coordinates.
In order to to this we must transform the term $\delta x^k\pa_k \Phi^i$ in (2.33)
to the Lagrangian frame.
If $\Phi^a=\Phi^i \pa y^a/\pa x^i$, then
 $(\delta x^k\pa_k \Phi^i) \pa y^a/\pa x^i=W^c\na_c \Phi^a$, where $\na_c$ is covariant
differentiation, see e.g. \cite{CL}, and (2.57) follows.
\qed\enddemo

\comment
$$\align
\leftalignspace
{\Cal L}_{D_t}^2 \delta x^i
&=D_t^2\delta x^i -\delta x^k\pa_k D_t^2 x^i-2(\pa_k V^i)
{\Cal L}_{D_t}\delta x^k \rightalignspace
\tag 2.40\\
\hat{\Cal L}_{D_t}^2 \delta x^i&={\Cal L}_{D_t}^2 \delta x^i
+2\div V {\Cal L}_{D_t}\delta x^i+\big(D_t\div V+(\div V)^2\big) \delta x^i
\tag 2.42
\endalign
$$
Since $D_t\div V=\div\Phi$ we can hence write
$$
L_1 \delta x^i=\hat{\Cal L}_{D_t}^2\delta x^i+A \delta x^i
+\big( 2(\pa_k V^i){\Cal L}_{D_t}\delta x^k-2\div V{\Cal L}_{D_t}\delta x^i
-(\div V)^2\delta x^i
+\pa_i \delta p_0\big) \tag 2.42
$$
where the divergence of of the first two terms vanish and,  in view of (2.34),
$$\multline
\pa_i\Big( 2(\pa_k V^i){\Cal L}_{D_t}\delta x^k-2\div V{\Cal L}_{D_t}\delta x^i
-(\div V)^2\delta x^i
+\pa_i \delta p_0\Big) \\
=(\div V)^2\div \delta x-2\div V\, \div \hat{\Cal L}_{D_t}\delta x
=-(\div V)^2\div \delta x-2\div V\, D_t\div\delta x
\endmultline \tag 2.43
$$
\endcomment

\head{3. The projection onto divergence free vector fields and the normal operator.
}\endhead
%\subheading{The orthogonal projection onto divergence free vector fields}
Let us now also define the projection $P$ onto divergence free
vector fields by
$$
PU^a=U^a-g^{ab}\pa_b p_{U},\qquad \triangle p_U=\div U,
\qquad p_U\big|_{\pa\Omega}=0\tag 3.1
$$
(Here $\triangle q=\kappa^{-1}\pa_a\big( \kappa g^{ab}\pa_b q\big)$. )
$P$ is the orthogonal projection in the inner product
$$
\langle U,W\rangle =\int_\Omega g_{ab} U^a W^b\kappa dy\tag 3.2
$$
and its operator norm is one:
$$
\|PW\|\leq \|W\|,\qquad\quad\text{where}\qquad \|W\|=\langle W,W\rangle^{1/2}.\tag 3.3
$$

%\subheading{The normal operator}
For a function $f$ that vanishes on the boundary define
$A_f W^a=g^{ab}\underline{A}_f W_b $, where
$$
\underline{A}_f W_a=-\pa_a\big( (\pa_c f) W^c-q\big),\qquad
\triangle \big( (\pa_c f) W^c-q\big)=0,
\qquad q\big|_{\pa\Omega}=0,\tag 3.4
$$
i.e. $A_f W$ is the projection of $-g^{ab}\pa_b\big( (\pa_c f) W^c\big)$.
This is defined for general vector fields but it is only symmetric in the
divergence free class. We have
$$
\langle U,A_f W\rangle =\int_{\pa\Omega} n_a \, U^a (-\pa_c f) W^c\, dS,
\qquad\quad\text{if}\qquad \div U=\div W=0,
\tag 3.5
$$
where $n$ is the unit conormal. If $f\big|_{\pa\Omega}=0$ then
$-\pa_c f\big|_{\pa\Omega} =(-\na_N f) n_c$. It follows that
$A_f$ is a symmetric operator on divergence free vector fields, and
in particular, the normal operator in (2.50)
$$
A=A_p\tag 3.6
$$
is positive since we assumed that $-\na_N p\geq c>0$
on the boundary. We have
$$
|\langle U,A_f W\rangle|\leq
\| \na_N f/\na_N p\|_{L^\infty(\pa\Omega)}\langle U, A U\rangle^{1/2}
\langle W, A W\rangle^{1/2},
\qquad\quad\text{if}\qquad \div U=\div W=0. \tag 3.7
$$
Since the projection has norm one it follows from (3.4) that
$$
\| A_f W\|\leq \|\pa^2 f\|_{L^\infty(\Omega)} \|W\|+
\|\pa f\|_{L^\infty(\Omega)} \|\pa W\|.\tag 3.8
$$
Note also that $A_f$ acting on divergence free vector fields by (3.5)
depends only on $\na_N f\big|_{\pa\Omega}$, i.e. $A_{\tilde{f}}=A_f$
if $\na_N \tilde{f}\big|_{\pa\Omega}=\na_N f\big|_{\pa\Omega}$.
We can therefore replace $f$ by the Taylor expansion of order one in the distance
to the boundary in polar coordinates multiplied by a smooth function that is one
close to the boundary and vanishes close to the origin. It follows that
$$
\|A_f W\|\leq C\sum_{S\in{\Cal S}}\|\na_N Sf\|_{L^\infty(\pa\Omega)}\|W\|
+ C\|\na_N f\|_{L^\infty(\pa\Omega)}(\|\pa W\|+\|W\|),
\qquad\text{if}\quad \div W\!=0,\tag 3.9
$$
where ${\Cal S}$ is a set of vector fields that span the tangent space of
$\pa\Omega$, see section 4.

%\subheading{The smoothed out normal operator}
In order to prove existence for the linearized equations
we in \cite{L1} replaced the normal operator $A$ by a smoothed out bounded operator that still has the
same positive properties as $A$ and commutators with Lie derivatives, and which also
has vanishing divergence and curl away from the boundary.
This makes possible to pass to the limit and obtain existence for the linearized equations.
The smoothed out normal operator is defined as follows.
Let $\rho=\rho(d)$ be a smooth out version of the distance
function to the boundary $d(y)=\dist(y,\pa\Omega)=1-|y|$ in the
standard Euclidean metric $\delta_{ij} d y^i d y^j$  in the $y$ coordinates,
$\rho^\prime\geq 0$, $\rho(d)=d$, when $d\leq 1/4$ and
$\rho(d)=1/2$ when $d\geq 3/4$. Then we can alternatively express $A_f$ as
$$
\underline{A}_f W_a=-\pa_a\big( (f/\rho)(\pa_c \rho) W^c-q\big),\qquad
\triangle \big( (f/\rho)(\pa_c \rho) W^c-q\big)=0,
\qquad q\big|_{\pa\Omega}=0\tag 3.10
$$
Let $\chi(\rho)$ be a smooth function such that $\chi^\prime\geq 0$,
$\chi(\rho)=0$ when $\rho\leq 1/4$, $\chi(\rho)=1$ when $\rho\geq 3/4$.
$A_f$ is unbounded so we now define an approximation that
is a bounded operator: $A_f^\varepsilon W^a=g^{ab} \underline{A}_f^\varepsilon W_b$,
where
$$
\underline{A}_f^\varepsilon W_a=-\chi_\varepsilon\pa_a\big(
(f/\rho)(\pa_c \rho) W^c\big) +\pa_a q ,\quad
\triangle q
=\kappa^{-1}\pa_a \big( g^{ab} \kappa\chi_\varepsilon
\pa_b\big( (f/\rho)(\pa_c \rho) W^c\big)\big),
\quad q\big|_{\pa\Omega}\!=0\tag 3.11
$$
where $\chi_\varepsilon(\rho)=\chi(\rho/\varepsilon)$.
We have
$$
\langle U,A_f^\varepsilon W\rangle
=\int_{\Omega} (f/\rho) \chi_\varepsilon^\prime (\pa_a\rho) \, U^a (\pa_c \rho) W^c\,
\kappa dy,\qquad\quad\text{if}\qquad \div U=\div W=0,
\tag 3.12
$$
from which it follows that $A_f^\varepsilon$ is also symmetric.
And in particular $A^\varepsilon=A^\varepsilon_p$ is positive since we
assumed that $p\geq 0$, at least close to the boundary.
We have
$$
|\langle U,A^\varepsilon_f W\rangle|\leq
\| f/p\|_{L^\infty(\Omega\setminus\Omega_{\varepsilon/4})}\langle U, A^\varepsilon U\rangle^{1/2}
\langle W, A^\varepsilon W\rangle^{1/2},
\qquad\quad\text{if}\qquad \div U=\div W=0, \tag 3.13
$$
where $\Omega_\varepsilon=\{y\in\Omega; d(y,\pa\Omega)<\varepsilon\}$.
It also follows from (3.12) that another expression for $\underline{A}_f^\varepsilon$ is
$$
\underline{A}_f^\varepsilon W_a=
 (f/\rho) \chi_\varepsilon^\prime (\pa_a\rho) \,  (\pa_c \rho) W^c\!-\pa_a q,
\qquad \triangle q=\kappa^{-1}\pa_a\big( \kappa g^{ab}
(f/\rho) \chi_\varepsilon^\prime (\pa_b\rho) \,  (\pa_c \rho) W^c\big),
\quad q\big|_{\pa\Omega}\!=0\tag 3.14
$$
acting on divergence free vector fields.
Furthermore, by (3.12)
$$
\|D_t^k A^\varepsilon W\|_r\leq C_{\varepsilon} \sum_{j=0}^k
\|D_t^j W\|_r,\qquad\quad\text{where}\qquad \|W\|_r=\sum_{|\alpha|\leq r}\|\pa_y^\alpha W(t,\cdot)\|_{L^2(\Omega)}
\tag 3.15
$$

%\subheading{The multiplication operators}
Let us also define the projected multiplication operators $M_\beta$ with a two form
$\beta$ by
$$
\underline{M}_\beta W_a=\underline{P}(\beta_{ab} W^b)\tag 3.16
$$
Since the projection has norm one it follows that
$$
\|M_\beta W\|\leq \|\beta\|_{\infty}\|W\|\tag 3.17
$$
Furthermore we define the
operator taking vector fields to one forms
$$
\underline{G} W_a=\underline{M}_{g} W_a=P(g_{ab} W^b)\tag 3.18
$$
Then $G$ acting on divergence free vector fields is just the identity $I$.

Let $L_1$ be the modified linearized operator in (2.49) and
let $\dot{W}=\hat{D}_t W=D_t W+(\div V) W =\kappa^{-1} D_t(\kappa W)$,
$\ddot{W}=\hat{D}_t^2 W$.
We want to prove existence
of a solution $W$ to
 $$
 L_1 W=\ddot{W}+AW-B_0 W-B_1\dot{W}=F,\qquad \quad W\big|_{t=0}=\dot{W}\big|_{t=0}=0
\tag 3.19
 $$
for general vector fields $F$ that are not
necessarily divergence free. To do this we first subtract of a vector field
$W_1$ that picks up the divergence and then solve (3.19) in the divergence free class.
Let us decompose a vector field into a divergence free part and a gradient
using the orthogonal projection:
$$
W=W_0+W_1,\qquad W_0=PW, \qquad W_1^a=g^{ab}\pa_b q_1,\qquad
q_1\big|_{\pa\Omega}=0.\tag 3.20
$$
Then if $\dot{g}_{ab}=\check{D}_t g_{ab}$, where
$\check{D}_t=D_t-\dot{\sigma} $, we have $ \pa_a D_t \,q_1=D_t(
g_{ab} W_1^b)=\dot{g}_{ab} W_1^b+g_{ab}\dot{W}_1^b$ and $\pa_a D_t^2
q_1= \ddot{g}_{ab} W_1^b +2\dot{g}_{ab}
\dot{W}_1^b+g_{ab}\ddot{W}_1^b$, where $\ddot{g}_{ab}=\check{D}_t^2 g_{ab}$. Hence
$$
\ddot{W}_1^a =g^{ab}\pa_b D_t^2 q_1-2 g^{ab}\dot{g}_{bc}
\dot{W}_1^c-g^{ab}\ddot{g}_{bc} W_1^c\tag 3.21
$$
Since $D_t^2 q_1\big|_{\pa\Omega}=0$ and
the projection of a gradient of a function that vanishes on
the boundary vanishes
$$
P\ddot{W}_1^a=B_2(W_1,\dot{W}_1)^a,\qquad\text{where}
\qquad B_2(W_1,\dot{W}_1)^a=-P\big(2 g^{ab}\dot{g}_{bc}
\dot{W}_1^c+g^{ab}\ddot{g}_{bc} W_1^c\big)\tag 3.22
$$
Since $\div{W}_0=0$ it follows that $\div\dot{W}_0=\div\ddot{W}_0=0$
and hence by Lemma {2.}5
$$
\align
P L_{1} W_0&=L_1 W_0=\ddot{W}_0+A W_0 -B_{1} \dot{W}_0-B_{0} W_0\tag 3.23\\
PL_{1} W_1&=AW_1-B_{11} \dot{W}_1-B_{01} W_1 \tag 3.24
\endalign
$$
where
$$
 B_{11}\dot{W}^a\!=PB_{1} \dot{W}^a\! +2 P\big(
g^{ab}\dot{g}_{bc} \dot{W}^c)\,\qquad
 B_{01} W^a\!=PB_{0} W^a\!
+P\big(g^{ab}\ddot{g}_{bc} W^c\big).\tag 3.25
$$
Hence projection of (3.19) gives
$$
L_{1} W_{0}=-PL_{1}W_1+PF=-AW_1+B_{11} \dot{W}_1+B_{01} W_1+PF,\tag 3.26
$$
Here, by (2.54)
$$
W_1^{a}=g^{ab}\pa_b q_1,\qquad \triangle\,
q_1=\varphi,\qquad q_1\big|_{\pa\Omega}=0,\tag
3.27
$$
where
$$
D_t^2\varphi +\ddot{\sigma}\, \varphi=\div F.\tag 3.28
$$
By (3.23)-(3.24) we also have
$$\align
(I-P) L_1 W_0&=0\tag 3.29\\
(I-P) L_1 W_1&=\ddot{W}_1-B_2(W_1,\dot{W}_1)-(I-P)B_0 W+(I-P)B_1\dot{W}_1\tag 3.30
\endalign
$$
Summing up, we have proven:
\proclaim{Lemma {3.}1} Suppose that $W$ satisfies $L_1 W=F$. Let
$W_0=PW$, $W_1=(I-P)W$, $F_0=PF$ and $F_1=(I-P)F$. Then
$$\align
L_1 W_0&=F_0-A W_1+B_{11} \dot{W}_1+B_{01} W_1\tag 3.31\\
\ddot{W}_{11}&=F_1+B_2(W_1,\dot{W_1})+(I-P) B_0 W_1+(I-P) B_1\dot{W}_1\tag 3.32
\endalign
$$
where $B_{01}$ and $B_{11}$ are given by (3.25), $B_2$ is given by (3.22)
and $B_0$, $B_1$ is as in (2.51), (2.52).
Furthermore
$$
D_t^2 \div W_1+\ddot{\sigma}\div W_1=\div F\tag 3.33
$$
\endproclaim
We now find a solution of (3.19) by first solving the ordinary differential
equation (3.28) and then solving the Dirichlet problem for $q_1$ and defining
$W_1$ by (3.27).
Finally we solve (3.26) for $W_0$ within the divergence free class.
This gives existence of solutions for (3.19) for general vector fields $F$ once
we can solve it for divergence free vector fields.
However, we also need estimates for (3.19) that do not loose regularity going from
$F$ to $W$ in order to show existence also for the linearized
equations (2.57):
$$
L_0 W=L_1 W-B_3 W=F, \qquad\quad W\big|_{t=0}=\dot{W}\big|_{t=0}=0,\tag 3.34
$$
by iteration. It seems like there is a loss
of regularity in the term $-AW_1$ in (3.26).
However, $\curl AW_1=0$ and there is an improved estimate for (3.19) when
$\div F=0$ and $\curl F=0$, obtained by differentiating with respect to time
and using that an estimate for two time derivatives also gives an estimate
for the operator $A$ through the equation (3.19). We can estimate
any first order derivative of a vector field in terms
of the curl, the divergence and the normal operator $A$ and
there is an  identity for the curl.

%\subheading{The basic energy estimate}
Let us now also derive the basic energy estimate which will be used to prove
existence and estimates for (3.19) within the divergence free class:
$$
\ddot{W}+AW=H,\qquad\quad W\big|_{t=0}=\dot{W}\big|_{t=0}=0,\qquad
\quad   \div H=0\tag 3.35
$$
where $A$ is the normal operator or the smoothed version.
For any symmetric operator $B$ we have
$$
\frac{d}{dt}\langle W, BW\rangle =\frac{d}{dt}
\int_{\Omega} \kappa W^a \underline{B}W_{\!a }\, dy=2\langle \dot{W}, BW\rangle
+\langle W,\dot{B}W\rangle\tag 3.36
$$
where $\dot{W}=\kappa^{-1}D_t(\kappa W)$ and
$\dot{B}$ is the time derivative of the operator $B$
considered as an operator from the divergence free vector fields
to the one forms corresponding to divergence free vector fields:
$$
\dot{B} W^a=P\big(g^{ab}( D_t \underline{B}W_b -\underline{B} \dot{W}_b )\big),
\qquad \quad \underline{B}W_b=g_{bc} B W^c,
\tag 3.37
$$
see section 4. The projection
comes up here since we take the inner product with a
divergence free vector field in (3.37).
 Let the lowest order
energy $E_0=E(W)$ be defined by
$$
E(W)=\langle \dot{W},\dot{W}\rangle +\langle W,(A+I)W\rangle\tag 3.38
$$
Since $\langle W,W\rangle =\langle W,GW\rangle$, where $G$ is the projection onto
divergence free vector fields given by (3.18),
it follows that
$$
\dot{E}_0=2\langle \dot{W},\ddot{W}+(A+I)W\rangle
+\langle \dot{W},\dot{G}\dot{W}\rangle+\langle W,(\dot{A}+\dot{G})W\rangle
\tag 3.39
$$
 In particular it follows from (3.4) or (3.10)
respectively (3.16) and (3.18) that
$$
\dot{A}_f=A_{\dot{f}},\qquad
\dot{G}=M_{\dot{g}},\qquad\text{where}\qquad
\dot{f}=\kappa D_t(\kappa^{-1} f)\quad\text{and}\quad
\dot{g}=\kappa D_t(\kappa^{-1} g).
\tag 3.40
$$
In fact the time derivate of an operator,
as defined by (3.37), commutes with the projection since
$D_{t\,} \pa_a q=\pa_a D_{t\,} q$, where
$D_{t\,} q\big|_{\pa\Omega}\!=0$ if $q\big|_{\pa\Omega}\!=0$,
and the projection of the gradient of function that vanishes on the
boundary vanishes.
It therefore follows from (3.7) or (3.12) and (3.17) that
 $$
 |\langle
W,\dot{A}W\rangle| \leq \|\dot{p}/p\|_\infty \langle W,AW\rangle
,\qquad |\langle W,\dot{G}W\rangle |\leq \|\dot{g}\|_\infty
\langle W,W\rangle\tag 3.41
$$
The last two terms in (3.38) are hence bounded by a constant times the
energy so it follows that
$$
|\dot{E}_0|\leq \sqrt{E_0}\big( 2\|H\| +c\sqrt{E_0}\big), \qquad
c=\|\dot{p}/p\|_\infty+\|\dot{g}\|_\infty+2 \tag 3.42
$$
from which a bound for the lowest order energy follows.

Similarly, we get higher order energy estimates for vector fields that are
tangential at the boundary, see section 10. Once we have these estimates
we use that any derivative of a vector field
can be bounded by tangential derivatives
and derivatives of the divergence and the curl, see section 5.
The divergence vanishes and we can get estimates for the
curl as follows. Let $w_a=g_{ab} W^b$, $\dot{w}_a=g_{ab} \dot{W}^b$ and
$\ddot{w}_a=g_{ab}\ddot{W}^b$.
Then $D_t w_a= \dot{g}_{ab}W^b+\dot{w}_a$ and $ D_t
\dot{w}_a=\dot{g}_{ab}\dot{W}^b+\ddot{w}_a$
where $\dot{g}_{ab}=\check{D}_t g_{ab}=\kappa D_t(\kappa g_{ab})$.
Since
$$
\ddot{w}+\underline{A} W=\underline{H},\qquad H=B_0
W+B_1\dot{W}+F\tag 3.43
$$
where $\curl \underline{A}W=0$ it follows that
$$
|D_t\curl w|+
|D_t \curl \dot{w}|\leq C\big(|\pa W|+|W| +|\pa \dot{W}|
+|\dot{W}|+|\curl \underline{F}|\big)\tag 3.44
$$
Note that the estimate for the curl is actually very strong. The
higher order operator $A$ vanishes so there is no loss of
regularity anymore and furthermore the estimate is point wise.
This crude estimate suffices for the most part. However, there is
an additional cancellation, whereas one would not need to assume
estimate for $|\pa\dot{W}|$ in the right hand side of (3.41). The
improved estimate is for $\dot{w}_a$ replaced by
$\tilde{w}_a=\dot{w}_a-\omega_{ab}W^b$, where $\omega_{ab}=\pa_a
v_b-\pa_b v_a$. It follows from Lemma {2.}5 that
$$
|D_t\curl w|+|D_t \curl \tilde{w}|\leq C\big(|\curl\tilde{w}|+
|\pa W|+|W|+|\curl \underline{F}|\big),\qquad
|\curl (\tilde{w}-\dot{w})|\leq C\big(|W|+|\pa W|\big)\tag 3.45
$$

\head 4. The tangential vector fields, Lie derivatives
and commutators.\endhead
%\subheading{The tangential vector fields}
Following \cite{L1}, we now construct the tangential vector fields, that are
time independent expressed in the Lagrangian coordinates, i.e. that commute
with $D_t$. This means that in the Lagrangian coordinates they are of the form
$S^a(y)\pa/\pa y^a$. Furthermore, they will satisfy,
$$
\pa_a S^a=0,\tag 4.1
$$
Since $\Omega$ is the unit ball in $\bold{R}^n$
the vector fields can be explicitly given.
The vector fields
$$
y^a\pa/\pa y^b-y^b\pa/\pa y^a\tag 4.2
$$
corresponding to rotations, span the tangent space of
the boundary and are divergence free in the interior.
Furthermore they span the tangent space of the level sets of the
distance function from the boundary in the Lagrangian coordinates
$$
d(y)=\dist{(y,\pa \Omega)}=1-|y|\tag 4.3
$$
away from the origin $y\neq 0$.
We will denote this set of vector fields by ${\Cal S}_0$
We also construct a set of divergence free
vector fields that span the full tangent space
at distance $d(y)\geq d_0$ and that are compactly supported
in the interior at a fixed distance $d_0/2$ from the boundary.
The basic one is
$$
h(y^3,...,y^n)\Big(f(y^1)g^\prime(y^2)\pa/\pa y^1-
f^\prime(y^1)g(y^2)\pa/\pa y^2\Big),\tag 4.4
$$
which satisfies (4.1).
Furthermore we can choose $f,g,h$ such that
it is equal to $\pa/\pa y^1 $ when $|y^i|\leq 1/4$, for $i=1,...,n$
and so that it is $0$ when $|y^i|\geq 1/2$ for some $i$.
In fact let $f$ and $g$ be smooth functions such that
$f(s)=1$ when $|s|\leq 1/4$ and $f(s)=0$ when $|s|\geq 1/2$
and $g^\prime(s)=1$ when $|s|\leq 1/4$ and $g(s)=0$ when $|s|\geq 1/2$.
Finally let $h(y^3,...,y^n)=f(y^3)\cdot\cdot\cdot f(y^n)$.
By scaling, translation and rotation
of these vector fields we can obviously construct a finite set
of vector fields that span the tangent space when $d\geq d_0$ and are
compactly supported in the set where $d\geq d_0/2$.
We will denote this set of vector fields by ${\Cal S}_1$.
Let ${\Cal S}={\Cal S}_0\cup {\Cal S}_1$ denote
the family of tangential space vector fields  and let
${\Cal T}={\Cal S}\cup \{D_t\}$ denote the family of space time
tangential vector fields.

Let the radial vector field be
$$
R=y^a\pa/\pa y^a. \tag 4.5
$$
Now,
$$
\pa_a R^a=n\tag 4.6
$$
 is not $0$ but for our purposes it suffices that
it is constant. Let ${\Cal R}={\Cal S}\cup\{R\}$.
Note that ${\Cal R}$ span the full tangent space of the space everywhere.
Let ${\Cal U}={\Cal S}\cup \{R\}\cup\{D_t\}$ denote the family of all vector fields.
Note also that the radial vector field commutes with the rotations;
$$
[R,S]=0,\qquad S\in {\Cal S}_0\tag 4.7
$$
Furthermore, the commutators of two vector fields in ${\Cal S}_0$
is just $\pm$ another vector field in ${\Cal S}_0$.
Therefore, for $i=0,1$, let ${\Cal R}_i={\Cal S}_i\cup\{R\}$, ${\Cal T}_i={\Cal S}_i\cup\{D_t\}$
and ${\Cal U}_i={\Cal S}_i\cup\{R\}\cup\{ D_t\}$.

Let us now introduce the Lie derivative of the vector field $W$
with respect to the vector field $T$;
$$
{\Cal L}_T W^a=TW^a -(\pa_c T^a) W^c\tag 4.8
$$
We will only deal with Lie derivatives with respect to the vector fields $T$
constructed above. For those vector fields $T$ we have
$$
 [D_t,T],\qquad\text{and}\qquad [ D_t , {\Cal L}_{T} ]=0\tag 4.9
$$
The Lie derivative of a one form is defined by
$$
{\Cal L}_T \alpha_a =T\alpha_a+(\pa_a T^c) \alpha_c,\tag 4.10
$$
The Lie derivative also commute with exterior differentiation,
$[{\Cal L}_T, d]=0$ so
$$
{\Cal L}_T \pa_a q=\pa_a T q\tag 4.11
$$
if $q$ is a function.
The Lie derivative of a two form is given by
$$
{\Cal L}_T \beta_{ab} =T\beta_{ab}
+(\pa_a T^c) \beta_{cb}+(\pa_b T^c) \beta_{ac}\tag 4.12
$$
Furthermore if $w$ is a one form and
$\curl w_{ab}= dw_{ab}=\pa_a w_b-\pa_b w_a$
then since the Lie derivative commutes with exterior differentiation:
$$
{\Cal L}_T \curl w_{ab}=\curl {\Cal L}_T w_{ab}\tag 4.13
$$
We will also use that the Lie derivative satisfies
Leibniz rule, e.g.
$$
{\Cal L}_T (\alpha_{c} W^c)=
({\Cal L}_T \alpha_{c}) W^c+\alpha_{c} {\Cal L}_T W^c,\qquad
{\Cal L}_T (\beta_{ac} W^c)=
({\Cal L}_T\beta_{ac}) W^c+\beta_{ac} {\Cal L}_T W^c. \tag 4.14
$$
Furthermore, we will also treat $D_t$ as if it was a Lie derivative and set
$$
{\Cal L}_{D_t}=D_t\tag 4.15
$$
Now of course this is not a space Lie derivative. It can however be interpreted
as a space time Lie derivative restricted to the space components.
What we use is that it satisfies the same properties (4.9)-(4.14)
as the other Lie derivatives we are considering.
 The reason we want to call it ${\Cal L}_{D_t}$
is simply a matter of that we will apply products of Lie derivatives and $D_t$
and since they behave in exactly the same way
it is more efficient to have one notation for them.

%\subheading{Lie Derivatives}

%\subheading{The modified Lie derivative and the divergence free condition}

The modification of the Lie derivative
$$
\tilde{\Cal L}_U W={\Cal L}_U W+
(\div U) W,\qquad \tag 4.16
$$
preserves the divergence free condition:
$$
\div \tilde{\Cal L}_U W=\tilde{U}\div W,
\qquad\text{where}\qquad \tilde{U}f=Uf+(\div U)f. \tag 4.17
$$
if $f$ is a function. (4.16) is invariant and (4.17) holds for any vector field $U$.
However, since we are considering Lie derivatives only with respect to
the vector fields constructed above and only expressed in the Lagrangian coordinates it is
simpler to use the modification
$$
\hat{\Cal L}_U W=\kappa^{-1} {\Cal L}_U (\kappa W)= {\Cal L}_U  W+(U\sigma)W,
\qquad \text{where}\qquad \sigma=\ln \kappa\tag 4.18
$$
Due to (4.1), $\div S=\kappa^{-1}\pa_a(\kappa S^a)=S\sigma$,
if $S$ is any of the tangential vector fields and
$\div R=R\sigma+n $, if $R$ is the radial vector field. For any of out tangential vector fields
it follows that
$$
\div \hat{\Cal L}_U W=\hat{U}\div W,
\qquad\text{where}\qquad \hat{U}f=Uf+(U\sigma)f=\kappa^{-1}U(\kappa f). \tag 4.19
$$
This has several advantages. The commutators satisfy $[\hat{\Cal L}_U,\hat{\Cal L}_T]
=\hat{\Cal L}_{[U,T]}$, since this is true for the usual Lie derivative.
Furthermore, this definition is constant with our previous definition of $\hat{D}_t$.

However, when applied to one forms we want to use the regular definition of
the Lie derivative.
Also, when applied to two forms most of the time we use the regular definition:
However, when applied to two forms it turns out to be sometimes
convenient to use the opposite  modification:
$$
\check{\Cal L}_T \beta_{ab}={\Cal L}_T \beta_{ab}-(U\sigma) \beta_{ab},\tag 4.20
$$
We will most of the time apply the Lie derivative to
products of the form $\alpha_a=\beta_{ab} W^b$:
$$
{\Cal L}_T \big(\beta_{ab} W^b\big)=
(\check{\Cal L}_T \beta_{ab})W^b+\beta_{ab} \hat{\Cal L}_T W\tag 4.21
$$
since the usual Lie derivative satisfies Leibniz rule.
 Using the modified Lie derivative
we indicated in \cite{L2} how to extend
the existence theorem in \cite{L1}
to the case when $\kappa $ is no longer constant,
i.e. $D_t\sigma=\div V\neq 0$.
This will be carried out in more detail here.

Let ${\Cal U}=\{ U_i\}_{i=1}^M $ be some labeling of our family of vector fields.
We will also use multindices $I=(i_1,...,i_r)$ of length $|I|=r$.
Let ${U}^I={U}_{i_1}\cdot\cdot\cdot {U}_{i_r}$ and
${\Cal L}_U^I={\Cal L}_{U_{i_1}}\cdot\cdot\cdot {\Cal L}_{U_{i_r}}$,
where ${\Cal L}_U$ is the Lie derivative. Similarly let
$\hat{U}^I f=\hat{U}_{i_1}\cdot\cdot\cdot \hat{U}_{i_r}f =\kappa^{-1} U^I(\kappa f)$ and
$\hat{\Cal L}_U^I W=\hat{\Cal L}_{U_{i_1}}\cdot\cdot\cdot \hat{\Cal L}_{U_{i_r}}W=
\kappa^{-1}{\Cal L}_U^I(\kappa W)$,
where $\hat{\Cal L}_U$ is the modified Lie derivative.
Sometimes we will also write ${\Cal L}_U^I$, where $U\in {\Cal S}_0$ or $I\in {\Cal S}_0$,
meaning that $U_{i_k}\in {\Cal S}_0$ for all of the indices in $I$.

%\subheading{Commutators with Lie derivatives with respect to tangential vector fields}
We will now calculate commutator between Lie derivatives and the operator
defined in section 3, i.e. the normal operator and the projected
multiplication operators. It is easier to calculate the commutator with Lie derivatives of
these operators considered as operators with values in the one forms.
The one form $w$ corresponding to the vector fields
$W$ is given by lowering the indices
$$
w_a=\underline{W}_a=g_{ab}W^b\tag 4.22
$$
For an operator $B$ on vector fields we denote the corresponding
operator with values in the one forms by $\underline{B}$. These
are related by
$$
\underline{B} W_a=g_{ab} BW^b,\qquad BW^a=g^{ab}\underline{B}_a\tag 4.23
$$
Most operators that
we consider will map onto the divergence free vector fields so
we will project the result afterwards to stay in this class.
Furthermore, in order to preserve the divergence free condition we will
use the modified Lie derivative. If the modified Lie derivative is applied to
a divergence free vector field then the result is divergence free so projecting
after commuting does not change the result.
As pointed out above, for our operators it is easier to commute Lie derivatives
with the corresponding operators from the divergence free vector fields to
the one forms.
Let $B_T$ be defined by
$$
{B}_T W^a=P\big(g^{ab}\big({\Cal L}_T \underline{B}W_b
-\underline{B}_b\hat {\Cal L}_{T}W\big)\big)\tag 4.24
$$
In particular if $B$ is a multiplication operator
 $\underline{B}_a W=P(\beta_{ab} W^b)=\beta_{ab}W^b-\pa_a q$, where $q$ vanishes
 on the boundary is chosen so that $\div BW=0$ then
$$
{\Cal L}_T \underline{B}_{a} W= \beta_{ab}\hat{\Cal L}_T W^b+
(\check{\Cal L}_T \beta_{ab})W^b+\pa_a Tq\tag 4.25
$$
and if we project to the divergence free vector fields then the term $\pa_a Tq$
vanishes since if $T$ is a tangential vector field then $Tq=0$ as well.
It therefore follows that $B_T$ is another multiplication operator:
$$
\underline{B}_T W_a =P\big((\check{\Cal L}_T \beta_{ab})W^b\big)\tag 4.26
$$
In particular, we will denote the time derivative of an operator
by $\dot{B}=B_{D_t}$ and for a multiplication operator this is
$$
\dot{B}W=B_{D_t} W= P (\big( \check{D}_t \beta_{ab}) W^b\big)\tag
4.27
$$

If $B$ maps on to the divergence free vector fields
$$
\hat{\Cal L}_T BW^a=\hat{\Cal L}_T (g^{ab}\underline{B}_a W)
=(\hat{\Cal L}_T g^{ab})\underline{B}_a W+ g^{ab}{\Cal L}_T
\underline{B}_a W\tag 4.28
$$
Here $\hat{\Cal L}_T g^{ab}=-g^{ac} g^{bd}\check{\Cal L}_T g_{cd}$.
If $B$ maps onto the divergence free vector fields then $\hat{\Cal L}_T B$
is also divergence free so the left hand side is unchanged if we project:
$$
\hat{\Cal L}_T BW^a=-P\big(g^{ab}(\check{\Cal L}_T
g_{bc})\underline{B}W^c\big) +P\big(g^{ab}\big({\Cal L}_T
\underline{B}_a W-\underline{B}_a \hat{\Cal L}_T W\big) \big)
+B\hat{\Cal L}_T W^a \tag 4.29
$$
By (4.26) applied the $\underline{G}_{ab}=P(g_{ab}W^b)$ we see
that $G_T W=P\big((g^{ab}\check{\Cal L}_T g_{bc}) W^c)\big)$ so
the first term in the right of (4.29) is $G_T BW^a$. The second
term is by definition (4.24) $B_T W$ so we get
$$
\hat{\Cal L}_T BW=B\hat{\Cal L}_T W +B_T W -G_T BW\tag 4.30
$$

The most important property of the projection is that it almost
commutes with Lie derivatives with respect to tangential vector fields.
If $\underline{P} u_a=u_a-\pa_a p_U$ then
$$
\underline{P}{\Cal L}_T \underline{P} u_a=\underline{P} {\Cal
L}_{T} u_a \tag 4.31
$$
since ${\Cal L}_T \pa_a p_{U}= \pa_a T p_{U}$ vanishes when we project again
since $T p_U$ vanishes on the boundary.
We have just used this fact above.
We have already calculated commutators between Lie derivatives and the
multiplication operators so let us now also calculate the commutator
between the Lie derivative with respect to tangential vector fields and
the normal operator.
Recall that the normal operator is defined by
$A_f W^a=g^{ab}\underline{A}_f W_b $, where
$$
\underline{A}_f W_a=-\pa_a\big( (\pa_c f) W^c-q\big),\qquad
\triangle \big( (\pa_c f) W^c-q\big)=0, \qquad
q\big|_{\pa\Omega}=0\tag 4.32
$$
and $f$ was function that vanished on the boundary.
Since the Lie derivative commutes with exterior differentiation it follows that
$$
{\Cal L}_T \underline{A}_f W_a= -\pa_a\big( (\pa_c f)\hat{\Cal
L}_T W^c + (\pa_c \check{T} f)W^c+(\pa_c T\sigma) f
W^c-Tq\big)\tag 4.33
$$
Since $q$ vanishes on the boundary it follows that $Tq$  also
vanish on the boundary and so does $(\pa_c T\sigma) f W^c$.
Therefore the last two terms vanish when we project again so we get
$$
P\big(g^{ab}{\Cal L}_T \underline{A}_f W_b\big)
=P\big(g^{ab}\underline{A}_f \hat{\Cal L}_T W_b\big) +
P\big(g^{ab}\underline{A}_{\check{T}f} W_b\big)\tag 4.34
$$
Let us now change notation so $A=A_p$, where $p$ is the pressure.
Then we have just calculated $A_T$ defined by (4.24) to be
$A_T=A_{\check{T}p}$, i.e.
$$
A_T=A_{\check{T} p},\qquad\text{if}\qquad A=A_{p}\tag 4.35
$$
In particular, if $T=D_t$ is the time derivative we will use the
notation $\dot{A}=A_{D_t}$ which then is
 $$
 \dot{A} W=A_{D_t }W= A_{ \check{D}_t p}W \tag 4.36
$$
Exactly the same formulas hold for $A_f^\varepsilon$. By (3.14)
$$
\underline{A}_f^\varepsilon W_a=
 (f/\rho) \chi_\varepsilon^\prime (\pa_a\rho) \,  (\pa_c \rho) W^c\!-\pa_a q,
\qquad \triangle q=\kappa^{-1}\pa_a\big( \kappa g^{ab}
(f/\rho) \chi_\varepsilon^\prime (\pa_b\rho) \,  (\pa_c \rho) W^c\big),
\quad q\big|_{\pa\Omega}\!=0\tag 4.37
$$
where $\rho=\rho(d)$, $d(y)=\operatorname{dist}{(y,\pa\Omega)}$. It follows that
$T\rho=$, if $T\in {\Cal T_0}$. Furthermore $S\in {\Cal S}_1={\Cal S}\setminus
{\Cal S}_0$  vanishes close to the boundary when $d(y)\leq d_0/2$
and $\chi^\prime_\varepsilon =0$ when $d(y)\geq \varepsilon$ so it follows that
$$
{\Cal L}_T \underline{A}_f^\varepsilon W_a=
 ( (\check{T}f)/\rho) \chi_\varepsilon^\prime
 (\pa_a\rho) \,  (\pa_c \rho) W^c\!-
 (f/\rho) \chi_\varepsilon^\prime (\pa_a\rho) \,  (\pa_c \rho) \hat{\Cal L}_T W^c\!
-\pa_a T q . \tag 4.38
$$
Hence
$$
P\big(g^{ab}{\Cal L}_T \underline{A}_f^\varepsilon W_b\big)
=P\big(g^{ab}\underline{A}_f^\varepsilon \hat{\Cal L}_T W_b\big) +
P\big(g^{ab}\underline{A}_{\check{T}f}^\varepsilon W_b\big)\tag 4.39
$$

We can now also calculate higher order commutators:
\demo{Definition {4.}1} If $T$ is a vector fields let $B_T$ be defined by (4.24).
If $T$ and $S$ are two tangential vector fields we define $B_{TS}=(B_S)_T$
to be the operator obtained by first using (4.24) to define $B_S$ and then
define $(B_S)_T $ to be the operator obtained from (4.24) with $B_S$ in place
of $B$. Similarly if $S^I=S^{i_2}\cdot\cdot\cdot S^{i_r}$ is a product
of $r=|I|$ vector fields then we define
$$
B_I= \big(\cdot\cdot\cdot
(B_{S^{i_1}})\cdot\cdot\cdot\big)_{S^{i_r}}\tag 4.40
$$
\enddemo

If $B$ is a projected multiplication operator
$BW^a=P\big(g^{ab}\beta_{bc}W^c\big)$ then
$$
B_I W=P\big(g^{ab}(\check{\Cal L}_T^I \beta_{bc}) W^c\big).\tag
4.41
$$
In particular if
$GW^a=P\big(g^{ab}g_{bc}W^c\big)$ then
$$
G_I W=P\big(g^{ab}(\check{\Cal L}_T^I g_{bc}) W^c\big).\tag 4.42
$$
If $A$ is the normal operator then
$$
{A}_I W^a=P\big( g^{ab}\pa_b\big( (\pa_c \check{T}^I p)
W^c\big)\,\big)\tag 4.43
$$

With $B_T$ as in (4.4) we have proven that if $B$ maps onto the divergence
free vector fields then
$$
\hat{\Cal L}_T BW=B W_T+B_T W-G_T BW,\qquad W_T=\hat{\Cal L}_T
W\tag 4.44
$$
Repeating this, gives for a product of modified Lie derivatives:
$$
\hat{\Cal L}_T^{I} BW=c_I^{\, I_1...I_k} G_{I_3}\cdot\cdot\cdot
G_{I_k} B_{I_1} W_{I_2} \qquad W_J=\hat{\Cal L}_T^J W\tag 4.45
$$
where the sum is over all combinations of $I=I_1+...+I_k$,
and $c_I^{\, I_1...I_k}$ are some constants, such that
$c^{\, I_1...I_k}_I=1$ if $I_1+I_2=I$. Let us then also introduce the notation
$$
G_I^{\, I_1 I_2}=c_I^{\, I_1...I_k} G_{I_3}\cdot\cdot\cdot G_{I_k},\tag
4.46
$$
where the sum is over all combination such that $I_3+...I_k=I-I_1-I_2$.
With this notation we can write (4.41)
$$
\hat{\Cal L}_T^{I} BW=G_I^{\, I_1 I_2} B_{I_1} W_{I_2}\tag 4.47
$$
where again $G_I^{\, I_1 I_2}=1$ if $I_1+I_2=I$. Also let
$$
\tilde{G}_I^{\,I_1...I_k}=0,\quad \text{if} \quad I_2=I, \qquad
\text{and}\qquad \tilde{G}_I^{\,I_1...I_k}=G_I^{\,I_1...I_k},
\quad\text{otherwise}. \tag 4.48
$$
Then we also have
$$
\hat{\Cal L}_T^{I} BW=B W_I+\tilde{G}_I^{\,I_1 I_2} B_{I_1}
W_{I_2}\tag 4.49
$$

\head{5. Estimating derivatives of a vector field in terms of the curl,
the divergence and tangential
derivatives or the normal operator}\endhead
The first part of the lemma below says that one can get a point wise estimate
of any first order derivative of a vector field by the curl, the divergence
and derivatives that are tangential at the boundary.
The second part say that one can get $L^2$ estimates with
a normal derivative instead of tangential derivatives.
The last part says that we can get the estimate for the normal derivative
from the normal operator. The lemma is formulated
in the Eulerian frame, i.e. in terms the Euclidean coordinates.
Later we will reformulate it in the Lagrangian frame and get
similar estimates for higher derivatives.

%\subheading{Estimates of a derivative of a vector field in terms of
% the curl, the divergence and tangential derivatives or the normal operator}
\proclaim{Lemma {5.}1} Let $\tilde{{\Cal N}}$ be a vector field that is
equal to the  normal ${\Cal N}$
at the boundary $\pa{\Cal D}_t$ and satisfies
$ |\tilde{{\Cal N}}|\leq 1$ and $|\pa\tilde{{\Cal N}}|\leq K$. Let
$q^{ij}=\delta^{ij}-\tilde{\Cal N}^i\tilde{\Cal N}^j$.
Then
$$\align
\!\!\!\!\!\!\!\!\!\!\!\!\! |\pa \beta|^2 &\leq
C\big( q^{kl}\delta^{ij} \pa_k \beta_i\, \pa_l \beta_j
+|\curl \beta|^2 +|\div \beta|^2 \big)
\!\!\!\!\!\!\!\!\!\!\!\!\!\!\!\!\!\!\!\!\!\!\!\!\!\tag 5.1\\
\!\!\!\!\!\!\!\!\!\!\!\!\!\int_\dt\! |\pa \beta|^2 dx&\leq
C\int_\dt\!\!\big( \delta^{ij} \tilde{\Cal N}^k \tilde{\Cal
N}^l \pa_i \beta_k\, \pa_j \beta_l
+|\curl \beta|^2\!+|\div \beta|^2\!+K^2|\beta|^2\big)\, dx
\!\!\!\!\!\!\!\!\!\!\!\!\!\!\!\!\!\!\!\!\!\!\!\!\!\!\!\!\!\!\!\!\!\!\!\!
\!\!\!\!\!\!\!\!\!\tag 5.2
\endalign
$$
Suppose that $\delta^{ij} \alpha_j$ is another vector field that is normal
at the boundary and let $A\beta_i=\pa_i( \alpha_k \beta^k-q)$ and $q$ is chosen so that
$\div A\beta=0$ and $q|_{\pa\Omega}=0$. Then
$$
\int_\dt\!\! \delta^{ij} \alpha_k \alpha_l\,
 \pa_i \beta^k\, \pa_j \beta^l\, dx
\leq
C\int_\dt\!\!\big( \delta^{ij} A\beta_i\, A\beta_j
+|\alpha|^2\big(|\curl \beta|^2\!+|\div \beta|^2\big)
\!+|\pa\alpha|^2
|\beta|^2\big)\, dx\tag 5.3
$$
\endproclaim
\demo{Proof}
(5.1) follows from the point wise estimate
$$\align
\delta^{ij}\delta^{kl} w_{ki}w_{lj}&\leq
C\big(\delta^{ij}q^{kl}w_{ki}w_{lj}+|\hat{w}|^2+(\tr{w})^2\big)\tag 5.4\\
\delta^{ij}\delta^{kl} w_{ki}w_{lj}&\leq
C\big(\tilde{\Cal N}^i\tilde{\Cal N}^j \delta^{kl}w_{ki}w_{lj}
+(q^{ij}q^{kl}-q^{ik}q^{jl})w_{ki}w_{lj}
+|\hat{w}|^2+(\tr{w})^2\big)\!\!\!\!\!\!\!\!\!\tag 5.5
\endalign
$$
where $\hat{w}_{ij}=w_{ij}-w_{ji}$ is the antisymmetric part and
$\tr w=\delta^{ij} w_{ij}$ is the trace. To prove (5.4)-(5.5) we may assume that
$w$ is symmetric and traceless. Writing
$\delta^{ij}=q^{ij}+\tilde{\Cal N}^i \tilde{\Cal N}^j$ we see that (5.4) for such
tensors follows from the estimate
$\tilde{\Cal N}^i \tilde{\Cal N}^j \tilde{\Cal N}^k \tilde{\Cal N}^l w_{ki} w_{lj}
=(\tilde{\Cal N}^i \tilde{\Cal N}^k w_{ki})^2=(q^{ik} w_{ki})^2\leq
n q^{ij} q^{kl} w_{ki} w_{lj}$. (This says that
$(\tr(QW))^2\! \leq\! n \tr(QWQW)$ which is obvious if one writes it out and use the
symmetry.) (5.5) follows since
$(\delta^{ij}q^{kl}\!-\tilde{\Cal N}^i\tilde{\Cal N}^j \delta^{kl}) w_{ki}w_{lj}\!=
(q^{ij}q^{kl}\!-\tilde{\Cal N}^i\tilde{\Cal N}^j\tilde{\Cal N}^k\tilde{\Cal N}^l)
 w_{ki}w_{lj}=(q^{ij}q^{kl}-q^{ik}q^{jl} )w_{ki}w_{lj}$.
(5.2) follows from (5.5) and
integration by parts using that the boundary terms vanish, since we assumed that
$\tilde{\Cal N}={\Cal N}$ there, and that
$(q^{ij}q^{kl}-q^{ik}q^{jl}) \beta_i \, \pa_k\pa_j\beta_l=0$:
$$
\int_\dt  (q^{ij}q^{kl}-q^{ik}q^{jl})\pa_k \beta_i \, \pa_j\beta_l \, dx
=-\int_\dt  \pa_k(q^{ij}q^{kl}-q^{ik}q^{jl})\, \beta_i \, \pa_j\beta_l\, dx\tag 5.6
$$
We have $A\beta_i =(\pa_i\alpha_k)\beta^k+\alpha_k\pa_i\beta^k-\pa_i q$
so to prove (5.3) we must estimate $\|\pa q\|_{L^2}$.
Since $0=\pa_i A\beta^i=\triangle (\alpha_k\beta^k)-\triangle q$ it
follows that $\triangle q=\triangle\big(\alpha_k\beta^k\big)
=2\pa_i\big( (\pa^i\alpha_k)\beta^k\big)+\alpha_k\triangle \beta^k
-(\triangle \alpha_k)\beta^k $
and $\alpha_k\triangle \beta^k
=\pa_i\big(\alpha^i\div\beta+\alpha_k\curl \beta^{ik}\big)
-\div\alpha\,\div\beta -(\pa_k \alpha_i) \pa^k\beta^i+(\pa_k\alpha_i)\pa^i\beta^k$,
and hence
$\triangle (\alpha_k\beta^k)=
\pa_i\Big( 2(\pa^i\alpha_k)\beta^k+\alpha^i\div \beta
+\alpha_k\curl \beta^{ik}-\div\alpha \,\beta^i
-\curl \alpha^i_{\,\,k}\, \beta^k\Big)$. It follows that
$$
\int_\Omega |\pa q|^2\, dx=-\int_\Omega q\triangle q\, dx
=-\int_\Omega q \pa_i \big(\alpha^i \div\beta
+\alpha_k \curl \beta^{ik}
+(\pa^i  \alpha_k+\pa_k\alpha^i)\beta^k
-\div\alpha\, \beta^i\big)\, dx\tag 5.7
$$
and integrating by parts again gives
$\|\pa q\|_{L^2}\leq C\big(\|\,|\alpha|\div\beta\|_{L^2}
+\|\,|\alpha|\curl\beta\|_{L^2}+\|\,|\pa\alpha|\beta\|_{L^2}\big)$.
\qed\enddemo

%\subheading{Bounds of Lie derivatives of a vector field by
%the curl, the divergence and tangential derivatives}
\demo{Definition {5.}1} For ${\Cal V}$ any of the family of vector fields introduced
in \cite{L1} and for $\beta$ a two form, a one form, a function or a vector field
we define
$$
|\beta|_r^{\Cal V}=\sum_{|I|\leq r ,\, I\in {\Cal V}}
|{\Cal L}_U^I \beta\,|,\qquad\quad
[\beta]_r^{\Cal V}=\sum_{r_1+...r_k\leq r,\, r_i\geq 1}
|\beta|_{r_1}^{\Cal V}\cdot\cdot\cdot |\beta|_{r_k}^{\Cal V}
\tag 5.8
$$
and $[\beta]_0^{\Cal V}=1$.
Furthermore let
$$
|\beta|_r=\sum_{|\alpha|\leq r}|\pa_y^\alpha \beta|\tag 5.9
$$
\enddemo
If $\beta$ is a function then ${\Cal L}_U \beta =U\beta$ and
in general it is equal to this plus terms proportional to $\beta$.
Hence (5.8) is equivalent to just the sum $\sum_{|I|\leq r,\,I\in {\Cal V}}
|U^I \beta| $ . In particular if ${\Cal R}$ denotes the family of space
vector fields then $|\beta|_r^{\Cal R}$ is equivalent to
$|\beta|_r$ with a constant of equivalence independent of the metric.
Note also that if $\beta$ is the one form $\beta_a=\pa_a q$ then ${\Cal L}_U^I\beta=
\pa U^I q$ so $|\pa q|_r^{\Cal V}=\sum_{|I|\leq r,\, I\in {\Cal V}}|\pa U^I q|$.

\demo{Definition {5.}2 } Let $c_1$ be a constant such that
$$
|\pa x/\pa y|^2+|\pa y/\pa x|^2\leq c_1^2,\qquad\quad
\sum_{a,b=1}^n\big(|g_{ab}|+|g^{ab}|\big)\leq n c_1^2,\tag 5.10
$$
and let $K_1$ denote a continuous function of $c_1$.
\enddemo
We note that the second condition in (5.10) follows from the first and the
first follows from the second with a larger constant.
We remark that this condition is fulfilled initially since we are composing with a
diffeomorphism. Furthermore, for solution of Euler's equations, $\div V=0$,
so the volume form $\kappa$ is preserved and hence an upper bound for the metric also implies
a lower bounded for the eigenvalues and an upper bound for the  inverse of
the metric follows.

In what follows it will be convenient to consider the norms of
$\hat{\Cal L}_U^I W=\kappa^{-1}{\Cal L}_U^I (\kappa W)$ if $W$ is a vector field
and of $\check{\Cal L}_U^I g=\kappa {\Cal L}_U^I(\kappa^{-1} g)$, if $g$ is the metric.
The reason for this is simply that $\div\, (\hat{\Cal L}_U^I W)
=\hat{U}^I \div W$ and ${\Cal L}_U^I\curl w=\curl\, ({\Cal L}_U^I w)$
and when we lower indices $w_a=g_{ab} W^b=(\kappa^{-1} g_{ab})(\kappa W^b)$ and
apply the Lie derivative to the product we get ${\Cal L}_U w_a
=(\check{\Cal L}_U g_{ab})W^b +g_{ab}\hat{\Cal L}_U W^b$.

\proclaim{Lemma {5.}2} Let $W$ be a vector field and let
$w_a=g_{ab} W^b$ be the corresponding one form.
Let $\kappa=\det{(\pa x/\pa y)}=
\sqrt{\det{g}}$. Then
$$
|\kappa|+|\kappa^{-1}|\leq K_1,
\qquad |U^I\kappa|+|U^I \kappa^{-1}|\leq K_1 c^{I_1...I_k} |U^{I_1}g|\cdot\cdot\cdot
|U^{I_k} g|\tag 5.11
$$
where the sum is over all $I_1+...I_k=I$.

With notation as in Definition {5.}1 and section 4 we have
$$
|\kappa W|_r^{\Cal R}
\leq K_1\big( |\curl w|_{r-1}^{\Cal R}+|\kappa\div W|_{r-1}^{\Cal R}
+|\kappa W|_r^{\Cal S}+\sum_{s=0}^{r-1} |g/\kappa|_{r-s}^{\Cal R}
 | \kappa W|_{s}^{\Cal R}\big)\tag 5.12
$$
We also have
$$
|\kappa W|_{r}^{\Cal R}\leq K_1\sum_{s=0}^{r}\,\, [g/\kappa]_{s}^{\Cal R}
\big(|\curl w|_{r-1-s}^{\Cal R}+|\kappa\div W|_{r-1-s}^{\Cal R}
+|\kappa W|_{r-s}^{\Cal S}\big),\tag 5.13
$$
where for $s=r$ we use the convention that
$|\curl w|_{-1}^{\Cal V}=|\kappa\div W|_{-1}^{\Cal V}=0$.
Furthermore (5.12)-(5.13) holds without the factors
 $\kappa$ and $1/\kappa$, i.e.
$$
|W|_{r}^{\Cal R}\leq K_1\sum_{s=0}^{r} \,\, [g]_{s}^{\Cal R}
\big(|\curl w|_{r-1-s}^{\Cal R}+|\div W|_{r-1-s}^{\Cal R}
+|W|_{r-s}^{\Cal S}\big),\tag 5.14
$$
(5.12)-(5.13) also holds for the vector field $W$ replaced by a
one form $w$, i.e.
$$
|w|_{r}^{\Cal R}\leq K_1\sum_{s=0}^{r} \,\, [g]_{s}^{\Cal R}
\big(|\curl w|_{r-1-s}^{\Cal R}+|\div W|_{r-1-s}^{\Cal R}
+|w|_{r-s}^{\Cal S}\big),\tag 5.15
$$
Moreover, the inequalities (5.12)-(5.15) also hold with $({\Cal R},{\Cal S})$
replaced by $({\Cal U},{\Cal T})$.
\endproclaim
\demo{Proof}If $\sigma=\ln \kappa=(\ln\det{g})/2$ then
 $U\sigma=\tr {\Cal L}_U g/2=g^{ab}{\Cal L}_U g_{ab}/2$ and
${\Cal L}_U g^{ab}=-g^{ac}g^{bd}{\Cal L}_U g_{cd}$.
An easy consequence of Lemma {5.}1, see \cite{L1}, is:
In the Lagrangian frame we have, with
$w_a=\underline{W}_a=g_{ab} W^b$,
$$\align
|\hat{\Cal L}_U W|&\leq K_1\Big(|\curl\,\underline{W}\,|+|\div {W}|
+\tsize{\sum_{S\in{\Cal S}}} |\hat{\Cal L}_S W|+[g]_1|W|\Big),
\qquad U\in{\Cal R},\tag 5.16\\
|\hat{\Cal L}_U W|&\leq K_1\Big(|\curl\, \underline{W}\,|+|\div {W}|
+\tsize{\sum_{T\in{\Cal T}}} |\hat{\Cal L}_T W|+[g]_1|W|\Big),
\qquad U\in{\Cal U}.\tag 5.17
\endalign
$$
where $[g]_1=1+|\pa g|$. Furthermore
$$
|\pa W|\leq K_1\Big(|\hat{\Cal L}_R W|
+\tsize{\sum_{S\in{\Cal S}}} |\hat{\Cal L}_S W|+[g]_1 |W|\Big)
\tag 5.18
$$
When $d(y)\leq d_0$ we may replace the sums over ${\Cal S}$
by the sums over ${\Cal S}_0$ and the sum over ${\Cal T}$
by the sum over ${\Cal T}_0$.
In \cite{L1} this was proven for $\hat{\Cal L}_U$ replaced
by ${\Cal L}_U$, but the difference is just a lower order term.

We claim that
$$
\sum_{|I|=r,U\in{\Cal R}}\!\!
|\hat{\Cal L}_U^I W|\leq K_1\!\!\!\!\sum_{|J|=r-1, U\in{\Cal R} }
\!\!\!\big(|\curl\underline{\hat{\Cal L}_U^J W}|
+|\div\hat{\Cal L}_U^J W|+ [g]_1 |\hat{\Cal L}_U^J W|\big)
+K_1\!\!\!\!\sum_{|I|=r,S\in {\Cal S}}\!\!|\hat{\Cal L}_S^I W|\tag 5.19
$$
First we note that there is noting to prove if $d(y)\geq d_0$ since then
${\Cal S}$ span the full tangent space. Therefore, it suffices to prove (5.19)
when $d(y)\leq d_0$ and with ${\Cal S}$ replaced by ${\Cal S}_0$ and
${\Cal R}$ replaced by ${\Cal R}_0$.
Then (5.19) follows from (5.16) if $r=1$ and assuming that its true for
$r$ replaced by $r\!-\!1$ we will prove that it holds for $r$.
If we apply (5.16) to $\hat{\Cal L}_{U}^J W$, where $|J|=r\!-\!1$, we get
$$
|\hat{\Cal L}_U \hat{\Cal L}_U^J W|
\leq K_1\big(|\curl\underline{\hat{\Cal L}_U^J W}|
+|\div\hat{\Cal L}_U^J W|+\sum_{S\in{\Cal S}}|\hat{\Cal L}_S \hat{\Cal L}_U^JW|
+[g]_1|\hat{\Cal L}_U^J W|\big).\tag 5.20
$$
If $\hat{\Cal L}_U^J$ consist of all tangential derivatives then
it follows that $|\hat{\Cal L}_U \hat{\Cal L}_U^J W|$ is bounded by the right hand
side of (5.19). If $\hat{\Cal L}_U^J$ does not consist of only
tangential derivatives then, since $[\hat{\Cal L}_R, \hat{\Cal L}_S]=\hat{\Cal L}_{[R,S]}=0$,
if $S\in {\Cal S}_0$, we can write
$\hat{\Cal L}_S\hat{\Cal L}_U^J W=\hat{\Cal L}_U^K \hat{\Cal L}_{S^\prime}W$,
for some $S^\prime\in {\Cal S}_0$.
If we now apply (5.19) with $r$ replaced by $r-1$ to $\hat{\Cal L}_{S^\prime} W$,
(5.19) follows also for $r$.

In the Lemma we have ${\Cal L}_U^I\curl w=\curl{\Cal L}_U^I w$
which however is different from $\curl\, \underline{\hat{\Cal L}_U^I W}$,
We have:
$$
{\Cal L}_U^J w_a={\Cal L}_U^J (g_{ab}W^b)=-g_{ab}\hat{\Cal L}_U^J W^b
+\tilde{c}^{\,J}_{\,J_1 J_2} g^{J_1}_{ab}\hat{\Cal L}_U^{J_2} W^b,
\qquad\text{where}\quad g^{J}_{ab}=\check{\Cal L}_U^J g_{ab}\tag 5.21
$$
where the sum is over all $J_1+J_2=J$ and $\tilde{c}^{\,J}_{\, J_1 J_2}=1$ for $|J_2|<|J|$
$c^{\,J}_{J_1 J_2}=0$ if $J_2=J$.
It follows that
$$
|\curl \underline{\hat{\Cal L}_U^J W}-\curl{\Cal L}_U^J w|\leq
2\tilde{c}^J_{\, J_1 J_2} \big( |\pa g^{J_1}| |\hat{\Cal L}_U^{J_2} W|
+|g^{J_1}| |\pa\hat{\Cal L}_U^{J_2} W|\big),\qquad
|J_2|<|J|,\tag 5.22
$$
where the partial derivative can be estimated by Lie derivatives.
Furthermore, in the Lemma we have $|U^I(\kappa \div W)|=
\kappa^{-1}|\hat{U}^I \div W|=\kappa^{-1}| \div\hat{\Cal L}_U^I W|$.
(5.13) follows by induction from (5.12).
\qed\enddemo

\demo{Definition {5.}3} For ${\Cal V}$ any of the family of vector
fields introduced in \cite{L1} let
$$
\| W\|_{r}^{\Cal V}=\sum_{|I|\leq r,\, I\in{\Cal V}} \|{\Cal L}_U^I W\|,
\qquad\qquad \| W\|_{r,\infty}^{\Cal V}=
\sum_{|I|\leq r,\, I\in{\Cal V}} \|{\Cal L}_U^I W\|_{\infty}
\tag 5.23
$$
and let
$$
\|W\|_r=\sum_{|\alpha|\leq r}\|\pa_y^\alpha W\|\qquad\qquad
\|W\|_{r,\infty}=\sum_{|\alpha|\leq r}\|\pa_y^\alpha W\|_\infty
\tag 5.24
$$
where $\|W\|=\|W\|_{L^2(\Omega)}$, $\|W\|_\infty=\|W\|_{L^\infty(\Omega)}$.
\enddemo
It follows from the discussion after Definition {5.}1 and (5.11) that $\|W\|_r$
is equivalent to $\|W\|_r^{\Cal R}$ with a constant of equivalence
independent of the metric. As with the point wise estimates it will
sometimes be convenient to instead use $\|\hat{\Cal L}_U^I W\|
=\|\kappa^{-1}{\Cal L}_U^I(\kappa W)\|$. This in particular true for the
family of space tangential vector fields ${\Cal S}$.
However instead of introducing a special notation we then write
$\|\kappa W\|_r^{\Cal S}$. Since $\kappa$ is bounded from above and below by
a constant $K_1$ this is equivalent with a constant of equivalence $K_1$.
Furthermore, by interpolation
$\|\kappa W\|^{\Cal S}_r\leq K_1 (\|g\|_r\|W\|+\|W\|_r^{\Cal S})$
and  $\|W\|^{\Cal S}_r\leq K_1 (\|g\|_r\|W\|+\|\kappa W\|_r^{\Cal S})$,
and our inequalities anyway contain lower order terms of this form,
so the inequalities below are true either with or without $\kappa$.

\proclaim{Lemma {5.}3}  We have with a constant $K_1$ as in Definition {5.}1:
$$
\|W\|_{r}\leq K_1\big(\|\curl w\|_{r-1}+\|\kappa \div W\|_{r-1}+\|\kappa W\|_{r}^{\Cal S}
+K_1 \sum_{s=0}^{r-1}\|g\|_{r-s,\infty}\|W\|_s\big)\tag 5.25
$$
and, with the convention that $\|\curl w\|_{-1}+\|\div W\|_{-1}=0$,
$$
\|W\|_{r}\leq K_1\sum_{s=0}^r \|g\|_{r-s,\infty}
\big(\|\curl w\|_{s-1}+\|\kappa\div W\|_{s-1}+\|\kappa W\|_{s}^{\Cal S}\big)\tag 5.26
$$
\endproclaim
\demo{Proof} This follows from Lemma {5.}2 and the interpolation inequalities below
in Lemma {6.}2.
\qed\enddemo

We can also bound derivatives of a vector field by
the curl, the divergence and the normal operator:
\proclaim{Lemma {5.}4} Let $c_0>0$ be a constant such that
$|\na_N p|\geq c_0>0$, let $K_2$ and $K_3$ be constants such that
$\|\na_N p\|_{L^\infty(\pa\Omega)}\leq K_2$
and $\sum_{S\in{\Cal S}}\|\na_N  S p\|_{L^\infty(\pa\Omega)}\leq K_3$.
Then
$$
c_0\| \pa W\|\leq
C\big( \|AW\|+ K_2(\|\curl w\|+\|\div W\|)+(K_3+[g]_1)\|W\|\big)\tag 5.27
$$
\endproclaim
\demo{Proof} We want to express (5.2) and (5.3)
in the Lagrangian frame. We also want to pick an extension of the normal
to the interior. If $d(y)$ be the distance to the boundary
in the Lagrangian frame,
since $\Omega$ is the unit ball this is just $1-|y|$. Let $\chi_1(d)$
be a smooth function that is $1$ close to $0$ and $0$ when $d>1/2$.
If $u_c=\pa_c d$ then  $n_c=u_c/\sqrt{g^{ab}u_a u_b}$ is the unit
conormal at the boundary and $\tilde{n}_c=\chi_1(d) n_c$
defines an extension to the interior and $\tilde{N}^a=g^{ab}\tilde{n}_b$
is an extension of the unit normal to the interior.
Similarly, by the remarks in section 3, the normal operator only
depends on $\na_N p$ restricted to the boundary. Let us define
$\alpha_b=\chi_2(d) f \pa_b d$, where $f$ is a function that is equal
to $N^c\pa_c p=\na_N p$ at the boundary and extended to be constant along
rays through
the origin, and $\chi_2$ is a function that is $1$ on the support of $\chi_1$
and $0$ when $d>3/4$. Then $\underline{A} W^a=P(g^{ab}\pa_b( (\pa_c p) W^c))
=P(g^{ab}\pa_b ( \alpha_c W^c))$ by the remarks in section 3.
Now, in expressing (5.2) and (5.3) in the Lagrangian coordinates
partial differentiation becomes covariant differentiation so we will
pick up a constant coming from the Christoffel symbols, i.e. one derivative
of the metric $[g]_1=1+|\pa g|$. Similarly, one derivative of the normal
$N^a$ also gives rise to one derivative of the metric. Hence (5.2)
and (5.3) become
$$
\| \pa W\| \leq C\big( \| \chi_1 ( n_c \pa W^c)\|
+\|\curl w\|+\|\div W\|+ [g]_1 \| W\|\big)\tag 5.28
$$
and
$$
\|f  \chi_2 ( n_c \pa W^c)\|\leq
\big(\|A W\|+ \|f\curl w\|+\|f\div W\|+ [g]_1 \|f W\|+\|\,|\pa f|W\|\big)\tag 5.29
$$
Since $|f|\geq c_0$ and $\chi_2=1$ in a neighborhood of the support of $\chi_1$,
the lemma follows.
\qed\enddemo
By Lemma {5.}4 we have
$$
c_0\|\pa \hat{\Cal L}_S^J W\|
\leq K_3\big(\|\curl \underline{\hat{\Cal L}_S^J W}\|+\|\div\hat{\Cal L}_S^J W\|
+\| A\hat{\Cal L}_S^J W\|+\|\hat{\Cal L}_S^J W\|\big)\tag 5.30
$$
where $K_3$ is as in Definition {6.}1 and $c_0$ as in the physical
condition (1.6). Here, the curl of
 $(\underline{\hat{\Cal L}_{S}^J W})_{a}=g_{ab}
\hat{\Cal L}_S^J W^b$ is by (5.22) equal to the curl of ${\Cal L}_S^J w$
plus lower order terms.
In particular we see that we can get any space tangential derivative in this way so
we also get:
\proclaim{Lemma {5.}5} With $K_3$ as in Definition {6.}1 we have
$$
c_0 \|W\|_{r}\leq  K_3\big(\|\curl w\|_{r-1}
+\|\div W\|_{r-1}+\| W\|_{r-1,A}^{\Cal S}
+\sum_{s=0}^{r-1} \|g\|_{r-s,\infty}\|W\|_s
\big)\tag 5.31
$$
where
$$
\|W\|_{s,A}^{\Cal S}
=\sum_{|I|=s,I\in{\Cal S}}\| A\hat{\Cal L}_S^I W\|\tag 5.32
$$
\endproclaim

\head 6. Interpolation, the $L^\infty$
estimates for the pressure in terms of the coordinate
and the $L^\infty$ norms. \endhead
Let us now first state the interpolation inequalities that we will use:
 \proclaim{Lemma {6.}1} Let $\beta$ be a two
form, a function or a vector field. Let $\|\beta\|_r$ be
$L^2$-Sobolev norms and $\|\beta\|_{r,\infty}$ is the $C^k$ norms
on the unit ball $\Omega$ in $\bold{R}^n$.
Then
if $0\leq s\leq r$ and $j\geq 0$
$$\align
\|\beta\|_{j+s,\infty}&\leq
 C\|\beta\|_{j,\infty}^{1-s/r}\|\beta\|_{j+r,\infty}^{s/r}\tag 6.1 \\
\|\beta\|_{s}&\leq C\|\beta\|_{0}^{1-s/r}\|\beta\|_{r}^{s/r}\tag 6.2
\endalign
$$
\endproclaim
For a proof see e.g. \cite{H1,H2}, for the $L^\infty$ norm and
\cite{CL}, for the $L^2$ norms. ( (6.1) for $j>0$ follows from (6.1) for $j=0$
applied to $\pa_y^{\alpha}$ for $|\alpha|\leq j$. )A consequence is:
 \proclaim{Lemma {6.}2} With the same assumptions as in Lemma {6.}1 we have
$$
\align
\|\alpha\|_{j+r-s,\infty}\|\beta\|_{j+s,\infty}& \leq
\big(\|\alpha\|_{j,\infty}\|\beta\|_{j+r,\infty}
+\|\beta\|_{j,\infty}\|\alpha\|_{j+r,\infty}\big)\tag 6.3 \\
\|\beta\|_{r-s,\infty}\|W\|_{s}&\leq C\big(
\|\beta\|_{0,\infty}\|W\|_r +\|\beta\|_{r,\infty}\|W\|_0\big)\tag
6.4\\
\leftalignspace \|f_1\|_{j+s_1,\infty}\!\cdot\cdot\cdot
\|f_k\|_{j+s_k,\infty}&\leq
C\sum_{i=1}^k\|f_1\|_{j,\infty}\!\cdot\cdot\cdot
\|f_{i-1}\|_{j,\infty} \|f_i\|_{j+s_1+...+s_k,\infty}
\|f_{i+1}\|_{j,\infty}\!\cdot\cdot\cdot
\|f_k\|_{j,\infty}\rightalignspace \tag 6.5
\endalign
$$
\endproclaim
\demo{Proof} This follows from using Lemma {6.}1 on each factor
and the inequality $A^{s/r} B^{1-s/r}\leq A+B$, e.g.
 $$\multline
 \|\beta\|_{r-s,\infty}\|W\|_{s}\leq C\|\beta\|_{0,\infty}^{s/r}
 \|\beta\|_{r,\infty}^{1-s/r}\|W\|_0^{1-s/r} \|W\|_r^{s/r}
  \\
 =C \big(
\|\beta\|_{0,\infty}\|W\|_r\big)^{s/r}\big(\|\beta\|_{r,\infty}
\|W\|_0\big)^{1-s/r}\leq C \big(\|\beta\|_{0,\infty} \|W\|_r+
\|\beta\|_{r,\infty} \|W\|_0 \big). \endmultline \tag 6.6
$$
This proves (6.4). The proof of (6.3) is exactly the same,
(6.5) follows from (6.3) by induction.
 \qed\enddemo

Let us now introduce some notation to be used in subsequent sections.
We will derive tame estimate involving the higher norms of the coordinate $x$
with constants that are bounded if some lower
norms of the coordinate $x$ are bounded: Recall Definition {5.}2 of $c_1$:
$$
|\pa x/\pa y|^2+|\pa y/\pa x|^2\leq c_1^2,\qquad\quad
\sum_{a,b=1}^n\big(|g_{ab}|+|g^{ab}|\big)\leq n c_1^2,\tag 6.7
$$
and $K_1$ denotes a continuous function of $c_1$.
\demo{Definition {6.}1} Let $c_1$ be as in Definition {5.}2
and for $i=2,3,4$ let  $c_i\geq c_1$ be a constant such that
$$
\align
 \|x\|_{2,\infty}+\|\dot{x}\|_{1,\infty}&\leq c_2,\tag 6.8 \\
 \|x\|_{3,\infty}+\|\dot{x}\|_{2,\infty}+\|\ddot{x}\|_{1,\infty}
&\leq c_3,\tag 6.9\\
\|x\|_{4,\infty}+\|\dot{x}\|_{3,\infty}+\|\ddot{x}\|_{2,\infty}&\leq c_4
\tag 6.10
\endalign
$$
and let $K_i\geq 1$ denote a constant that depends continuously on $c_i$.
\enddemo

It now follows from using Lemma {6.}2:
\proclaim{Lemma {6.}3} With $K_1$ as in Definition {5.}1
$$
\| \pa y/\pa x\|_{r,\infty}\leq K_1 \| x\|_{r+1,\infty}.\tag 6.11
$$
If $\pa_i=\pa/\pa x^i=(\pa y^a/\pa x^i)\pa /\pa y^a$ and
$\alpha=(\alpha_1,...,\alpha_n)$ let $\pa^\alpha=\pa_{1}^{\alpha_1}\cdots
\pa_n^{\alpha_n}$. For any function $f$ we have
$$
\|\pa^\alpha f\|_{r,\infty}\leq
K_1 \big(\|f\|_{r+k,\infty} + \|x\|_{r+k,\infty}\|f\|_{1,\infty}\big),
\qquad k=|\alpha|. \tag 6.12
$$
Moreover,
$$
\align
\!\!\!\!\!\!\!\!\!\!
\|(\pa_{i_1} f_1)\!\cdots  (\pa_{i_n} f_n)&\|_{r,\infty}\leq
K_1\!\! \sum_{i=1}^n \|f_1\|_{1,\infty}\!\cdots \|f_{i-1}\|_{1,\infty}
\|f_i\|_{r+1,\infty}\|f_{i+1}\|_{1,\infty}\!\cdots \|f_{n}\|_{1,\infty}
\!\!\!\!\!\!\!\!\!\!\!\!\!\!\!\!\!\!\!\!\!\!\!\!\!\!\!\!\tag 6.13\\
 &\qquad\qquad\qquad
+K_1 \|x\|_{r+1,\infty}\|f_1\|_{1,\infty}\!\cdots\|f_n\|_{1,\infty},
\!\!\!\!\!\!\!\!\\
\!\!\!\!\!\!\!\!\!\! \|(\pa_{i_0}\pa_{i_1} f_1) (\pa_{i_2} f_2)\!\cdots
  (\pa_{i_n} f_n)&\|_{r,\infty}\leq
K_1\!\! \sum_{i=1}^n \|f_1\|_{1,\infty}\!\cdots \|f_{i-1}\|_{1,\infty}
\|f_i\|_{r+2,\infty}\|f_{i+1}\|_{1,\infty}\!\cdots \|f_{n}\|_{1,\infty}
\!\!\!\!\!\!\!\!\!\!\!\!\!\!\!\!\!\!\!\!\!\!\!\!\!\!\!\!\tag 6.14\\
 &\qquad\qquad\qquad
+K_1 \|x\|_{r+2,\infty}\|f_1\|_{1,\infty}\!\cdots\|f_n\|_{1,\infty},
\!\!\!\!\!\!\!\!
\endalign
$$
\endproclaim
\demo{Proof} Let $A$ be the matrix $\pa x^i/\pa y^a$. Using the {formula
for the derivative of a matrix $\pa_a A^{-1}=-A^{-1} (\pa_a A)A^{-1}$ we see
that $\pa_y^\alpha A^{-1}$ is a sum of terms of the form
$$
A^{-1}(\pa_y^{\alpha_1} A)A^{-1}\cdots (\pa_y^{\alpha_k} A)A^{-1},
\qquad |\alpha_1|+...+|\alpha_k|=|\alpha|=r\tag 6.15
$$
Since $|A^{-1}|\leq Cc_1$ we see from (6.5) with $j=0$ that this is bounded by
$K_1 |A|_{r,\infty}$ which proves (6.11). Now $\pa_y^{\gamma} \pa_x^\alpha f$
is sum of terms of the form
$$
A^{-1}(\pa_y^{\beta_1} A^{-1})\cdots(\pa_y^{\beta_{k-1}} A^{-1})
\pa_y^{\beta_k} (\pa_y f),\qquad |\beta_1|+...+|\beta_k|=|\gamma|+|\alpha|=r-1+k
\tag 6.16
$$
By (6.5) and what we used proved this is bounded by
$K_1\|\pa_y f\|_{r-1+k,\infty} +K_1\|A\|_{r-1+k,\infty} \|\pa_y f\|_{0,\infty}$
which proves (6.12). By (6.13) follows from (6.5)
with $j=0$ and (6.12) with $k=1$. (Note by (6.12)
$\|\pa f\|_{0,\infty}\leq K_1\|f\|_{1,\infty}$.) Similarly, by (6.5) with $j=0$
and (6.12) we can bound the left of (6.14) by
$$\multline
\big(\|f_1\|_{r+2,\infty}+\|x\|_{r+2,\infty}\big)
\|f_2\|_{1,\infty}\cdots\|f_n\|_{1,\infty}\\
+\sum_{i=2}^n (\|f_1\|_{2,\infty}+\|x\|_{2,\infty})
\|f_1\|_{1,\infty}\cdots \|f_{i-1}\|_{1,\infty}
\big( \|f_{i}\|_{r+1,\infty}+\|x\|_{r+1,\infty}\big)
\|f_{i+1}\|_{1,\infty}\cdots\|f_{n}\|_{1,\infty}
\endmultline \tag 6.17
$$
The first term is of the form in the right of (6.14). The terms in the sum becomes
a sum of four terms of the form $K_1\|h_1\|_{2,\infty}\|h_2\|_{r+1,\infty}$
multiplied by factors of the form $\|f_k\|_{1,\infty}$.
Using (6.5) with $j=1$ we can bound
$\|h_1\|_{2,\infty}\|h_2\|_{r+1,\infty}\leq
C\|h_1\|_{1,\infty}\|h_2\|_{r+2,\infty}+\|h_1\|_{r+1,\infty}\|h_2\|_{1,\infty}$.
This proves also (6.14).
\qed
\enddemo

\proclaim{Lemma {6.}4} Let $p$ be the solution of $\triangle p=-(\pa_i V^j)(\pa_j V^i)$,
where $v^i=D_t \, x^i$ and let $\dot{p}=D_t p$.
 Then for $r\geq 1$ we have
$$\align
\|p\|_{r,\infty}&\leq K_3(\|\dot{x}\|_{r,\infty}+\|{x}\|_{r+1,\infty})
\tag 6.18\\
\|\dot{p}\|_{r,\infty}&\leq K_3(\|\ddot{x}\|_{r,\infty}+\|\dot{x}\|_{r+1,\infty}
+\|{x}\|_{r+2,\infty}) \tag 6.19
\endalign
$$
\endproclaim
\demo{Proof} We just apply Proposition {7.}1 to
$$
 \triangle p=-(\pa_i V^j)\pa_j V^i,\quad\quad  v^i=D_t x^i,\qquad
p\big|_{\pa\Omega}=0  \tag 6.20
$$
using Lemma {6.}3 to estimate the product.
(Recall that $\|V\|_{2,\infty}\leq K_3$.)
Since $D_t=\pa_t+V^k\pa_k$, where $\pa_t=\pa_t\big|_{x=const}$, we have
$$
 \triangle \dot{p}=\triangle\big( (\pa_t+V^k\pa_k) p\big)
=D_t\triangle p+ (\triangle V^k)\pa_k p
+2\delta^{ij}(\pa_i V^k)\pa_j\pa_k p
\tag 6.21
$$
and
$$
D_t\, \triangle p=-(\pa_t+V^k\pa_k )\big( (\pa_i V^j)(\pa_j V^i)\big)
=-2(\pa_i V^j)(\pa_j \dot{V}^i)+2(\pa_i V^j)(\pa_j V^k)\pa_k V^i\tag 6.22
$$
so
$$
 \triangle \dot{p}=-2(\pa_i V^j)(\pa_j \dot{V}^i)+2(\pa_i V^j)(\pa_j V^k)\pa_k V^i
+ (\triangle V^k)\pa_k p
+2\delta^{ij}(\pa_i V^k)\pa_j\pa_k p\tag 6.23
$$
The second part of the lemma now follows from Proposition {7.}1 using
Lemma {6.}3 and the first part of the lemma.
 \qed
\comment
$$
 \triangle \dot{p}=-\kappa^{-1}\pa_b\big(( D_t (\kappa g^{ab}))\pa_b p\big)
-\hat{D}_t\big((\pa_i V^k)\pa_j V^i\big),
\tag 6.24
$$
\endcomment
\enddemo

Let us now introduce the $L^\infty$ norms that we will use:
\demo{Definition {6.}2}
$$
\align
m_s(t)&=\|x(t,\cdot)\|_{1+s,\infty}\tag 6.24\\
\dot{m}_s(t)&= \|x(t,\cdot)\|_{2+s,\infty}
+\|\dot{x}(t,\cdot)\|_{1+s,\infty},\tag 6.25\\
 \ddot{m}_s(t) &=  \|x(t,\cdot)\|_{3+s,\infty}
+\|\dot{x}(t,\cdot)\|_{2+s,\infty}
 +\|\ddot{x}(t,\cdot)\|_{1+s,\infty}, \tag 6.26\\
 {n}_s(t) &=  \|x(t,\cdot)\|_{4+s,\infty}
+\|\dot{x}(t,\cdot)\|_{3+s,\infty}
 +\|\ddot{x}(t,\cdot)\|_{2+s,\infty}, \tag 6.27
\endalign
$$
\enddemo

We remark that in Definition {5.}2 we made an assumption that the
inverse of $g$ and $\pa y/\pa x$ are bounded. This means that
$m_0$ etc are all bounded from below as well. We note that the
corresponding bounds for the metric $g_{ab}=\delta_{ij}(\pa
x^i/\pa y^a)(\pa x^j/\pa y^b)$ and $\omega_{ab}=(\curl v)_{ab}
=(\pa x^i/\pa y^a)(\pa x^j/\pa y^b) (\pa_i v_j-\pa_j v_i)$ follows
from the bounds for $x$, $\dot{x}$, and $\ddot{x}$:
$$
\|g\|_{r,\infty}\leq K_1 m_r,\qquad
\|\dot{g}\|_{r,\infty}+\|\omega\|_{r,\infty}\leq K_2 \dot{m}_r,
\qquad \|\ddot{g}\|_{r,\infty}+\|\dot{\omega}\|_{r,\infty}\leq K_3 \ddot{m}_r.
\tag 6.28
$$
The proof of this uses the interpolation inequality (6.5) in Lemma
{6.}2 applied to each term we get when we differentiate.
In view of the coordinate condition, see
Definition {5.}1, the same bounds also hold for $g$ replaced by
the inverse of $g$.

Furthermore, we will now prove that the corresponding bounds also
for the pressure follows from this. We will actually loose a
derivative when passing to the bounds for the pressure because we
will go over H\"older spaces, but this does not matter. In
Lemma {6.}4, we proved that
$$
\align
\|p(t,\cdot)\|_{r+1,\infty}&\leq K_3  \dot{m}_r(t),\tag 6.29 \\
 \|p(t,\cdot)\|_{r+2,\infty}+ \|\dot{p}(t,\cdot)\|_{r+1,\infty}&\leq
K_3 \ddot{m}_r(t),\tag 6.30 \\
 \|p(t,\cdot)\|_{r+3,\infty}+ \|\dot{p}(t,\cdot)\|_{r+2,\infty}&
\leq K_3 n_r(t)\tag 6.31
\endalign
$$
In particular
$$
\|p\|_{2,\infty}+\|\dot{p}\|_{1,\infty}\leq K_3,
\qquad
\|\dot{p}\|_{2,\infty}\leq K_4\tag 6.32
$$
We will frequently use interpolation, e.g.
$$
m_r\dot{m}_s\leq C( m_{r+s}\dot{m}_0+m_0\dot{m}_{r+s}) \leq K_2
\dot{m}_{r+s},\tag 6.33
$$
which follows from Lemma {6.}1 and the proof of Lemma {6.}2
applied to each term we get when multiplying any of the
expressions (6.24)-(6.27) together.

We must also ensure that if the
physical condition (2.7) and coordinate condition (2.8)
holds initially they will hold for some
small time $0\leq t\leq T$, with $c_0$ replaced by $c_0/2$
and $c_1$ replaced by $2c_1$. This will be proven in section 11,
and until then we will just assume that $T$ is so small that4
these conditions hold for
$0\leq t\leq T$. Furthermore, we will also assume that $T\leq c_0\leq 1$
since the estimates we will derive then will be independent of $T$ and $c_0$.

\head 7. The $L^\infty$ estimates for the Dirichlet problem.\endhead

In this section, we give tame H\"older estimates for the
solution of the Dirichlet problem:
$$
\triangle  q=F,\qquad\qquad q\big|_{\pa\Omega}=0.\tag 7.1
$$
Our H\"older estimates loose a derivative since we want to use them
for integer values.
This is not important and with an additional loss of regularity,
we could have avoided using H\"older estimates altogether and just
gotten the $C^k$ estimates from the Sobolev estimates, proved
in the next section, using Sobolev's lemma.
Apart from getting estimates for the solution of (7.1) we also need estimates
for time derivatives and variational derivatives. For this we need to know that
the solution of (7.1) depend smoothly on parameters if the metric and the
inhomogeneous term do.
We remark that the coordinate condition is critical
since it is needed in order to invert the Laplacian.

One can also use the results in section 5
 to get tame estimates for the solution of the
Dirichlet problem:
In fact if we take $W^a=g^{ab}\pa_b q$, and $w_a=\pa_q q$, then
$ \div W= \triangle q$ and
$\curl w=0$. Applying Lemma {5.}2
to $W$ therefore gives:
$$
|W|_{r}^{\Cal R}\leq K_1\sum_{s=0}^{r} \,\,[g]_{s}^{\Cal R}
\big(|\triangle q |_{r-1-s}^{\Cal R}
+|W|_{r-s}^{\Cal S}\big),\tag 7.2
$$
where for $s=r$ we should interpret $|\triangle q|_{-1}=0$,
and
$$
|\pa q|_{r}^{\Cal R}\leq K_1\sum_{s=0}^{r}\,\, [g]_{s}^{\Cal R}
\big(|\triangle q|_{r-1-s}^{\Cal R}
+|\pa q|_{r-s}^{\Cal S}\big)\tag 7.3
$$
%Here $|\pa q|_r^{\Cal V}=\sum_{|I|\leq r,\, I\in{\Cal V}} |\pa U^I q|$ and
%$[g]^{\Cal R}_s=\sum_{s_1+...s_1\leq s,\, s_i\geq 1} |g|_{s_1}\cdot\cdot\cdot|g|_{s_k}$.
Therefore it suffices to obtain estimates for tangential derivatives only
which is easier because the Dirichlet boundary condition is preserved by
tangential derivatives. If $q\big|_{\pa\Omega}=0$ then
$S^I q\big|_{\pa\Omega}=0$.
 The $L^\infty$ estimates uses the standard
Schauder estimates for the Dirichlet problem. Because we want to
have the final result in terms of $C^k$ norms instead of H\"older
norms these results loose a derivative.

\proclaim{Proposition {7.}1} If $q\big|_{\pa\Omega}=0$ then for $r\geq 1$
$$
\leftalignspace\|q\|_{r,\infty}\leq K_3 \big(\|\triangle
q\|_{r-1,\infty} +\|g\|_{r,\infty} \|\triangle q\|_{0,\infty}\big)\tag
7.4
$$
\endproclaim
\demo{Proof}
If we apply Lie derivatives $\hat{\Cal L}_S^I$
 to $W^a=g^{ab}\pa_b q$ we get
$$
W_I^a=g^{ab}\pa_b S^I q + \tilde{c}^{\,I}_{I_1 I_2} \hat{g}^{I_1 ab}\,\pa_b S^{I_2} q,
\qquad\quad \hat{g}^{I \,ab}=\hat{\Cal L}_S^I g^{ab}, \qquad
W_I=\hat{\Cal L}_S^I W\tag 7.5
$$
and the sum is over all combinations $I=I_1+I_2$,
$\tilde{c}^{\,I}_{I_1 I_2}$ are constants such that $\tilde{c}^{\,I}_{I_1 I_2}=0$
if $I_2=I$. Since $\div W_I=\div \hat{\Cal L}_S^I W=\hat{S}^I \div W
=\hat{S}^I \triangle q=\kappa^{-1}S^I(\kappa \triangle q)$ it follows from
taking the divergence of (7.5) that
$$
\triangle( S^I q)=\hat{S}^I \triangle q
 -\kappa^{-1}\pa_a\big( \tilde{c}^{\,I}_{I_1 I_2} \hat{g}^{I_1 ab}\,\pa_b S^{I_2} q\big),
\qquad\quad \hat{g}^{I \,ab}=\hat{\Cal L}_U^I g^{ab}\tag 7.6
$$
Let $\|u\|_{2+\alpha,\infty}$ denote H\"older norms, see section 17 and
Proposition {7.}2. By Proposition {7.}2 we have
$$
\|S^I q\|_{2+\alpha,\infty}\leq K_1  \big( \|\hat{S}^I \triangle
q\|_{1,\infty} +\tilde{c}^{\,I}_{I_1 I_2}( \|\hat{g}^{I_1}\|_{1,\infty}
\|S^{I_2} q\|_{2+\alpha,\infty} +\|\hat{g}^{I_1}\|_{2,\infty} \|S^{I_2}
q\|_{1+\alpha,\infty}) \big) \tag 7.7
$$
If we let $M_r=\sum_{|I|\leq r-2} \|S^I q\|_{2+\alpha,\infty}$, $r\geq
2$, $M_r=\|q\|_{r+\alpha,\infty}$ for $r=0,1$  it follows from Proposition
{7.}2 that $M_0+M_1\leq K_3 \|\triangle q\|$, $M_2\leq
K_3\|\triangle q\|_{1,\infty}$ and for $r\geq 3$ we have:
$$
M_r\leq K_3\big( \|\triangle q\|_{r-1,\infty}+\sum_{s=1}^{r-1}
\|{g}\|_{r+1-s,\infty} M_{s} \big)\tag 7.8
$$
Inductively it follows that
$$
M_r\leq K_3\big( \|\triangle q\|_{r-1,\infty}+\sum_{s=0}^{r-2}
\|{g}\|_{r-s,\infty} \|\triangle q\|_{s,\infty}\big) \leq K_3
\big(\|\triangle q\|_{r-1,\infty}+\|g\|_{r} \|\triangle
q\|_{0,\infty}\big) \tag 7.9
$$
where we used interpolation. With $I\in {\Cal S}$,
$|I|=r-2$ we have hence get from differentiating (7. ) and using what we used proved
$$
\|\pa W_I\|_{0,\infty}\leq K_3  \big(\|\triangle q\|_{r-1,\infty}
+\|g\|_{r} \|\triangle q\|_{0,\infty}\big) \tag 7.10
$$
However, once we have bound for the tangential components, the bound for all components
in terms of these and $\hat{R}^I\triangle p$.
follows from Lemma {5.}2.
\enddemo

Theorem 6.6 in \cite{GT}, together with Theorem 8.16,
and Theorem 8.33 in \cite{GT} in our setting it reads:
\proclaim{Proposition {7.}2} Suppose that $\|\phi\|_{k+\alpha,\infty}$ denotes the H\"older
norms and $0<\alpha<1$, and $k$ is an integer, see section 17.
$$
\triangle p=g^{ab}\pa_a\pa_b p +\kappa^{-1}(\pa_a (\kappa g^{ab}))\pa_b p
=\kappa^{-1}\pa_a \big(\kappa g^{ab}\pa_b p\big)\tag 7.11
$$
where
$$
\|g^{ab}\| _{0+\alpha,\infty}+\|\pa g^{ab}\| _{0+\alpha,\infty}\leq \Lambda,
\qquad \sum_{a,b}|g^{ab}|+|g_{ab}|\leq \lambda\tag 7.12
$$
Suppose that $p\big|_{\pa\Omega}=0$. Then
$$
\|p\|_{2+\alpha,\infty}\leq C\big(\|p\|_\infty+\|\triangle p\|_{0+\alpha,\infty}\big)
\tag 7.13
$$
where $C=C(n,\alpha,\lambda,\Lambda)$
and
$$
\|p\|_\infty\leq C \|\triangle p\|_\infty \tag 7.14
$$
and if $\triangle p=F+\kappa^{-1}\pa_a(\kappa G^a )$
$$
\|p\|_{1+\alpha,\infty}
\leq C\big(\|p\|_\infty + \|F\|_\infty +\|G\|_{0+\alpha,\infty}\big)
\tag 7.15
$$
Furthermore
$$
\|uv\|_{\alpha,\infty}\leq C\|u\|_{\gamma,\infty}
\|v\|_{\alpha,\infty},\qquad \gamma\geq \alpha,
\qquad
\| u v\|_{\alpha,\infty}\leq
 C\big(\|u\|_{0,\infty} \|v\|_{\alpha,\infty} +\|v\|_{0,\infty}
 \|u\|_{\alpha,\infty}\big)
\tag 7.16
$$
\endproclaim
Note that if we multiply by $\kappa$ then the
operator is also in the divergence form that \cite{GT} has in Theorem 10.33.
Anyway, in our case it is equivalent to a domain in ${\Cal D}_t$ with the
the usual metric.

Let us now prove that the solution of (7.1) depends smoothly on
parameters if the metric $g$ and the inhomogeneous term $F$ do.
Let us assume that the parameter is time $t$.
We have:
\proclaim{Lemma {7.}3} Let $\phi$ be the solution of
$$
\triangle \phi= \kappa^{-1}\pa_a \big( \kappa  g^{ab}\pa_b
\phi\big)=F, \qquad \phi\big|_{\pa\Omega}=0\tag 7.17
$$
where $\kappa=\sqrt{\det{g}}$, and $g$ satisfies the coordinate condition
(2.8) on $[0,T]$. Suppose that $g_{ab},F\in C^k([0,T],C^\infty(\overline{\Omega}))$.
Then $\phi\in  C^k([0,T],C^\infty(\overline{\Omega}))$.
\endproclaim
\demo{Proof}
Let us write $\phi_t$, $g_t$ $F_t$, and $\triangle_t=\triangle_{g_t}$,
to indicate the dependence of $t$. Our initial assumption is that
$g_t,F_t \in C^k([0,T],C^\infty(\overline{\Omega}))$. That
$$
\triangle_t \phi_t =F_t,\qquad \phi_t\big|_{\pa{\Omega}}=0\tag 7.18
$$
has a solution $\phi_t\in C^\infty(\overline{\Omega})$ if
$F_t,g_t\in C^\infty(\overline{\Omega})$ and the coordinate
condition is fulfilled is well known. We will prove that
$F_t,g_t\in C^1([0,T],C^\infty(\overline{\Omega}))$ implies that
$\phi_t \in C^1([0,T],C^\infty(\overline{\Omega}))$. If this is
the case, then $\dot{\phi}_t=D_t\phi_t$ satisfies
$$
\triangle_t \dot{\phi}_t= \dot{F}_t-\dot{\triangle}_t \phi_t,
\qquad\quad  \dot{\phi}_t\big|_{\pa{\Omega}}=0\tag 7.19
$$
where $\dot{\triangle}_t=[D_t, \triangle ]$ and $\dot{F}_t=D_t F_t$.
Since the right hand side of (7.19) is also in $ C^1([0,T],C^\infty(\overline{\Omega}))$
we can repeat argument to conclude that
 $\dot{\phi}_t\in C^1([0,T],C^\infty(\overline{\Omega}))$,
 i.e ${\phi}_t\in C^2([0,T],C^\infty(\overline{\Omega}))$
In general we can then use induction to conclude that $\phi_t\in
C^k([0,T],C^\infty(\overline{\Omega}))$ since
$$
\triangle_t D_t^k \phi_t=D_t^k F
-\sum_{j=0}^{k-1}c^k_{l}\triangle_t^{(k-j)} D_t^j \phi_t, \tag
7.20
$$
where $\triangle_t^{(k)}$ are the repeated commutators defined
inductively by $\triangle_t^{(k)}=[D_t,\triangle^{(k-1)}_t]$,
$\triangle_t^{(0)}= \triangle_t$.

First we will show that $F_t,g_t\in
C([0,T],C^\infty(\overline{\Omega}))$ implies that $\phi_t \in
C([0,T],C^\infty(\overline{\Omega}))$. We will only prove this for
$t=0$ since the proof in general is just a notational difference
from the proof for $t=0$.
$$
\triangle_t
(\phi_t-\phi_0)=F_t-F_0-(\triangle_t-\triangle_0)\phi_0\tag 7.21
$$
Since the $C^m(\overline{\Omega})$ or $H^m(\overline{\Omega})$
norm of the right hand side tends to $0$ as $t\to 0$ for any $m$
and since we have uniform bounds for $\triangle_t^{-1}$, in Lemma
{7.}3, it follows that  the $C^m(\overline{\Omega})$ or
$H^m(\overline{\Omega})$ norm of $\phi_t-\phi_0$ tends to $0$ as
$t\to 0$ for any $m$. Hence $\phi_t\in
C([0,T],C^\infty(\overline{\Omega}))$. Now, let $\dot{\phi}_t$ be
defined by (7.19). By the same argument it follows that
$\dot{\phi}_t\in C([0,T],C^\infty(\overline{\Omega}))$. It remains
to prove that $\phi_t$ is differentiable. We have
$$
\triangle_t \big(\phi_t-\phi_0-t\dot{\phi}_0\big)
=F_t-F_0-t \dot{F}_0
+\big(\triangle_t-\triangle_0-t\dot{\triangle}_0\big)\phi_0
+t(\triangle_t-\triangle_0)\dot{\phi}_0\tag 7.22
$$
Since $F_t$ and $g_t$ are differentiable as functions of $t$ it
follows that the $C^m(\overline{\Omega})$ or
$H^m(\overline{\Omega})$ norm of the right hand side divided by
$t$ tends to $0$ as $t\to 0$ for any $m$. Since we also have
bounds for the inverse of $\triangle_t$ that are uniform in $t$ we
conclude that any $C^m$ or $H^m$ norm of
$\phi_t-\phi_0-t\dot{\phi}_0$ divided by $t$ also tends to $0$ as
$t\to 0$ for any $m$. It follows that $\phi_t\in
C^1([0,T],C^\infty(\overline{\Omega}))$. \qed\enddemo

\head 8. The $L^2$ estimates for the Dirichlet problem. \endhead
In this section, we give tame $L^2$-Sobolev estimates for the
solution of the Dirichlet problem:
$$
\triangle  q=F,\qquad\qquad q\big|_{\pa\Omega}=0.\tag 8.1
$$

We also remark that the coordinate condition is critical in all the estimates
in this section since it is needed in order to invert the Laplacian $\triangle$.
As pointed out in the beginning of section 7, using the results from section 5
it suffices to obtain estimates for tangential derivatives only
which is easier because the Dirichlet boundary condition is preserved by
tangential derivatives. If $q\big|_{\pa\Omega}=0$ then
$S^I q\big|_{\pa\Omega}=0$.

\proclaim{Proposition {8.}1} Suppose that $q$ is a solution of the
Dirichlet problem, $q\big|_{\pa\Omega}=0$, and $W^a=g^{ab}\pa_b
q$. Then if $r\geq 0$ we have
$$
\|W\|_{r}
\leq K_1\sum_{s=0}^{r-1} \|g\|_{r-1-s,\infty}\|\triangle q\|_s
+K_1\|g\|_{r,\infty} \|W\|,
\tag 8.2
$$
and if $r\geq 1$ we have
$$
\|W\|_{r}+\|q\|_{r+1}
\leq K_1\sum_{s=0}^{r-1} \|g\|_{r-s,\infty}\|\triangle q\|_s
\tag 8.3
$$
Furthermore, for $i\leq 2$ and $r\geq 0$ we have
$$
\|\hat{D}_t^i {W}\|_{r}\leq K_3\sum_{s=0}^{r-1} \sum_{j+k \leq i}
 \| \check{D}_t^k g\|_{r-s,\infty}
\|\hat{D}_t^j\triangle q\|_s+ K_3\sum_{j+k\leq i}  \| \check{D}_t^{k} g\|_{r,\infty}
\|\hat{D}_t^j W\|\tag 8.4
$$
and if $r\geq 1$ we have
$$
\|\hat{D}_t^i {W}\|_{r}+\|D_t^i q\|_{r+1}\leq K_3\sum_{s=0}^{r-1}\sum_{j+k\leq i}
 \| \check{D}_t^k g\|_{r-s,\infty}
\|\hat{D}_t^j\triangle q\|_s\tag 8.5
$$
Moreover if $P$ is the orthogonal projection onto divergence free
vector fields and $W$ is any vector field then, for $r\geq 0$,
$$
\|PW\|_r\leq K_1 \sum_{s=0}^{r} \|g\|_{r-s,\infty}\|W\|_s\tag 8.6
$$
and, for $i\leq 2$,
$$
\|\hat{D}_t^i  PW\|_r\leq K_3 \sum_{s=0}^{r}
\sum_{j+k\leq i} \|\check{D}_t^k g\|_{r-s,\infty}\|\hat{D}_t^{j} W\|_s
\tag 8.7
$$
\endproclaim

First a useful lemma:
\proclaim{Lemma {8.}2} Suppose that $S\in {\Cal S}$ and $q\big|_{\pa\Omega}=0$,
and
$$
\hat{\Cal L}_S W^a=g^{ab}\pa_b q+F^a .\tag 8.8
$$
Then
$$
\| \hat{\Cal L}_S W\|\leq K_1(\| \div W\|+\|F\|)\tag 8.9
$$
\endproclaim
\demo{Proof} Let $W_S=\hat{\Cal L}_S W$.
$$
\int_\Omega g_{ab}  W_S^a W_S^b \, \kappa \, dy
=\int_\Omega W_S^a \pa_a  q\, \,\kappa dy
+\int_\Omega W_S^a  g_{ab} F^b\kappa dy \tag 8.10
$$
If we integrate by parts in the first integral on the right,
using that $q$ vanishes on the boundary we get
$$
\int_\Omega W_S^a\, \pa_a q\, \kappa\, dy
=-\int_\Omega \div\,( W_S)\, q \,\kappa \, dy\tag 8.11
$$
However $\div W_S=\hat{S} \div W$.  Then we
can integrate by parts in the angular direction. $S=S^a\pa_a$,
$\hat{S}=S+\div S$ so we get
$\int_\Omega ( \hat{S} f)\, \kappa dy=
\int_\Omega \pa_a \big( S^a f\kappa )\, dy=0$, where $\pa_a S^a=0$.
 Hence we get
$$
\int_\Omega W_S^a\, \pa_a q\, \kappa\, dy
=\int_\Omega \div\,( W)\, ({S} q )\kappa \, dy\tag 8.12
$$
Here $| {S}  q|\leq K_1 |\pa q|$ so it follows that
$$
\|W_S\|^2\leq K_1 \|\div W\|\big( \|W_S\| +\|F\| \big)
+K_1 \|W_S\| \| F\|\tag 8.13
$$
and so
$$
\|W_S\|\leq K_1 (\| \div W\|+\|F\|)\qed \tag 8.14
$$
\enddemo

\demo{Proof of Proposition {8.}1} If we apply ${\Cal L}_S^I$ to $w_a=g_{ab}W^b$ we get
$$
\pa_a S^I q=g_{ab}W_I^b+ \tilde{c}^{I_1 I_2} \check{g}_{I_2\,ab} W_{I_2}^b,
\qquad \quad
W_I=\hat{\Cal L}_S^I W, \qquad  \check{g}_{I\,ab}=\check{\Cal L}_S^{I} g_{ab}
\tag 8.15
$$
and the sum is over all combinations $I=I_1+I_2$,
$\tilde{c}^{I_1 I_2}$ are constants such that $\tilde{c}_I^{\,I_1 I_2}=0$
if $I_2=I$.
If we write $S^I=S S^J$, $W_I= \hat{\Cal L}_S W_J$,
and use Lemma {8.}2 we get since $\div W_J =\hat{S}^J \div W=
\hat{S}^J \triangle q=\kappa^{-1} S^J (\kappa \triangle q)$
$$
\|W_I\|\leq K_1\|\hat{S}^J \triangle q\|
+ K_1\tilde{c}_I^{\,I_1 I_2}\|\check{g}_{I_1}\|_\infty \|W_{I_2}\|\tag 8.16
$$
or if we sum over all of them and use interpolation
$$
\|\kappa W\|^{\Cal S}_{r}\leq K_1\|\kappa\triangle q\|_{r-1}
+K_1 \sum_{s=0}^{r-1}\|g\|_{r-s,\infty} \|\kappa W\|_s^{\Cal S},
\qquad r\geq 1\tag 8.17
$$

We now want to apply Lemma {5.}3 to $W^a=g^{ab}\pa_b q$. Then
$\curl w=0$ and $\div W=\triangle q$ so
$$
\| W\|_{r}\leq K_1\|\kappa\triangle q\|_{r-1}+
K_1\|\kappa  W\|_{r}^{\Cal S}+K_1\sum_{s=0}^{r-1}
\|g\|_{r-s,\infty} \| W\|_{s},\qquad r\geq 1\tag 8.18
$$
We now use (8.17) to replace the term  $K_1\| \kappa W\|_{r}^{\Cal S}$
by the terms of the form already in the right hand side of (8.18) and
we also replace $\kappa$ by $1$ which just produces more terms of the same form.
Using induction in $r$ and interpolation (8.2) follows.

We also need an estimate for the lowest order term:
$\|W\|^2=\int_\Omega g^{ab}(\pa_a q)(\pa_b q)\, \kappa \, dy
=-\int_{\Omega}( \triangle q) q\, \kappa \, dy$.
However, there is a constant depending just on the volume of $\Omega$,
i.e. $\int_{\Omega} \kappa \, dy$,
such that $\|q\|\leq K_1\|\triangle q\|$, see \cite{SY}.
Therefore in addition we have
$$
\|W\| +\|q\|+\|\pa q\|\leq K_1\|\triangle q\|,\qquad
\text{if}\quad W^a =g^{ab}\pa_b q. \tag 8.19
$$
Therefore we inductively get from also using interpolation:
$$
\|\kappa W\|_{r}\leq
K_1\sum_{s=0}^{r-1} \|g\|_{r-s,\infty}\|\triangle q\|_{s},\qquad r\geq 1
\tag 8.20
$$
We note that we can remove $\kappa$ in the left since doing so
just produces lower order terms of the same form.
This proves the estimate for $W$ in (8.3) and the estimate for $q$
follows from this since $W^a=g^{ab}\pa_b q$.
In fact by (8.15) with $S^I$ replaced by any space vector fields
$R^I$, $I\in{\Cal R}$, $\|\pa q\|_{r}\leq
K_1\sum_{s=0}^{r} \|g\|_{r-s,\infty}\|W\|_{s}$
and by (8.19) we also have an estimate for $\|q\|$.

It remains to prove the
estimate with time derivatives.
We can now repeat the argument with $i$ of the tangential derivatives
being the time derivative, $\hat{\Cal L}_{D_t}^i=\hat{D}_t^i$
and $\check{\Cal L}_{D_t}^i=\check{D}_t^i$. This gives
$$
\pa_a S^I D_t^i\, q=g_{ab}\hat{D}_t^i W_I^b+ \tilde{c}_{I i}^{\,I_1 i_1 I_2 i_2}
(\check{D}_t^{i_1}\check{g}_{I_2\,ab}) \hat{D}_t^{i_2} W_{I_2}^b,\tag 8.21
$$
where $\tilde{c}_I^{\,I_1 i_1  I_2 i_2 }=1$ for all $(I_1+I_2,i_1+i_2)=(I,i)$
such that $|I_2|+i_2<|I|+i$ and $0$ otherwise.
By Lemma {8.}2 again
$$
\|\hat{D}_t^i {W}_I\|\leq C\|\hat{S}^J\hat{D}_t^i \triangle q\|
+ C\tilde{c}_{I i}^{\,I_1 i_1 I_2 i_2}
\|\check{D}_t^{i_1}\check{g}_{I_1}\|_\infty \|\hat{D}_t^{i_2} W_{I_2}\|
\tag 8.22
$$
where $|J|=|I|-1$. Hence
$$
\|D_t^i (\kappa {W})\|_{r}^{\Cal S}\leq K_1\|{D}_t^i (\kappa\triangle q)\|_{r-1}
+K_1 \!\!\!\!\!\!\sum_{s\leq r,\, j\leq i,\, s+j<r+i} \!\!\!\!\!\!
\|D_t^{i-j}(\kappa^{-1}g)\|_{r-s,\infty}
\|D_t^j(\kappa {W})\|_s^{\Cal S}
\tag 8.23
$$
By Lemma {8.}3 below and (8.23) it follows that for $i\leq 2$
$$
\| D_t^i(\kappa {W})\|_{r}\leq K_3\Big(\|{D}_t^i(\kappa\triangle q)\|_{r-1}
+\!\!\!\!\!\!\sum_{s\leq r,\, j\leq i,\, s+j<r+i} \!\!\!\!\!\!
\|D_t^{i-j}(\kappa^{-1}g)\|_{r-s,\infty}
\|D_t^j(\kappa {W})\|_s\Big). \tag 8.24
$$
(8.4) now follows from this and interpolation.

Since $\hat{D}_t \triangle  q=\triangle D_t \, q
+ \kappa^{-1}\pa_a \big( \kappa(  \hat{D}_t g^{ab}) \pa_b q\big)$
we get
$$
\|\triangle D_t q\|\leq K_2 \big(\|\hat{D}_t\triangle q\|
+\|\check{D}_t g\|_{0,\infty}\|\pa^2 q\|+\|\check{D}_t g\|_{1,\infty}
\|\pa q\|\big)\leq K_3\big(\|\hat{D}_t\triangle q\|+\|\triangle q\|\big)
,\tag 8.25
$$
where we used (8.3) with $r=2$.
Therefore using (8.19) applied to $D_t q$ in place of $q$
we also get an estimate for the lowest order norm:
$$
\|\hat{D}_t W\|+\|\pa D_t q\|+\|D_t q\|\leq
K_3\big(\|\hat{D}_t\triangle q\|+\|\triangle q\|\big)\tag 8.26
$$
Using this, (8.4) and interpolation
gives (8.5) for one time derivative, apart from the estimate for
$\| D_t q\|_{r+1}$. By (8.21) with $S^I$ replaced by any space vector fields
$R^I$, $I\in{\Cal R}$, $\|\pa D_t \,q\|_{r}\leq
K_1\big(\sum_{s=0}^{r} \|g\|_{r-s,\infty}\|\dot{W}\|_{s}
+\sum_{s=0}^{r} \|\check{D}_t g\|_{r-s,\infty}\|{W}\|_{s}\big)$
and by (8.26) we also have an estimate for $\| D_t q\|$.

Since $\hat{D}_t^2 \triangle  q=\triangle D_t^2 \, q
+2 \kappa^{-1}\pa_a \big( \kappa(  \hat{D}_t g^{ab}) \pa_b D_t q\big)
+ \kappa^{-1}\pa_a \big( \kappa(  \hat{D}_t^2 g^{ab}) \pa_b q\big)$
we get
$$\multline
\|\triangle D_t^2 q\|\leq K_2 \big(\|\hat{D}_t^2\triangle q\|
+\|\check{D}_t^2 g\|_{0,\infty}\|\pa^2 q\|+\|\check{D}_t^2 g\|_{1,\infty}
\|\pa q\|+\|\check{D}_t g\|_{0,\infty}\|\pa^2 D_t\, q\|+
\|\check{D}_t g\|_{1,\infty}
\|\pa D_t\,q\|\big)\\
\leq K_3\big(\|\hat{D}_t^2\triangle q\|+\|\hat{D}_t\triangle q\|+\|\triangle q\|
\big)
\endmultline \tag 8.27
$$
where we used (8.5) for $i\leq 1$.
Therefore we also get an estimate for the lowest order norm:
$$
\|\hat{D}_t^2 W\|+\|\pa D_t^2 q\|+\|D_t^2 q\|\leq
K_3\big(\|\hat{D}_t^2\triangle q\|+\|\hat{D}_t\triangle q\|+\|\triangle q\|
\big)\tag 8.28
$$
Using this, (8.4) and interpolation
gives (8.5) also for two time derivative, apart from the estimate for
$\| D_t^2 q\|_{r+1}$. By (8.21) with $S^I$ replaced by any space vector fields
$R^I$, $I\in{\Cal R}$, $\|\pa D_t^2 \,q\|_{r}\leq
K_3 \sum_{s=0}^{r} \big(\|g\|_{r-s,\infty}\|\ddot{W}\|_{s}
+\|\check{D}_t g\|_{r-s,\infty}\|\dot{W}\|_{s}+
\|\check{D}_t^2 g\|_{r-s,\infty}\|{W}\|_{s}\big)$
and by (8.28) we also have an estimate for $\| D_t^2 q\|$.

It remains to prove the estimates for the projection (8.6)-(8.7). We
have $W=W_0+W_1$, where $W_0=PW$, and $W_1=g^{ab}\pa_b q$ where
$\triangle q=\div W$ and $q\big|_{\pa\Omega}=0$. Proving
(8.6)-(8.7) for $r\geq 1$ reduces to proving it for $r=0$
by using (8.6)-(8.7), since $\hat{R}^I \triangle
q=\div\,( \hat{\Cal L}_R^I W)$ and replacing $\kappa$ by $1$ just produces
more terms of the same form . (8.6)) for $r=0$ follows since
the projection has norm $1$, $\|PW\|\leq \|W\|$.  Since the
projection of $g^{ab}D_t^k w_{1b}=g^{ab}\pa_b D_t^k q$ vanishes we obtain from Lemma
{8.}3 below:
$$
\| P \hat{D}_t^i W_1\| \leq K_1 \sum_{j=0}^{i-1}
\|\check{D}_t^{i-j} g\|\,
\|\hat{D}_t^j W_1\|\tag 8.29
$$
Since also $P \hat{D}_t^i W_0=\hat{D}_t^i W_0$ we have
$$
(I-P)\hat{D}_t^i W_1=(I-P) D_t^i W\tag 8.30
$$
and hence since the projection has norm one
$$
\| \hat{D}_t^i W_1\|+\| \hat{D}_t^i W_0\|\leq K_1 \|\hat{D}_t^i W\|+ K_1
\sum_{j=0}^{i-1} \|\check{D}_t^{i-j} g\|\, \|\hat{D}_t^j W_1\|\tag 8.31
$$
Hence for $i=0,1,2$ it inductively follows that
$$
\|\hat{D}_t^i W_0\|+\|\hat{D}_t^i W_1\|\leq K_3 \sum_{j=0}^i
\|\hat{D}_t^{j} W\|,\qquad i\leq 2. \tag 8.32
$$
Since as before replacing $\kappa$ by $1$ just produces
more terms of the same form
this proves (8.7) also for $r=0$.
(6.4)-(6.5) follows from interpolation.
\qed\enddemo

\proclaim{Lemma {8.}3} Let $W^a=g^{ab} w_b $ Then
$$
\hat{D}_t^i W^a =g^{ab} D_t^i w_b-\sum_{j=0}^{i-1} \binom{i}{j}
g^{ab}(\check{D}_t^{i-j} g_{bc})\hat{D}_t^j W^c\tag 8.33
$$
Furthermore if $W^a=g^{ab} \pa_b q$ then
$$
|D_t^i(\kappa W)|_r^{\Cal R}
\leq K_1\Big( |D_t^i(\kappa\div W) |_{r-1}^{\Cal R}
+|\kappa W|_r^{\Cal S}+\sum_{s\leq r,\, j\leq i, s+j<r+i}
|D_t^{i-j} (\kappa g) |_{r-s}^{\Cal R}
 | D_t^j(\kappa W)|_{s}^{\Cal R}\Big)\tag 8.34
$$
\endproclaim
\demo{Proof} We have $D_t^i w_b=D_t^i\big( \kappa^{-1} g_{bc} \,
\kappa W^c\big)=\sum_{j=0}^i\binom{i}{j} \big( D_t^{i-j}
(\kappa^{-1} g_{bc})\big) \hat{D}_t^{j} (\kappa
W^{c})$, which proves (8.33).

(8.34) follows from (5.12) using (8.33) and that fact the
curl of $w_a=\pa_a$ vanishes to estimate
the curl of $g_{ab} \hat{D}_t^i W^b$.
 \qed\enddemo

\comment
Similarly, for two time derivatives we get
$$
\pa_a S^I D_t^2\, q=g_{ab}\ddot{W}_I^b+ 2{c}^{I_1 I_2}
(\check{D}_t\check{g}_{I_2\,ab}) \dot{W}_{I_2}^b
+ {c}_I^{\,I_1 I_2}
(\check{D}_t^2\check{g}_{I_2\,ab}) {W}_{I_2}^b
+ \tilde{c}_I^{\,I_1 I_2} \check{g}_{I_2\,ab} \ddot{W}_{I_2}^b,\tag 8.28
$$
and using Lemma {8.}2 gives
$$
\|\ddot{W}_I\|\leq C\|\hat{S}^J\hat{D}_t^2 \triangle q\|
+ C{c}_I^{\,I_1 I_2}\|\check{D}_t\check{g}_{I_1}\|_\infty \|\dot{W}_{I_2}\|
+C{c}_I^{\,I_1 I_2}\|\check{D}_t^2\check{g}_{I_1}\|_\infty \|{W}_{I_2}\|
+ C\tilde{c}_I^{\,I_1 I_2}\|\check{g}_{I_1}\|_\infty \|\ddot{W}_{I_2}\|
\tag 8.29
$$
Hence using interpolation we get
$$
\|\kappa\ddot{W}\|_{r}^{\Cal S}\leq K_3\|\hat{D}_t^2\triangle q\|_{r-1}
+K_3 \sum_{s=0}^{r-1}\|g\|_{r-s,\infty} \|\kappa \ddot{W}\|_s^{\Cal S}
+K_3 \sum_{s=0}^{r}\big(\|\check{D}_t^2 g\|_{r-s,\infty} \|\kappa W\|_s^{\Cal S}
+\|\check{D}_t g\|_{r-s,\infty} \|\kappa \dot{W}\|_s^{\Cal S}\big)
\tag 8.30
$$
We now also want to apply (5.12) in Lemma {5.}2 to $\ddot{W}^a$.
Then $\div \ddot{W}=0$, but the curl of
$\ddot{w}_a=g_{ab}\ddot{W}^b =\pa_a D_t^2\, q
-2(\check{D}_t g_{ab})\dot{W}^b-(\check{D}_t^2 g_{ab}){W}^b$ does not vanish
so using interpolation we get
$$
\| \ddot{W}\|_{r}\leq K_3\|\hat{D}_t^2\triangle q\|_{r-1}
+K_3 \|\ddot{W}  \|_{r}^{\Cal S}
+ \sum_{s=0}^{r}
\big(\|\check{D}_t g\|_{r-s}\|\dot{W } \|_{s}+\|\check{D}_t^2 g\|_{r-s}\|{W } \|_{s}\big)
+ \sum_{s=0}^{r-1}\|g\|_{r-s} \|\ddot{W } \|_{s}\big). \tag 8.31
$$

It remains to prove the estimates for the projection.
We have $W=W_0+W_1$, where $W_0=PW$, and $W_1=g^{ab}\pa_b q$
where $\triangle q=\div W$ and $q\big|_{\pa\Omega}=0$.
The estimates (8.6)-(8.7) follows from (8.2)-(8.4) if
$r\geq 1$, since  $\hat{R}^I \triangle q=\div\,( \hat{\Cal L}_R^I W)$.
Therefore it only remains to consider the case $r=0$.
(8.6), for $r=0$ follows since the projection has norm $1$,
$\|PW\|\leq \|W\|$.
We have $\hat{D}_t W_1^a=(\hat{D}_t g^{ab})g_{bc}
W_1^c+g^{ab}\pa_b D_t q$. Since $W_0$ is divergence free it
follows that $\hat{D}_t W_0$ is divergence free so
 $\hat{D}_t W_0^a =P\hat{D}_t W_0^a =P\hat{D}_t W^a-P\hat{D}_t W_1^a=
 P\hat{D}_t W^a-P\big( (\hat{D}_t g^{ab}) g_{bc}W_1^c\big)$.
 Hence $ \|\hat{D}_t W_0\|\leq K_3 ( \|\hat{D}_t W\|+\|W\|)$
 so (8.7) follows for $r=0$ and $k=1$, and we also have
$ \|\hat{D}_t W_1\|\leq K_3 ( \|\hat{D}_t W\|+\|W\|)$.
 Similarly, $\hat{D}_t^2 W_1^a= (\hat{D}_t^2 g^{ab}) g_{bc} W_1^c
 +2(\hat{D}_t g^{ab}) g_{bc}\big( \hat{D}_t W_1^c-(\hat{D}_t g^{cd} )
 g_{de} W_1^e\big)+g^{ab}\pa_b D_t^2 q$. It therefore follows that
 $\hat{D}_t^2 W_0^a=P\hat{D}_t^2 W_0^a=P\hat{D}_t^2 W^a-
 P\big((\hat{D}_t^2 g^{ab}) g_{bc} W_1^c
 +2(\hat{D}_t g^{ab}) g_{bc}\big( \hat{D}_t W_1^c-(\hat{D}_t g^{cd} )
 g_{de} W_1^e\big)$. Hence
 $\| \hat{D}_t^2 W_0\|\leq K_3\big( \|\hat{D}_t^2 W\|+\|\hat{D}_t
 W\|+ \|W\|\big)$, which proves (8.7) also for $k=2$ and $r=0$.
\endcomment

\head 9. The estimates for the curl.\endhead
We are studying an equation of the general form
$$
\ddot{W}+AW=H,\qquad H=B(W,\dot{W})+F\tag 9.1
$$
Here $B$ is a linear combination of multiplication operators.
Here $A$ is the normal operator and it
projects to the divergence free vector fields even if $W$ is not
divergence free. We have
$\curl AW=0$ and $\div AW=0$ so
$$
\div \ddot{W}=\div H,\qquad
\curl \ddot{w}=\curl \underline{H}\tag 9.2
$$
where we defined $\ddot{w}_a=g_{ab} \ddot{W}^b$.
Now recall that $\dot{w}_a=g_{ab} \dot{W}^b$ so it follows that
$$
|D_t \curl \dot{w}|+|D_t\curl w|\leq
C\big(|\pa D_t g||W|+|\pa g|\pa W|+|\pa\dot{W}|+|\curl\ddot{w}|\big)\tag 9.3
$$
Similarly
$$
|\hat{D}_t \div \dot{W|}+|\hat{D}_t \div W|
\leq C\big( |\div \dot{W}|+ |\div\ddot{W}|\big)\tag 9.4
$$
Hence
$$\multline
|D_t \curl \dot{w}|+|D_t\curl w|
+|\hat{D}_t\div \dot{W}|+|\hat{D}_t \div W|+|\curl \ddot{w}|+|\div \ddot{W}|\\
\leq
C\big( |\pa \dot{W}|+|\pa W|+|\pa g| |\dot{W}|+|\pa g| |W|
+|\div H|+|\curl\underline{H}|\big)
\endmultline\tag 9.5
$$
Since $B$ is of order one and in fact is a multiplication operator it follows
that the terms in $\curl B(W,\dot{W})$ and $\div B(W,\dot{W})$ are also going
to be of the form in the right hand side of (9.5).
However, we need to take a closer look on what the operator $B$ really is
because on the one hand it will give an improved estimates and on the other hand
we want to find out exactly what the constants above depend on:

\proclaim{Lemma {9.1}} Suppose that $L_1 W=F$, where
$L_1 W=\ddot{W}+AW-B(W,\dot{W})$ is given by Lemma {2.}4,
$\dot{W}^a=\hat{D}_t W^a$ and
$\ddot{W}^a=\hat{D}_t^2 W^a$.
Let $w_a=g_{ab} W^b$, $\dot{w}_a=g_{ab}\dot{W}^b$,
$\ddot{w}_a=g_{ab}\ddot{W}^b$. Then
$$
\align
D_t\curl w_{ab}&=\curl\dot{w}_{ab}
+\pa_a \big(\dot{g}_{bc}W^c\big)
-\pa_b \big( \dot{g}_{ac}W^c\big)\tag 9.6\\
D_t\curl \dot{w}_{ab}&=\curl\ddot{w}_{ab}
+\pa_a \big(\dot{g}_{bc}\dot{W}^c\big)
-\pa_b \big( \dot{g}_{ac}\dot{W}^c\big)\tag 9.7\\
\curl \ddot{w}_{ab}&=\curl \underline{F}_{ab}+\curl\underline{B}(\dot{W},W)_{ab}\tag 9.8
\endalign
$$
where $\dot{g}_{ab}=\check{D}_t g_{ab}=D_t\, g_{ab}-\dot{\sigma}g_{ab}$ and
$$
\underline{B}_a(W,\dot{W})=-\big(\dot{{g}}_{ac}-\omega_{ac} -\dot{\sigma} g_{ac}\big)
\dot{W}^c
+\dot{\sigma}\big(\dot{g}_{ac}-\omega_{ac}\big)W^c
-\pa_a q_0\tag 9.9
$$
and $L_1 W=F$.
On the other hand, if $\tilde{w}_a=\dot{w}_a-(\dot{\sigma} g_{ab}+\omega_{ab})W^b $
and $L_1$ is given by (2.54) then
$$
\align
D_t \curl w_{ab}&=\curl\tilde{w}_{ab}
+\pa_a\big((\dot{g}_{bc}+\omega_{bc}+\dot{\sigma} g_{bc})W^c\big)
-\pa_b\big((\dot{g}_{ac}+\omega_{ac}+\dot{\sigma} g_{ac})W^c \big)\tag 9.10\\
D_t\curl \tilde{w}_{ab}&=-\pa_a\big((D_t\,{\omega}_{bc}
 +\ddot{\sigma}g_{bc}) W^c\big)+\pa_b\big((D_t\,{\omega}_{ac}
 +\ddot{\sigma}g_{ac}) W^c\big)
+\curl \underline{F}_{ab}\tag 9.11\\
\curl \dot{w}_{ab}&=\curl\tilde{w}_{ab}
+\pa_a \big((\dot{\sigma} g_{bc}+\omega_{bc}) W^c)
-\pa_b\big((\dot{\sigma} g_{ac}+\omega_{ac}) W^c)\tag 9.12
\endalign
$$
\endproclaim
\demo{Proof} The proof uses Lemma {2.}5 and the identity
$D_t\, w_a=D_t (g_{ab}W^b)=\dot{g}_{ab}W^b+g_{ab}\dot{W}^b$
and (2.54).
\qed\enddemo

Now we want to commute with Lie derivatives ${\Cal L}_R^I$, since the Lie derivative
commutes with the curl:
${\Cal L}_R \curl w=\curl {\Cal L}_R w$.

Using Lemma {5.}2 it follows from Lemma {9.}1 and Lemma {6.}2:
\proclaim{Lemma {9.}2} With notation as in Lemma {9.}1 and
Definition {6.}1 we have
$$\multline
\|D_t \curl w\|_{r-1}+\|D_t\curl \dot{w}\|_{r-1}+\|\curl\ddot{w}\|_{r-1}\leq
2\|\curl\underline{F}\|_{r-1}\\
+K_2 \sum_{s=0}^{r}\big( \|x\|_{r+1-s,\infty}
+\|\dot{x}\|_{r+1-s,\infty}\big) (\|W\|_s+\|\dot{W}\|_s)
\endmultline\tag 9.13
$$
We also have
$$\multline
\|D_t \curl w\|_{r-1}+\|D_t\curl \tilde{w}\|_{r-1}\leq \|\curl \tilde{w}\|_{r-1}
+\|\curl\underline{F}\|_{r-1}\\
+K_3 \sum_{s=0}^{r}
\big( \|x\|_{r+1-s,\infty}
+\|\dot{x}\|_{r+1-s,\infty}
+\|\ddot{x}\|_{r+1-s,\infty}) \|W\|_s
\endmultline\tag 9.14
$$
and
$$
\big|\,\|\curl \dot{w}\|_{r-1}-\|\curl \tilde{w}\|_{r-1}\,\big|
\leq K_2 \sum_{s=0}^{r}\big( \|x\|_{r+1-s,\infty}
+\|\dot{x}\|_{r+1-s,\infty}\big)\|W\|_s\tag 9.15
$$
\endproclaim
\demo{Remark} The difference between on the one hand (9.13) and on the
other hand (9.14)-(9.15) is that the latter does not require estimates for
$\|\dot{W}\|_s$ but instead it requires an extra time derivative of the coordinate.
However, we do control two time derivatives of the coordinates.
\enddemo

\head{10. Existence for the inverse of the
modified linearized operator in the divergence free
class}\endhead
We now first want to show that
$$
L_1 W=\ddot{W}+AW-B_0 W-B_1\dot{W}=F, \qquad W\big|_{t=0}=\tilde{W}_0,
\quad \dot{W}\big|_{t=0}=\tilde{W}_1,\tag 10.1
$$
has a smooth solution $W$:
\proclaim{Theorem {10.}1} Suppose that $x$ and $p$ are smooth,
$p\big|_{\pa\Omega}=0$ and that the coordinate and physical condition
(2.8) and (2.7)
hold for $0\leq t\leq T$. Let $L_1$ be defined by (2.49) and suppose that
$\tilde{W}_0$, $\tilde{W}_1$ and $F$ are smooth and divergence free.
Then (10.1) has a smooth solution for $0\leq t\leq T$.
\endproclaim

In case, $\div V=0$ and $\div F=0$ existence for (10.1)
was proven in \cite{L1}. We now want to generalize this result to prove existence
when $\div V\neq 0$ and $\div F=0$. This is just minor modification of the proof
in \cite{L1},
mostly notational differences, multiplying with $\kappa=\det\,(\pa x/\pa y)$
and $\kappa^{-1}$ since $\div W=\kappa \pa_a (\kappa W^a)$.
We will just give an outline of the proof.

First we note that we can reduce to the case
with vanishing initial data and $F$ vanishing to all orders as $t\to 0$.
In fact, we can get all higher time derivatives by differentiating the equation
with an inhomogeneous term
$$
\hat{D}_t^{k+2} W=B_k\big( W,.,\hat{D}_t^{k+1}W,\pa W,...,
\pa \hat{D}_t^{k} W\big)+\hat{D}_t^k F,\tag 10.2
$$
where $B_k$ are some linear functions followed by projections,
see (10.8) with $I$ consisting of just time derivatives.
Let us therefore define functions of space only by
$$
\tilde{W}^{k+2}=B_k\big(\tilde{W}^0,...,\tilde{W}^{k+1},
\pa \tilde{W}^0,...,\pa \tilde{W}^k \big)\big|_{t=0}
+\hat{D}_t^k F\big|_{t=0},
\qquad k\geq 0   \tag 10.3
$$
Then
$$
\tilde{W}(t,y)=\frac{\kappa(0,y)}{\kappa(t,y)}
\sum_{k=0}^{m-1} \tilde{W}^k(y) t^k/k!\tag 10.4
$$
defines a formal power series solutions at $t=0$. What we are
doing is just expanding $\kappa W$ in a formal power series in
$t$, since $D_t(\kappa W)=\kappa \hat{D}_t W$. Since
$\div\tilde{W}^k=0$ it follows that $\div \tilde{W}=0$.
 We also note that if the initial data
are smooth  then we can construct a smooth approximate solution $\tilde{W}$
that satisfies the equation to all orders as $t\to 0$. This is obtained by multiplying
the $k^{th}$ term in (10.4) by a smooth cutoff $\chi(t/\varepsilon_k)$, to be chosen
below, and summing up the infinite series.
Here $\chi$ is smooth $\chi(s)=1$ for $|s|\leq 1/2$ and $\chi(s)=0$ for $|s|\geq 1$.
The sequence $\varepsilon_k>0$ can then be chosen small enough so that the series
converges in $C^n([0,T],H^m)$ for any $n$ and $m$
if  take $(\|\tilde{W}^k\|_k+1)\varepsilon_k\leq 1/2$.
By replacing $W$ by $W-\tilde{W}$ and $F$ by $F-L_1\tilde{W}$ in (10.1) we reduce
to the situation with vanishing initial data and an inhomogeneous term $F$
vanishing to all order as $t\to 0$.

We will therefore assume that initial data in (10.1) vanishes and that
$F$ is smooth, divergence free and vanishes to all order as $t\to 0$.
 Then existence of a solution $W_\varepsilon$ for the
equation where we replace the normal operator $A$ by the smoothed out
normal operator $A^\varepsilon$, $\varepsilon>0$, in (10.1)
$$
L_1^\varepsilon W_{\varepsilon}=\ddot{W}_{\varepsilon}+A^\varepsilon W_{\varepsilon}
-B_0W_{\varepsilon}-B_1 \dot{W}_{\varepsilon}=F\tag 10.5
$$
follows since all the operators are bounded on $H^r(\Omega)$, see (3.15),
 so its just
an ordinary differential equation in $H^r(\Omega)$, for any $r>0$.
Additional regularity in time follows by applying time derivatives.
This was proven in \cite{L1}.
Lowering the indices in (10.5):
$$
\underline{G}\ddot{W}_{\varepsilon}+\underline{A}^\varepsilon W_{\varepsilon}
-\underline{B}_0W_{\varepsilon}-\underline{B}_1 \dot{W}_{\varepsilon}
=\underline{G}F\tag 10.6
$$
Let $\hat{\Cal L}_T^I $, $I\in{\Cal T}$, stand for a product of modified Lie derivatives,
see section 4, of $|I|$ vector fields in ${\Cal T}$ and let ${W}_{{\varepsilon} I}=
\hat{\Cal L}_T^I W_{\varepsilon}$.
If we repeatedly apply Lie derivatives ${\Cal L}_{T}$ and the
projection, see section 4,
$$
c_I^{\, I_1 I_2} \big(\underline{G}_{I_1}\ddot{W}_{{\varepsilon} I_2}+
\underline{A}_{I_1}^{\varepsilon} W_{{\varepsilon}I_2}
-\underline{B}_{0 I_1}W_{{\varepsilon} I_2}
-\underline{B}_{1 I_1} \dot{W}_{{\varepsilon} I_2}
-\underline{G}_{I_1}F_{I_2}\big)=0
\tag 10.7
$$
where the sum is over all combination of $I_1+I_2=I$ and $c_I^{\, I_1 I_2}=1$.
If we raise the indices again we get
$$
\ddot{W}_{{\varepsilon} I}+
{A}^\varepsilon W_{{\varepsilon} I}
=-\tilde{c}_I^{\, I_1 I_2} \big({G}_{I_1}\ddot{W}_{{\varepsilon} I_2}+
{A}_{I_1}^\varepsilon W_{{} I_2}\big)+c_I^{\, I_1 I_2}\big(
{B}_{0 I_1}W_{{\varepsilon} I_2}
+{B}_{1 I_1} \dot{W}_{{\varepsilon} I_2}
+{G}_{I_1}F_{I_2}\big)
\tag 10.8
$$
where $\tilde{c}_I^{\, I_1 I_2}=1$, if $|I_2|<|I|$,
and $\tilde{c}_I^{\, I_1 I_2}=0$ if $|I_2|=|I|$.

%\subheading{Energy estimates for tangential derivatives}
Let us define energies
$$
E_I^{}=E(W_{ \varepsilon I})=\langle \dot{W}_{{\varepsilon} I},\dot{W}_{{\varepsilon} I}\rangle
+\langle {W}_{{\varepsilon} I},(A^\varepsilon +I){W}_{{\varepsilon} I}\rangle,
\qquad\quad E_s^{{\Cal T}{}}=\sum_{|I|\leq s,I\in{\Cal T}}
\sqrt{E_I^{}}
\tag 10.9
$$
Note that in the sum we also included all time derivatives
$\hat{\Cal L}_{D_t}$. The reason for this is that when calculating commutators
second order time derivatives show up in the first term on the right in (10.7).
As for (3.38) we get by
differentiating (10.9)
$$
\dot{E}_I^{}=
2\langle \dot{W}_{{\varepsilon} I}, \ddot{W}_{{\varepsilon} I}
+(A^\varepsilon+I)W_{{\varepsilon} I}\rangle
+\langle \dot{W}_{{\varepsilon} I}, \dot{G}\dot{W}_{{\varepsilon} I}\rangle
+\langle {W}_{{\varepsilon} I},
(\dot{A}^\varepsilon+\dot{G}){W}_{{\varepsilon} I}\rangle\tag 10.10
$$
Now, $\dot{G}$ is a bounded operator by (3.17).
The last term can be bounded by
$\langle {W}_{{} I},({A}^\varepsilon+I){W}_{{\varepsilon} I}\rangle$
using (4.43) which also holds for $A^\varepsilon$ by (4.37),
and (3.13). Therefore, the last two terms are bounded by $E_r^{{\Cal T}{}}$,
where $r=|I|$. Using (10.8) to estimate the first term we see that
the $L^2$ norm of the last term on the right of (10.8) is bounded by
a constant times
$E_r^{{\Cal T}{}}$ plus
$\|F\|_r^{\Cal T}$ where
$\|F\|_r^{\Cal T}=\sum_{|I|\leq r,\, I\in{\Cal T}} \|\hat{\Cal L}_T^I F\|$,
and $\|F\|=\langle F,F\rangle^{1/2}$. The same is true with the first part
of the first term in on the right in (10.8) since $|I_2|<|I|$ there and
since we have included all time derivatives in the definition of
$ E_s^{{\Cal T}{}}$. It only remains to deal with the second part
of the first term on the right of (10.8).
This term comes from the commutators of $\hat{\Cal L}_T^I$ and $A^\varepsilon$
and we add an additional term to the energy in order to pick it up. Let
$$
D_I^{}=2\tilde{c}_I^{\,I_1 I_2}
\langle W_{{\varepsilon} I},A_{I_1}^\varepsilon
W_{{\varepsilon} I_2}\rangle\tag 10.11
$$
where the sum is over all $I_1+I_2=I$, $|I_2|<|I|$ and $\tilde{c}_I^{\, I_1 I_2}=1$.
This term is lower order, it is again bounded by using (3.13) by the energies
$CE_r^{{\Cal T}{}}E_{r-1}^{{\Cal T}{}}$.
Furthermore
$$
\dot{D}_I^{}=2\tilde{c}_I^{\, I_1 I_2}
\langle \dot{W}_{{\varepsilon} I},A_{I_1}^\varepsilon
W_{{\varepsilon} I_2}\rangle
+\langle W_{{\varepsilon} I},\dot{A}_{I_1}^\varepsilon
{W}_{{\varepsilon} I_2}\rangle
+\langle W_{{\varepsilon} I},{A}_{I_1}^\varepsilon
\dot{W}_{{\varepsilon} I_2}\rangle\tag 10.12
$$
where, by (3.13) the second to last term is bounded by
$CE_r^{{\Cal T}{}}E_{r-1}^{{\Cal T}{}}$
and the last term is bounded by
$CE_r^{{\Cal T}{}}E_r^{{\Cal T}{}}$,
since we have included all time derivatives in the definition (10.9).
Hence, we have proven that
$$
|\dot{E}_I^{}+\dot{D}_I^{}|\leq
C E_r^{{\Cal T}{}} \big(
E_r^{{\Cal T}{}}+\|F\|_r^{\Cal T}\big),
\qquad\quad
|{D}_I^{}|\leq
CE_r^{{\Cal T}{}}E_{r-1}^{{\Cal T}{}},
\qquad r\geq 0, \quad E_{-1}^{\Cal T}=0.\tag 10.13
$$
Using induction and a Gr\"onwall type of argument, see \cite{L1},
it now follows that:
\proclaim{Lemma {10.}2}  There is a constant $C$ depending only on
$t$ and on $x(t,y)$ but not on ${\varepsilon}$ such that
$$
E_r^{{\Cal T}{}}\leq C\int_0^t\|F\|_r^{\Cal T}\, d\tau
\tag 10.14
$$
\endproclaim
In fact, integrating the first inequality in (10.13) from $0$ to $t$,
using that $E_I(0)=D_I(0)=0$, summing over $|I|\leq r$,
and using the second inequality gives
$(E_r^{\Cal T})^2\leq CE_r^{\Cal T} E_{r-1}^{\Cal T}+
C\int_0^t E_r^{\Cal T}(E_r^{\Cal T}+\|F\|_r^{\Cal T})\, d\tau$.
Hence
$$
\overline{E}_r\leq C\overline{E}_{r-1}+C\int_0^t (\overline{E}_r
+\|F\|_r^{\Cal T})\, d\tau,\qquad\quad\text{where}\qquad
\overline{E}_r(t)=\sup_{0\leq \tau\leq t}E_r(\tau).\tag 10.15
$$
Introducing $M_r=\int_0^t\overline{E}_r\, d\tau$, gives
$\dot{M}_r-CM_r\leq C \overline{E}_{r-1}+C\int_0^t\|F\|_r^{\Cal T}\, d\tau$.
Multiplying by the integrating factor $e^{-Ct}$, we see that
$M_r$ is bounded by some constant depending on $t$ times
$C \overline{E}_{r-1}+C\int_0^t\|F\|_r^{\Cal T}\, d\tau$.
Hence for some other constant
$\overline{E}_r\leq C\overline{E}_{r-1}+C\int_0^t\|F\|_r^{\Cal T}\, d\tau$
and (10.14) follows by induction.

From the uniform energy bounds in Lemma {10.}2 it follows that $\|W_\varepsilon\|\leq C$,
where $C$ is independent of $\varepsilon$ so we can choose a weakly convergence subsequence
 $W_{{\varepsilon}_n}$ that converges weakly in the inner product
to $W$ which is also in that space.
Let $U$ be a smooth divergence free vector field which
$0<t<T$ in the support. Then
$$
\int_0^T  \int_\Omega
g_{ab}(L_{1}^{\varepsilon} W_{\varepsilon}^a) U^b\, \kappa dy\, d\tau
=\int_0^T   \int_\Omega g_{ab} W_{\varepsilon}^a (L_1^{{\varepsilon} * }U^b)\,
\kappa dy\, d\tau\tag 10.16
$$
where $L_1^{{\varepsilon} *}$ is the space time adjoint.
 The only term that
depends on $\varepsilon$ in $L_1^{\varepsilon *}$
is $A^{\varepsilon}$, since $A^{\varepsilon}$ is self adjoint.
Since the projection is
continuous it also follows that
$A^{\varepsilon}U\to A U$, as $\varepsilon\to 0$,
 strongly in $L^2$ if $U\in H^1$. Then right hand side of (10.16) therefore
 convergences so we get
$$
 \int_0^T   \int_\Omega g_{ab} W^a (L_1^{ * }U^b)\,\kappa dy\, d\tau
=\int_0^T  \int_\Omega g_{ab} F^a U^b\, \kappa dy\, d\tau\tag 10.17
$$
where now $W$ is the weak limit. Hence $W$ is a weak solution of the equation.
Furthermore $W_\varepsilon$ is divergence free so it follows that
$W$ is weakly divergence free, i.e.
$$
\int_0^T   \int_\Omega  W^a (\pa_a q) \, \kappa\, dy \, d\tau=0\tag 10.18
$$
for all functions smooth functions $q$ that vanishes on the boundary.
We now want to conclude that $W$ has additional regularity so we can integrate by parts
back and conclude that $W$ is a regular solution.

Note that since the curl of a gradient vanishes
$$
\curl \underline{A}^\varepsilon W_\varepsilon=0,\qquad
\text{when}\quad d(y)\geq \varepsilon,\tag 10.19
$$
It follows that the formulas in Lemma {9.}1 hold true for $d(y)\geq \varepsilon$ and
hence
\proclaim{Lemma {10.}3} When $d(y)\geq \varepsilon$ we have
$$\align
|D_t \curl \dot{w}_\varepsilon|_{r-1}^{\Cal U} &\leq
C
\big(|W_\varepsilon|_{r}^{\Cal U}+|\dot{W}_\varepsilon|_{r}^{\Cal U}\big)
+|\curl \underline{F}|_{r-1}^{\Cal U}
\tag 10.20\\
|D_t \curl w_\varepsilon|_{r-1}^{\Cal U}&\leq C
\big(|W_\varepsilon|_{r}^{\Cal U}+|\dot{W}_\varepsilon|_{r}^{\Cal U}\big)\tag 10.21
\endalign
$$
\endproclaim
and by Lemma {5.}2, see the last statement:
\proclaim{Lemma {10.}4}
$$\align
|W|_r^{\Cal U}&\leq C
 \big( |\curl w|_{r-1}^{\Cal U}+|\div W|_{r-1}^{\Cal U}+
|W|_r^{\Cal T}\big)\tag 10.22\\
|\dot{W}|_r^{\, \Cal U}&\leq C\big( |\curl \dot{w}|_{r-1}^{\Cal U}
+|\div \dot{W}|_{r-1}^{\Cal U}+
|\dot{W}|_r^{\Cal T}\big)\tag 10.23
\endalign
$$
\endproclaim
Let $C_0^{\,{\Cal U},\varepsilon}=0$ and for $r\geq 1$ let
$$
C_r^{\, {\Cal U},\varepsilon}=\|\curl\dot{w}_\varepsilon\|_{{\Cal U}^{r-1}(\Omega_\varepsilon)}+
\|\curl w_\varepsilon\|_{{\Cal U}^{r-1}(\Omega_\varepsilon)},
\qquad\text{where}\qquad
\|\beta\|_{\, {\Cal U}^r(\Omega_\varepsilon)}
=\sqrt{\int_{\Omega_\varepsilon}(|\beta|_r^{\Cal U})^2\, \kappa dy}\tag 10.24
$$
and $\Omega_\varepsilon=\{y\in \Omega; d(y)>\varepsilon\}$.
Since $\div W=\div\dot{W}=0$ and since $d(y)\geq \varepsilon$ and
 in the domain of integration in (10.24) it
therefore follows from Lemma {10.}3 and Lemma {10.}4 that
$$
|\dot{C}_r^{\,{\Cal U},\varepsilon}|
\leq C( C_r^{\, {\Cal U},\varepsilon}+E_r^{{\Cal T}})+
C\|\underline{F}\|_{r}^{\Cal U}
\tag 10.25
$$
where $C$ depends on $t$ and $x(t,y)$ but is independent of $\varepsilon$.
This together with Lemma {10.}2 and Lemma {10.}4 now gives us uniform
bounds:
\proclaim{Lemma {10.}5}
$$
\|\dot{W}_\varepsilon\|_{{\Cal U}^r(\Omega_\varepsilon)}
+\|W_\varepsilon\|_{{\Cal U}^r(\Omega_\varepsilon)}
+E_r^{{\Cal T}}
\leq C\int_0^t \|F\|_r^{\Cal U}\, d\tau
\tag 10.26
$$
\endproclaim
We can hence pass to the limit and conclude that
the limit $W$ also satisfies the same estimate and therefore
if we integrate by parts in  (10.17) and (10.18) conclude that $W$ in fact is a
smooth solution:

%\subheading{Existence of regular solutions}

\proclaim{Proposition {10.}6} Suppose $x(t,y)$ is smooth and (2.7) and (2.8) hold
for $0\leq t\leq T$.
Suppose also that $F$ is smooth for $0\leq t\leq T$,
$\div F=0$ and $F$ vanishes to all orders as $t\to 0$.
Then the modified linearized equation (10.1) with vanishing initial data, $\tilde{W}_0=\tilde{W}_1=0$, have a
smooth solution $W$ for $0\leq t\leq T$, satisfying $\div W=0$. Furthermore, the solution
satisfies the estimate:
$$
\|\dot{W}\|_{r}^{\Cal U}
+\|W\|^{{\Cal U}}_r
+E_r^{{\Cal T}}
\leq C_r\int_0^t  \|F\|_r^{\Cal U} \, d\tau,
\qquad r\geq 0.
\tag 10.27
$$
\endproclaim

\head {11. Estimates for the inverse of the modified
linearized operator in the divergence free class}\endhead
%\subheading{The energy estimate for tangential space derivatives}
We will now give improved estimates for the modified linearized
equation
$$
L_1 W=\ddot{W}+AW-B_0 W-B_1\dot{W}=F, \tag 11.1
$$
within the divergence free class. We have
\proclaim{Theorem {11.}1}
Suppose that $x$ and $p$ are smooth,
$p\big|_{\pa\Omega}=0$ and that the coordinate and physical conditions
(2.8) and (2.7)
hold for $0\leq t\leq T$. Let $L_1$ be defined by (2.49) and suppose that
$W$ and $F$ are smooth and divergence free satisfying (11.1) for $0\leq t\leq T$.
 Then if $W\big|_{t=0}=\dot{W}\big|_{t=0}=0$ we have
$$\align
\|\dot{W}\|_{r}+\|W\|_{r}
&\leq K_3 e^{K_3(1+c_0^{-1})T}\, \sum_{s=0}^{r} \underline{{n}}_{r-1-s}\!
\int_0^t\!\!\|{F}\|_{s}\, d\tau,\tag 11.2\\
\|\ddot{W}\|_{r-1}&\leq K_3 e^{K_3(1+c_0^{-1})T}\, \Big(\sum_{s=0}^{r} \underline{{n}}_{r-1-s}\!
\int_0^t\!\!\|{F}\|_{s}\, d\tau+\sum_{s=0}^r \underline{n}_{r-1-s}\|F\|_{s}
\Big)\tag 11.3
\endalign
$$
for $0\leq t\leq T$. If in addition $F\big|_{t=0}=0$ then
for $r\geq 1$ and $0\leq t\leq T$ we have
$$
\|\ddot{W}\|_{r-1}+\|\dot{W}\|_{r-1}+\|W\|_{r-1}+c_0\|W\|_{r}
\leq K_4e^{K_4(1+c_0^{-1})T}\,  \sum_{s=0}^{r-1} \underline{{n}}_{r-1-s}
\int_0^t\!\big(\|\dot{F}\|_{s}+\|F\|_{s}
+\|\curl{F}\|_s\big)\, d\tau\tag 11.4
$$
Here $c_0>0$ is the constant in (2.7), where
$$\align
\underline{n}_{s}&=\sup_{0\leq t\leq T} n_s(t), \tag 11.5\\
n_s(t)&= \|x(t,\cdot)\|_{4+s,\infty}
+\|\dot{x}(t,\cdot)\|_{3+s,\infty}
 +\|\ddot{x}(t,\cdot)\|_{2+s,\infty}   \tag 11.6
\endalign
$$
and $K_3$ is a constant, which depends on  $\underline{n}_{-1}+c_1$,
where $c_1$ is the constant in (2.8).
\endproclaim
For $r=0$ (11.2) is the basic energy estimate
in section 3.
For $r\geq 1$ it follows from first applying Lie derivatives with
respect to space tangential vector fields to the equation and estimating the energy
for these as well as using the evolution equation for the curl and
the fact that we can estimate any derivative by the curl,
the divergence and tangential derivatives. The difference between (11.2) and
(10.27) is, apart from that we keep track of
how the constants depend on the solution we linearize around,
that we only have space derivatives in the norms in (11.2).
The commutators in the energy estimate
are now estimated using the curl as well as the energies of tangential derivatives.
(11.3) follows from (11.2) using (11.1) to estimate
$\ddot{W}$. (11.2) and (11.3) follows from Proposition {11.}4 below.

The importance of (11.4) is that one gets control of an additional space
derivative $c_0\|{W}\|_r$ by only controlling an additional time derivative
and the curl of the right hand side.
(11.4) without the term $c_0\|{W}\|_r$ in the left and
$\|\curl F\|_s$ in the right,
in principle follows from (11.2)
applied to the equation one gets for $\dot{W}$
by commuting $\hat{D}_t$ through $L_1$ in (11.1).
The commutator term $\dot{A}W$
can in principle be controlled by the energy of the same order
but in order not to get constants depending on $\ddot{A}$
we will bound it using an additional space derivative.
$c_0\|{W}\|_r$ can be controlled as follows.
Using the estimate without this term in (11.1) gives control of
$\|A W\|_{r-1}$. By Lemma {5.}4 this gives us control
of $c_0\|{W}\|_r$ if we also control $\|\curl w\|_{r-1}$.
We then use that there is an improved evolution equation for the curl
which only requires control of $\|W\|_r$ and not $\|\dot{W}\|_r$,
by Lemma {9.}2. (11.4) follows from Proposition {11.}9.

Let us rewrite (11.1) slightly
$$
\ddot{W}+AW=H,\qquad\text{where}\qquad  H=B_0 W+B_1\dot{W}+F\tag 11.7
$$
and by (2.51)-(2.52) the operators
$By_1$ and $B_0$ are on divergence free vector fields
$$
B_1\dot{W}^a=-P\big(g^{ab}\big(D_t \,{{g}}_{bc}-\omega_{bc} -2\dot{\sigma} g_{bc}\big)
\dot{W}^c\big),\qquad
B_0 W^a=P\big(g^{ab}\dot{\sigma}\big(D_t\,{g}_{bc}-\omega_{bc}-\dot{\sigma}g_{bc}\big)W^c\big)
\tag 11.8
$$

It follows from (4.47) and (4.49) that
$$
\ddot{W}_I+AW_I=H_{I}+K_I,\qquad
K_I=\tilde{G}_I^{\, I_1 I_2}A_{I_{1}} W_{I_2},\quad
H_{I}=G_I^{\, I_1 I_2}\big(B_{0I_{1}} W_{I_2}+B_{1I_{1}} \dot{W}_{I_2}\big)+F_I
\tag 11.9
$$
where $W_I={\Cal L}_S^I W$, $F_I={\Cal L}_S^I F$,
and $A_{I}$ and $B_{iI}$ are given by (4.41) and (4.43).
Let
$$
E_I=E({W}_I)=\langle \dot{W}_I,\dot{W}_I\rangle+\langle W_I,(A+I)W_I\rangle
\tag 11.10
$$
Then
$$\multline
\dot{E}_I=2\langle \dot{W}_I,\ddot{W}_I+(A+I)W_I\rangle
+\langle \dot{W}_I,\dot{G}\dot{W}_I\rangle
+\langle {W}_I,(\dot{A}+\dot{G}){W}_I\rangle\\
=2\langle \dot{W}_I,K_I+H_I\rangle
+\langle \dot{W}_I,{W}_I\rangle
+\langle \dot{W}_I,\dot{G}\dot{W}_I\rangle
+\langle {W}_I,(\dot{A}+\dot{G}){W}_I\rangle
\endmultline\tag 11.11
$$
The last there terms can be estimated by $E_I$, by (3.42), so we get
$$
|\dot{E}_I|\leq 2\sqrt{E_I} \|K_I+H_I\|
+K_3(1+ c_0^{-1}) E_I\tag 11.12
$$
However,  in order to estimate
the first term we must estimate $\|K_I\|+\|H_I\|$:
\proclaim{Lemma {11.}2}
Let $c_i$, $K_i$, for $i=1,2,3$, $m_s$ and $\dot{m}_s$
 be as in Definitions {5.}2 and {7.}1. We have
$$\align
\|G_I^{\, I_1 I_1}W\|&\leq K_1 m_{s}\|W\|,\qquad s=|I|-|I_1|-|I|_2\tag 11.13\\
\|B_{i I_1} W\|&\leq K_2 \dot{m}_s\|W\|,\qquad s=|I_1|,\quad i=0,1\tag 11.14\\
\|A_{I_1} W\| &\leq K_3 \big(
\dot{m}_{s+1}\|W\|+\dot{m}_{s}\|W\|_1\big),\qquad s=|I_1|\tag 11.15
\endalign
$$
and if $r=|I|$ then
$$\align
\|K_I\|&\leq K_3\sum_{s=0}^{r} \dot{m}_{r+1-s}\|W\|_s\tag 11.16\\
\|H_I\|&\leq K_2 \sum_{s=0}^{r} \dot{m}_{r-s}(\|W\|_s+\|\dot{W}\|_s\big)+\|F\|_r
\tag 11.17
\endalign
$$
\endproclaim
\demo{Proof} The proof of (11.14) and (11.15) uses (4.41) and (4.43)
for $A_I$ and $B_I$, the bounds (3.9) and (3.3) and (7.10)-(7.11)
to estimate the pressure in terms of the coordinate. The proof (11.13)
also uses the interpolation inequalities in Lemma {6.}2. (11.16) and
(11.17) is just a combination of (11.14) and (11.15) with (11.13) and the
interpolation inequalities in Lemma {6.}2. Note the remark after
Definition {5.}2 that $\|W_I\|\leq K_1(\|W\|_s+\|g\|_s \|W\|)$ if
$s=|I|$. A remark about the estimate (11.15) for $A_{I_1}$ is
required. By (4.43) $A_I=A_{\check{S^I} p}$ which can be estimated by
(3.9) if we control $\|\na_N S \check{S}^I p\|_{L^\infty(\pa\Omega)}$.
In (11.17) we claim that this will involve at most $s+2$ space derivatives of
the metric.
In fact,  $\check{S}^I p=S^I p+C^{I_1...I_k}
(S^{I_1}\sigma)\cdot\cdot\cdot (S^{I_{k-1}} \sigma) S^{I_k} p)$
and $S^{I_k}p=0$ on the boundary so it follows that the normal
derivative must fall on $S^{I_k }p$ so the factors $S^{I_j}
\sigma$ never gets differentiated by $\na_N$.
\qed\enddemo

\demo{Definition {11.}1.} Let
$$
E_r^{\Cal S}=\big(\!\!\!\sum_{|I|\leq r, S\in{\Cal S}}\!\!\!\!\!\!  E_I\big)^{1/2}\!\!\!,
\qquad\quad C_r^{\Cal R}=\|\curl w\|_{r-1}^{\Cal R}
+\|\curl\dot{w}\|_{r-1}^{\Cal R},
\qquad \langle W\rangle_{A,r}=\!\!\!\!\!\sum_{|I|\leq r, I\in {\Cal S}}\!\!\!\!\!
\langle W_I,AW_I\rangle  \tag 11.18
$$
where $C_0^{\Cal R}$ should be interpreted as $0$.\enddemo
Summing up the results in Lemma {11.}2, Lemma {9.}2 and Lemma {5.}3 we have:
\proclaim{Lemma {11.}3} We have
$$
|\dot{E}_r^{\Cal S}|
\leq K_3(1+c_0^{-1}) E_r^{\Cal S} +
\sum_{s=0}^r \big(K_2 \dot{m}_{r-s}\|\dot{W}\|_s+K_3\dot{m}_{r+1-s}\|W\|_s
 \big)
+\|F\|_r\tag 11.19
$$
and
$$
|\dot{C}_r^{\Cal R}|+\|\curl\ddot{w}\|_{r-1}\leq K_2\sum_{s=0}^r \dot{m}_{r-s}(\|\dot{W}\|_s+\|W\|_s \big)
+K_1 \sum_{s=0}^r m_{r-s} \|F\|_s\tag 11.20
$$
and
$$
\|\dot{W}\|_r+\|W\|_r+\langle W\rangle_{A,r}
\leq K_1\sum_{s=0}^r m_{r-s}\big(C_s^{\Cal R} +E_s^{\Cal S}\big)\tag 11.21
$$
$$
C_r^{\Cal R} +E_r^{\Cal S}\leq K_1\sum_{s=0}^r m_{r-s}
\big(\|\dot{W}\|_s+\|W\|_s\big)+\langle W\rangle_{A,s}\tag 11.22
$$
\endproclaim
\demo{Proof} The first inequality follows from (11.012) and Lemma {11.}2, the second
from Lemma {9.}2 and the third from Lemma {5.}3.
The last inequality is just due to that $E_r^{\Cal S}$ contains
$\|\kappa W\|_s^{\Cal S}$ and differentiating $\kappa$ produces lower
order terms.
\qed\enddemo

\proclaim{Proposition {11.}4} Let $c_0>0$ and $c_1<\infty$
be such that (2.7) and (2.8) hold and $x$ is smooth for $0\leq t\leq T$.
Let $\underline{\dot{m}}_s=\sup_{0\leq t\leq T} \dot{m}_s(t)$,
where $\dot{m}_s$ is as in Definition {6.}2, and set
$$
E_r=\|\dot{W}\|_r+\|W\|_r+\langle W\rangle_{A,r}\tag 11.23
$$
Then for $r\geq 0$, there is $K_3$, as in Definition {6.}1,
such that, for $0\leq t\leq T$,
$$
E_r(t)\leq K_3 e^{K_3(1+c_0^{-1})T}\,\,
\sum_{s=0}^r \underline{\dot{m}}_{r+1-s}
\Big(E_s(0)+\int_0^t\|F\|_{s}\, d\tau\Big),\tag 11.24
$$
and for $r\geq 1$
 $$
 \|\ddot{W}\|_{r-1}\leq K_3e^{K_3(1+c_0^{-1})T}\Big(
\sum_{s=0}^r \underline{\dot{m}}_{r+1-s}
\Big(E_s(0)+\int_0^t\|F\|_{s}\, d\tau\Big)
+\sum_{s=0}^{r-1} m_{r-1-s}\|F\|_{s}\Big).\tag 11.25
$$
\endproclaim
\demo{Proof} We will prove the estimate for
 $\tilde{E}_r=E_r^{\Cal S}+C_r^{\Cal R}$, in place of $E_r$,
and in view of Lemma {11.}3 and interpolation, $\dot{m}_r m_{s}\leq K_1\dot{m}_{r+s}$,
the estimate for $E_r$ follows from this. By Lemma {11.}3 and
interpolation, $\dot{m}_{r}m_s\leq K_1 \dot{m}_{r+s}$,
we also have
$$
\frac{d\tilde{E}_r}{dt}\leq K_3(1+c_0^{-1}) \tilde{E}_r+
K_3 \sum_{s=0}^{r-1} \dot{m}_{r+1-s} \tilde{E}_s +
K_1 \sum_{s=0}^r m_{r-s} \|F\|_s\tag 11.26
$$
where we also used that $\dot{m}_1\leq c_3$. Let $a=K_3(1+c_0^{-1})$.
Multiplying by the integrating factor we get
$$
(\tilde{E}_r e^{-at})^\prime
\leq e^{-at} K_3
 \big(\sum_{s=0}^{r-1} \dot{m}_{r+1-s}(\tilde{E}_s+\|F\|_s)+m_0\|F\|_r\big),\tag 11.27
$$
 where this is to be interpreted as that the sum is absent if $r=0$.
 Integration of this gives that
 $$
 \tilde{E}_r(t)\leq K_3e^{aT} \Big(\tilde{E}_r(0)
 +\int_0^t \big(\sum_{s=0}^{r-1} \dot{m}_{r+1-s}(\tilde{E}_s+\|F\|_s)
 +n_0\|F\|_r\big) \,d\tau\Big),\qquad t\leq T,
 \tag 11.28
$$
where the sum is to be interpreted as absent if $r=0$.
The proof of the estimate with $\tilde{E}_r$ in place of $E_r$ is
by induction. Since the sum is  absent if $r=0$ it follows for $r=0$.
In general we use the interpolation:
$\dot{m}_{r+1-s}\dot{m}_{s+1-t}\leq C\dot{m}_1 \dot{m}_{r+1-t}\leq
K_3\dot{m}_{r+1-t}$.

To prove the estimate for $\|\ddot{W}\|_{r-1}$ we note that by Lemma {5.}3
it is bounded by the curl and the tangential components:
$$
\|\ddot{W}\|_{r}\leq K_1 \sum_{s=0}^{r} m_{r-s}\big( \|\curl \ddot{w}\|_{s-1}
+\sum_{|I|=s,I\in{\Cal S}} \|\ddot{W}_I\|\big)\tag 11.29
$$
where the curl was estimated in Lemma {11.}3 and the tangential components can be
estimated using the equation $\ddot{W}_I=AW_I+K_I+H_I$ and Lemma {11.}2:
$$
\sum_{|I|\leq r,I\in{\Cal S}}
\|\ddot{W}_I\|
\leq \sum_{s=0}^{r}
\big(K_2 \dot{m}_{r-s}\|\dot{W}\|_s+ K_3\dot{m}_{r+1-s}\|W\|_s\big)+\dot{m}_0 \|W\|_{r+1}
+\|F\|_{r}\tag 11.30
$$
Hence by (11.29), (11.20) and (11.30)
$$
\|\ddot{W}\|_r\leq K_2\sum_{s=0}^{r} \dot{m}_{r-s}\|\dot{W}\|_s
+K_3\sum_{s=0}^{r+1}
\dot{m}_{r+1-s}\|W\|_s+K_1\sum_{s=0}^{r} m_{r-s}\|F\|_{s}\tag 11.31
$$
(11.25) follows from this with $r$ replaced by $r-1$.
\qed\enddemo

%\subheading{Estimates for additional time derivatives}

We now want to get estimates for an additional time derivatives by differentiating
the equation. This gives an estimate for the
normal operator through the equation and this together with estimates for the curl
gives the estimate for an additional space derivative that we are seeking.
Recall that
$$
\ddot{W}+AW=H,\qquad \text{where}\qquad
H=B(W,\dot{W})+F.\tag 11.32
$$
where $B$ given by (2.49) or (2.63) is
$$
\underline{B}_a(W,\dot{W})=-\big(\dot{{g}}_{ac}-\omega_{ac} -\dot{\sigma} g_{ac}\big)
\dot{W}^c
+\dot{\sigma}\big(\dot{g}_{ac}-\omega_{ac}\big)W^c
-\pa_a q_0\tag 11.33
$$
where $\dot{g}_{ab}=\check{D}_t g_{ab}$.
In order to differentiate with respect to time let us now write this in the form
$\underline{G}\ddot{W}+\underline{A} W=\underline{H}$:
$$
g_{ac}\ddot{W}^c+\underline{A}_a W=\underline{B}_a(W,\dot{W})+g_{ac} F^c
\tag 11.34
$$

Differentiating the equation gives
$$
g_{ac}\dddot{W}^c+\underline{A}_a \dot{W}+\dot{\underline{A}}_a W
=D_t \,\underline{B}_a(\dot{W},W)-\dot{g}_{ac}\ddot{W}^c
+\dot{g}_{ac}F^c+g_{ac}\dot{F}^c\tag 11.35
$$
We have
$$
\multline
D_t \,\underline{B}_a(W,\dot{W})=
-\big(\dot{{g}}_{ac}-\omega_{ac} -\dot{\sigma} g_{ac}\big)
\ddot{W}^c
-\big(\ddot{g}_{ac}-\dot{\omega}_{ac}+\dot{\sigma}\omega_{ab}
-\dot{\sigma}\dot{g}_{ac}-\ddot{\sigma} g_{ac}\big)\dot{W}^c\\
+\dot{\sigma}\big(\dot{{g}}_{ac}-\omega_{ac}\big)\dot{W}^c
+\big(\dot{\sigma}(\ddot{g}_{ac}-\dot{\omega}_{ac}+\dot{\sigma}\omega_{ac})
-\ddot{\sigma}(\dot{g}_{ac}-\omega_{ac})\big)W^c-\pa_q D_t q_0
\endmultline\tag 11.36
$$

In conclusion we get
$$
\dddot{W^{}}+A\dot{W}+\dot{A}W=H_1
\qquad\text{where}\qquad
 H_1=B_9 \ddot{W}+B_8\dot{W}+B_7 W+\dot{G}F+\dot{F}\tag 11.37
$$
where
$$\align
{B}_9\ddot{W}^b&=-P\big(g^{ba}\big(2\dot{{g}}_{ac}-\omega_{ac}
-\dot{\sigma} g_{ac}\big)\ddot{W}^c\big)\tag 11.38\\
{B}_8\dot{W}^b&=P\big(g^{ba}
\big(2\dot{\sigma}(\dot{g}_{ac}-\omega_{ac})
-\ddot{g}_{ac}+\dot{\omega}_{ac}+\ddot{\sigma}g_{ac}\big)\dot{W}^c\big)
\tag 11.39\\
{B}_7{W}^b&=P\big(g^{ba}\big(\ddot{\sigma} (\dot{g}_{ac}-\omega_{ac})
+\dot{\sigma}(\ddot{g}_{ac}-\dot{\omega}_{ac}+\dot{\sigma}\omega_{ac})
\big)W^c\big)\rightalignspace \tag 11.40
\endalign
$$

Applying vector fields to (11.37) gives
$$
\dddot{W^{}}_{\!\!I}+A\dot{W}_I + \dot{A} W_I=H_{1 I}+K_{1I},\qquad\text{where}\quad
K_{1I}=-\tilde{G}_I^{\,I_1 I_2} \big(A_{I_1} \dot{W}_{I_2}
+\dot{A}_{I_1}W_{I_2}\big)\tag 11.41
$$
and
$$
{H}_{1 I}=G_I^{\,I_1 I_2}
\big(B_{6 I_1} W_{I_2}\!+B_{7 I_1} \dot{W}_{I_2}\!
+B_{8 I_1} \ddot{W}_{I_2}\big)
+\dot{F}_I +G_I^{\,I_1 I_2} \dot{G}_{I_1} F_{I_2}\tag 11.42
$$
Let
$$
E_{1I}=E(\dot{W}_I)=\langle \ddot{W}_I,\ddot{W}_I\rangle+
\langle \dot{W}_I,(A+I)\dot{W}_I\rangle\tag 11.43
$$
Then
$$\multline
\dot{E}_{1I}=2\langle \ddot{W}_I,\dddot{W^{}}_{\!\! I}+(A+I)\dot{W}_I\rangle
+\langle \ddot{W}_I,\dot{G}\ddot{W}_I\rangle
+\langle \dot{W}_I,(\dot{A}+\dot{G})\dot{W}_I\rangle\\
=-2\langle \ddot{W}_I,\dot{A} W_I\rangle
+2\langle \ddot{W}_I,K_{1I}+H_{1I}\rangle
+\langle \ddot{W}_I,\dot{W}_I\rangle+\langle \ddot{W}_I,\dot{G}\ddot{W}_I\rangle
+\langle \dot{W}_I,(\dot{A}+\dot{G})\dot{W}_I\rangle
\endmultline\tag 11.44
$$
The last three terms are estimated by $E_{1I}$ and the second term is
estimated as before by lower energies:
$$
|\dot{E}_{1I}|\leq 2\sqrt{E}_{1I} \|\dot{A} W_I\|+ 2\sqrt{E_{1I}} \|K_{1I}+H_{1I}\|
+K_3(1+c_0^{-1}) E_I\tag 11.45
$$
However, the estimate of the first term $-2\langle \ddot{W}_I,\dot{A} W_I\rangle$
requires a couple of new observations.
This term could be absorbed by adding
$2\langle \dot{W}_I,\dot{A}W_I\rangle$
to the energy, which instead would produce $2\langle \dot{W}_I,\dot{A}W_I\rangle$
and $2\langle \dot{W}_I,\ddot{A}W_I\rangle$.
However, we want to have an estimate that only requires
two time derivatives of the coordinate and this would require an estimate for
$\ddot{A}$, which requires three time derivatives of the coordinates.
Instead we will use that by Lemma {5.}4 and Lemma {5.}5
$$
\|\dot{A}W_I\|\leq K_3(\|\pa W_I\|+\| W_I\|)\leq K_3(1+c_0^{-1})
\big( \|A W_I\|+\|\curl\underline{W_I}\|+\|W_I\|\big)\tag 11.46
$$
Then there appears to be a loss of regularity, but remember that we now have
an estimate also for $\|\ddot{W}_I\|$ in the energy and by the equation
(11.2) we can estimate $\|AW_I\|\leq \|\ddot{W}_I\|+\|K_I\|+\|H_I\|$,
where the last two terms were estimated in Lemma {11.}2.
$\curl \underline{W_I}$ is by (5.22) up to terms of lower order equal to
${\Cal L}_S^I \curl w$. At this point we have to use that we have an improved
evolution equation for the curl.

\proclaim{Lemma {11.}5} We have
$$\align
\|B_{i I_1} W\|&\leq K_3 \ddot{m}_s\|W\|,\qquad s=|I_1|,\quad i=7,8,9,\tag 11.47\\
\|\dot{A}_{I_1} W\| &\leq K_3 \big(\ddot{m}_{s+1}\|W\|+\ddot{m}_{s}\|W\|_1\big),
\qquad s=|I_1|\tag 11.48
\endalign
$$
and if $r=|I|$ then
$$\align
\|K_{1I}\|&\leq K_3\sum_{s=0}^{r}\ddot{m}_{r+1-s}
\big(\|\dot{W}\|_s+\|W\|_s\big)\tag 11.49\\
\|H_{1I}\|&\leq K_3 \sum_{s=0}^{r}\ddot{m}_{r-s}
\big(\|\ddot{W}\|_s+\|\dot{W}\|_s+\|W\|_s\big)+\|\dot{F}\|_r
+K_2\sum_{s=0}^r \dot{m}_{r-s}\| F\|_s\tag 11.50
\endalign
$$
\endproclaim

\demo{Definition {11.}2} Let
$$
E_{r,1}^{\Cal S}=\big(\sum_{|I|\leq r,S\in{\Cal S}} E_{I,1}\big)^{1/2},
\qquad C_{r,1}^{\Cal R}=\|\curl w\|_r^{\Cal R}+\|\curl \tilde{w}\|_r^{\Cal R}
\tag 11.51
$$
where $\tilde{w}$ is as in Lemma {9.}1.
\enddemo
Summing up the results in Lemma {11.}5, Lemma {9.}2 and Lemma {5.}5 we have:
\proclaim{Lemma {11.}6} We have
$$\multline
\dot{E}_{r,1}^{\Cal S}
\leq K_3(1+c_0^{-1}) E_{r,1}^{\Cal S}\\
 +K_3\sum_{s=0}^r \ddot{m}_{r-s}\big(\|\ddot{W}\|_s+\|\dot{W}\|_s\big)
+K_3\sum_{s=0}^{r+1} \ddot{m}_{r+1-s}\|W\|_s
+\|\dot{F}\|_r+K_2\sum_{s=0}^r \dot{m}_{r-s}\| F\|_s
\endmultline \tag 11.52
$$
and
$$
\dot{C}_{r,1}^{\Cal R}\leq {C}_{r,1}^{\Cal R}+ K_2\sum_{s=0}^{r+1} \ddot{m}_{r+1-s}\|W\|_s
+K_1\|\curl \underline{F}\|_r\tag 11.53
$$
and
$$
\|\ddot{W}\|_r\leq K_1 \sum_{s=0}^r m_{r-s} E_{s,1}^{\Cal S}
+K_3 \sum_{s=0}^{r}\dot{m}_{r-s} \big(\|\dot{W}\|_s+\|W\|_s\big)
+K_1 \sum_{s=1}^{r} m_{r-s}\|F\|_s \tag 11.54
$$
and
$$
c_0\|W\|_{r+1}\leq K_3 \big(C_{r,1}^{\Cal R}+E_{r,1}^{\Cal S}
+K_2 \sum_{s=0}^{r} \dot{m}_{r+1-s}(\|\dot{W}\|_s+\|W\|_s)\big)\tag 11.55
$$
\endproclaim
\demo{Proof} The energy estimate (11.52)
follows from the energy estimate (11.47c)
and the estimates in Lemma {11.}5. The estimate for the curl (11.53)
follows from Lemma {9.}2.
The estimate for $\ddot{W}$ (11.54) uses Lemma {5.}3:
$$
\|\ddot{W}\|_r\leq K_1\sum_{s=1}^r m_{r-s}
\big( \|\curl\ddot{w}\|_{s-1}+E_{s,1}^{\Cal S}\big)+ m_r E_{0,1}^{\Cal S}
\tag 11.56
$$
where the estimate for the curl follows from Lemma {11.}3,
and we use interpolation, $m_s\dot{m}_r\leq K_3 \dot{m}_{s+r}$.
Let us now prove the additional space regularity (11.56).
We have from the equation, $AW_I=-\ddot{W}_I+H_I+K_I$, and Lemma {11.}2
$$
\|W\|_{r,A}^{\Cal S}\leq E_{r,1}^{\Cal S}
+K_2 \sum_{s=0}^{r} \dot{m}_{r+1-s}(\|W\|_s+\|\dot{W}\|_s\big)+\|F\|_r,
\qquad
\|W\|_{s,A}^{\Cal S}
=\!\!\!\!\sum_{|I|=s,I\in{\Cal S}}\!\!\| A\hat{\Cal L}_S^I W\|\tag 11.57
$$
and (11.56) follows since by Lemma {5.}5:
$$
c_0\|W\|_{r+1}\leq K_3\big(C_{r,1}^{\Cal R}
+\| W\|_{r,A}^{\Cal S}
+\sum_{s=0}^{r} m_{r+1-s}\|W\|_s
\big)\qed\tag 11.58
$$
\enddemo

\proclaim{Proposition {11.}7}
Let $c_0>0$ and $c_1<\infty$ be such that (2.7) and (2.8)
hold and $x$ is smooth for $0\leq t\leq T$.
Let $\underline{\ddot{m}}_s=\sup_{0\leq t\leq T} \ddot{m}_s(t)$,
where $\ddot{m}_s$ is as in Definition {6.}2, and set
$$
E_{r,1}=\|\ddot{W}\|_r+\|\dot{W}\|_r+\langle \dot{W}\rangle_{A,r}
+\|W\|_{r}+\langle {W}\rangle_{A,r}
+\|\curl\tilde{w}\|_r+\|\curl w\|_r+c_0\|W\|_{r+1}.\tag 11.59
$$
where $\tilde{w}$ is as in Lemma {9.}1.
Then for $r\geq 0$ there is $K_4$,
as in Definition {6.}1, such that, for $0\leq t\leq T$,
$$
E_{r,1}(t)\leq K_4 e^{K_4(1+c_0^{-1})T}\, \sum_{s=0}^r \underline{\ddot{m}}_{r+1-s}
\Big(E_{s,1}(0)+\int_0^t\big(\|\dot{F}\|_{s}+\|F\|_{s}
+\|\curl{F}\|_s\big)\, d\tau\Big)\tag 11.60
$$
and for $r\geq 1$
 $$\multline
 \|\dddot{W}\|_{r-1}\leq K_4 e^{K_4(1+c_0^{-1})T}\,\sum_{s=0}^r \underline{\ddot{m}}_{r+1-s}
\Big(E_{s,1}(0)+\int_0^t\big(\|\dot{F}\|_s+\|F\|_{s}+\|\curl F\|_s\big)\, d\tau\Big)\\
+K_4 e^{K_4(1+1/c_0)T}\,\sum_{s=0}^{r-1}\underline{\dot{m}}_{r-1-s}(\|F\|_s+\|\dot{F}\|_{s}).
\endmultline\tag 11.61
$$
\endproclaim
\demo{Proof} The proof would be the same as the proof of Proposition {11.}4 apart
from that we must worry more about the possibility of the constant $c_0$ being small.
As in the proof of Proposition {11.}4 the estimate for $E_{r,1}$ would follow
from the same estimate for $\tilde{E}_{r,1}=E_{r,1}^{\Cal S}+C_{r,1}^{\Cal R}
+\tilde{E}_{r}$, where $\tilde{E}_r=E_r^{\Cal S}+C_r^{\Cal R}$ is as in the proof
of Proposition {11.}4. The critical term with $c_0$ is by Lemma {11.}6 and Lemma {11.}3
bounded by the other terms plus lower order terms. Note that by Lemma {11.}3 and interpolation
$\sum_{s=0}^{r}\ddot{m}_{r+1-s}\tilde{E}_r$ bounds
the lower order terms with $\|\dot{W}\|_s$ and $\|W\|_s$, for $s\leq r$ in Lemma {11.}6.
By Lemma {11.}6 and the proof of Proposition {11.}4 we have
$$
\frac{d\tilde{E}_{r,1}}{dt}\leq K_4(1+c_0^{-1}) \tilde{E}_{r,1}+
K_3(1+c_0^{-1}) \sum_{s=0}^{r-1} \ddot{m}_{r+1-s} \tilde{E}_{s,1} +
K_1 \sum_{s=0}^r m_{r-s} \|F\|_s+C\|\dot{F}\|_r \tag 11.62
$$
where we also used that $\ddot{m}_1\leq c_4$. Let
$a= K_4(1+c_0^{-1}) $.
Multiplying by the integrating factor we get
$$
(\tilde{E}_{r,1} e^{-at})^\prime
\leq e^{-at} K_4
 \Big((1+c_0^{-1}) \sum_{s=0}^{r-1} \ddot{m}_{r+1-s}\tilde{E}_{s,1}
+\sum_{s=0}^r m_{r-s} \|F\|_s+\|\dot{F}\|_r\Big),
\tag 11.63
$$
where this is to be interpreted as that the sum is absent if $r=0$.
Integration of this gives that
$$
 \tilde{E}_{r,1}(t)\leq K_4e^{a T} \Big(\tilde{E}_{r,1}(0)
 +\int_0^t \big((1+c_0^{-1})\sum_{s=0}^{r-1} \ddot{m}_{r+1-s}\tilde{E}_{s,1}
  +\sum_{s=0}^r m_{r-s} \|F\|_s
+\|\dot{F}\|_r\big) \,d\tau\Big),
 \tag 11.64
$$
for $t\leq T$, where the sum is to be interpreted as absent if $r=0$.
The proof of the estimate (11.60) with $\tilde{E}_{r,1}$ in place of $E_{r,1}$
follows by induction from (11.64).
Since the sum in (11.64) is absent if $r=0$ it follows that it is true for $r=0$.
In general we use interpolation,
$\ddot{m}_{r+1-s}\ddot{m}_{s+1-t}\leq C\ddot{m}_1 \dot{m}_{r+1-t}\leq
K_4\ddot{m}_{r+1-t}$. (11.61) follows as in the proof of (11.25).
\qed\enddemo

\comment
Similarly if we apply another time derivative to
$\underline{G}\dooot{W}+\underline{A}\dot{W}+\underline{\dot{A}}W
=\underline{G}H_1$ we get
$$
\doooot{W}+A\ddot{W}+2\dot{A}\dot{W}+\ddot{A}W=H_2,
\qquad\text{where}\quad H_2=B_{20} W+...+B_{23}\dooot{W}+G_{20}F+G_{21}\dot{F}
+\ddot{F}\tag 11.65
$$
for some bounded multiplication operators $B_{2i}$ and $G_{2i}$.

\demo{Definition {11.}3} Let $c_5$ be a constant such that
$$
\|x\|_{5,\infty}+\|\dot{x}\|_{4,\infty}+\|\ddot{x}\|_{3,\infty}
+\|\dooot{x}\|_{2,\infty}\leq c_5
\tag 11.66
$$
and let $K_5$ denote a constant which is a continuous function of $c_5$.
Furthermore, let
$$
\ddot{n}_s(t)= \|x(t,\cdot)\|_{5+s,\infty}+\|\dot{x}(t,\cdot)\|_{4+s,\infty}
 +\|\ddot{x}(t,\cdot)\|_{3+s,\infty}+ \|\dooot{x}(t,\cdot)\|_{2+s,\infty},
\tag 11.67
$$
\enddemo

Then we have
\proclaim{Theorem {11.}7} Let $0< T\leq c_0\leq 1$ be such that (1.6)
hold for $0\leq t\leq T$.
Let $\underline{\ddot{n}}_s=\sup_{0\leq t\leq T} \ddot{n}_s(t)$,
where $\ddot{n}_s$ is as in Definition {11.}3, and set
$$
E_{r,2}=\|\dddot{W}\|_r+\|\ddot{W}\|_r++\|\dot{W}\|_{r}
+\|W\|_{r}+\langle \ddot{W}\rangle_{A,r}
+\langle \dot{W}\rangle_{A,r}+\langle {W}\rangle_{A,r}
+c_0\|\dot{W}\|_{r+1}+c_0\|W\|_{r+1}\tag 11.68
$$
Then for $r\geq 0$ there is a constant $K_5$,
as in Definition {11.}3, such that, for $0\leq t\leq T$,
$$
E_{r,2}(t)\leq K_5 \sum_{s=0}^r \underline{\ddot{n}}_{r+1-s}
\Big(E_{s,2}(0)+\int_0^t\big(\|\ddot{F}\|_{s}+\|\dot{F}\|_{s}
+\|F\|_{s}+\|\curl{F}\|_s\big)\, d\tau\Big),\tag 11.69
$$
and
 $$
 \|\doooot{W}\|_{r-1}\leq K_5\sum_{s=0}^r \underline{\ddot{n}}_{r+1-s}
\Big(E_{s,2}(0)+\int_0^t\big(\|\ddot{F}\|_s+\|\dot{F}\|_s+\|F\|_{s}+\|\curl F\|_s\big)\, d\tau\Big)
+\|F\|_{r-1}.\tag 11.60
$$
\endproclaim
We do not use this theorem so we do not give the proof here but state it just
to show that there is a similar gain in regularity for the time derivative.
But this requires three time derivatives of the coordinates so it will not
be of use for us.
\endcomment

\head{12. Existence and estimates for the inverse of the
modified linearized operator for general vector fields}\endhead
In this section we will show that the modified linearized operator
 can be solved for general
vector fields outside the divergence free class, i.e. we solve
$$
L_1 W=F,\qquad \quad W\big|_{t=0}=\dot{W}\big|_{t=0}=0\tag 12.1
$$
when $F$ is not necessarily divergence free.
Below we give estimates for the solution of (12.1)
that are good enough that the linearized operator can be considered as
a lower order modification of (12.1); In the next section we will use these
to prove existence and estimates
also for the inverse of the linearized operator by iteration.
One gets a new iterate
by substituting the previous iterate into the right hand side of (12.3)
and solving for the new iterate in the left hand side. We want estimates
that are good enough that we get the same regularity for the new iterate
so we need estimates for (12.1) that do not loose regularity going from $F$ to $W$.
We have:
\proclaim{Theorem {12.}1}  Let $0<T\leq c_0\leq 1$ and $0<c_1<\infty$
be such that (2.7)-(2.8) hold and $x$ is smooth
for $0\leq t\leq T$. Let $\underline{{n}}_s=\sup_{0\leq t\leq T} {n}_s(t)$,
where ${n}_s$ is as in Definition {6.}2. Then the equation (12.1),
 with $F$ smooth,
 has a smooth solution $W$, for $0\leq t\leq T$. Furthermore, there is $K_4$
as in Definition {6.}1, such that, for $0\leq t\leq T$,
$$
\|\dot{W}\|_{r-1}+\|W\|_{r}
\leq K_4\sum_{s=1}^r \underline{{n}}_{r-s}
\int_0^t\|F\|_{s}\, d\tau, \qquad r\geq 1\tag 12.2
$$
and
$$
\|\ddot{W}\|_{r-1}\leq  K_4\sum_{s=1}^r \underline{{n}}_{r-s}
\int_0^t\|F\|_{s}\, d\tau + K_4\sum_{s=0}^{r-1} \underline{{n}}_{r-1-s}\|F\|_s,
\qquad r\geq 1\tag 12.3
$$
\endproclaim

As in section 3 we can decompose $W=W_0+W_1$
where $W_0$ is divergence free and $W_1$ is the gradient of a function vanishing
at the boundary. By (3.26) $W_0$ satisfies
$$
L_{1} W_{0}=-AW_1+B_{11} \dot{W}_1+B_{01} W_1+PF,
\qquad \quad W_0\big|_{t=0}=\dot{W}_0\big|_{t=0}=0\tag 12.4
$$
where all the terms in the right hand side are divergence free
and $B_{01}$ and $B_{11}$ are bounded operators given by (3.25).
By (3.27)-(3.28) $W_1$ satisfies
$$
W_1^{a}=g^{ab}\pa_b q_1,\qquad \triangle\,
q_1=\varphi,\qquad q_1\big|_{\pa\Omega}=0,\tag 12.5
$$
where
$$
D_t^2\varphi +\ddot{\sigma} \varphi=\div F,\qquad\qquad
 \varphi\big|_{t=0} =D_t \varphi\big|_{t=0} =0
\tag 12.6
$$
The solution of (12.6) is a smooth function if $F$ is smooth so it follows
that that $W_1$ is smooth and hence (12.4) has a smooth solution
 $W_0$ by Theorem {10.}1.
Therefore, we have proven that the modified linearized operator (12.1) has a smooth
solution $W$ if $F$ and $x$ are smooth and the coordinate and physical conditions
are satisfied for $0\leq t\leq T$.
However, in the right hand side of
(12.4) the term $AW_1$ looses space regularity since $A$ is order
one. If we just use Proposition {11.}4 and Proposition {12.}3 below
we are going to get an estimate that looses space regularity
going from $F$ to $W$ in (12.1).  However, because the curl of $AW_1$
vanishes we can use the improved estimate in Proposition {11.}7 that gains an extra
space derivative to handle the term $-AW_1$.
Let us first prove the estimate for (12.5)-(12.6):
\proclaim{Lemma {12}.2} Suppose that
$$
D_t^2\varphi +\ddot{\sigma}\varphi= \hat{D}_t^2
\varphi-2\dot{\sigma}\hat{D}_t \varphi +\dot{\sigma}^2\varphi=f,
\tag 12.7
$$
Let $T<1$ and set
$\underline{\ddot{m}}_s=\sup_{0\leq t\leq T}\ddot{m}_s(t) $, where $\ddot{m}_s$
is as in Definition {6.}2.
Then, there is $K_3$, as in Definition {6.}1, such that, for $0\leq t\leq T$
and $r\geq 1$ we have
 $$
 \align
 \|\dot{\varphi}\|_{r-1}+\|\varphi\|_{r-1}&\leq K_3\sum_{s=0}^{r-1}
  \underline{\ddot{m}}_{r-1-s}\int_0^t
 \!\!\! \|f\|_s\, d\tau,\tag 12.8\\
 \|\ddot{\varphi}\|_{r-1}&\leq K_3\sum_{s=0}^{r-1}
\underline{\ddot{m}}_{r-1-s}\Big(
 \|f\|_s+\! \int_0^t\!\!\! \|f\|_s\, d\tau\Big)\tag 12.9
 \endalign
 $$
\endproclaim
\demo{Proof}
(12.4) is just an ode for each space coordinate, however one just has to
make sure to integrate it up in such a way that we do not get more than
two time derivatives on the metric.
$$
D_t\big( (\hat{D}_t \varphi)^2 +\dot{\sigma}^2 \varphi^2\big)=
2\dot{\sigma}(\hat{D}_t\varphi)^2+2\dot{\sigma}(\ddot{\sigma}-\dot{\sigma}^2)\varphi^2
+ 2(\hat{D}_t\varphi)f,\qquad \hat{D}_t^2 \varphi
-2\dot{\sigma}\hat{D}_t \varphi+ \dot{\sigma}^2 \varphi=f  \tag 12.10
$$
Integrating this in time and space gives the lowest
order energy estimate in (12.8) and the lowest order estimate in (12.9)
follows from this since once we have estimates for the $\varphi$
and $\hat{D}_t\varphi$ we get an estimate for $\hat{D}_t^2 \varphi$ from the
equation.

In order to get (12.8) and (12.9) for higher derivatives we commute through
$\hat{R}^I$, defined in section 4 by $\hat{R}^I f=\kappa^{-1} R^I (\kappa f)$,
where $I=(i_1,...,i_r)$ is a multiindex and $R^I=R_{i_1}\cdot\cdot\cdot R_{i_r}$
is a product of the vector fields in ${\Cal R}$ defined in section 4.
Then $[\hat{D}_t,\hat{R}^I]=$ and with $\varphi_I=\hat{R}^I\varphi$ and
$\dot{\sigma}_I=\hat{R}^I\dot{\sigma}$, we obtain
$$
 \hat{D}_t^2 \varphi_I -2\dot{\sigma}\hat{D}_t \varphi_I+
\dot{\sigma}^2 \varphi_I=f_I,\qquad
f_I=2\tilde{c}^{I_1 I_2}\dot{\sigma}_{I_1} \hat{D}_t \varphi_{I_2}
-\tilde{d}^{I_0 I_1  I_2} \dot{\sigma}_{I_0}\dot{\sigma}_{I_1}
 \varphi_{I_2}+\hat{R}^I f\tag 12.11
$$
were the sums are over all combinations of $I_1+I_2=I$ respectively
$I_0+I_1+I_2=I$ and $\tilde{c}^{I_1 I_2}=1$ and $\tilde{d}^{I_0 I_1  I_2}=1$
unless $I_2=I$ in which case they are $0$.
We can now use (12.10) applied to $f_I$ in place of $f$
and $\varphi_I$ in place of $\varphi$.
Here the terms in $f_I$ are lower order.
\qed\enddemo

Once we get the corresponding bounds for $\varphi$ in terms
of $\div F$, the bounds for $W_1$ follows from Proposition {6.}1:
\proclaim{Proposition {12.}3} Suppose that $W_1^a=g^{ab}\pa_b q$, where
$q\big|_{\pa\Omega}=0$ and $\triangle q=\varphi$ where $\varphi$ satisfies
$$
D_t^2\varphi +\ddot{\sigma}\varphi=\div F,
\tag 12.12
$$
Let $T<1$ and set
$\underline{\ddot{m}}_s=\sup_{0\leq t\leq T}\ddot{m}_s(t) $, where $\ddot{m}_s$
is as in Definition {6.}2.
Then, there is $K_3$, as in Definition {6.}1, such that, for $0\leq t\leq T$,
$$
\|\dot{W}_1(t)\|_r+\|W_1(t)\|_r\leq K_3\sum_{s=1}^r \underline{\ddot{m}}_{r-s}
\big( \|\dot{W}_1(0)\|_s+\|W_1(0)\|_s
+\int_0^t\|F\|_{s}\, d\tau\big),\qquad r\geq 1\tag 12.13
$$
and
$$
\|\ddot{W}_1(t)\|_r\leq  K_3\sum_{s=1}^r \underline{\ddot{m}}_{r-s}
\big( \|\dot{W}_1(0)\|_s+\|W_1(0)\|_s
+\int_0^t\|F\|_{s}\, d\tau+\|F(t)\|_s\big),\qquad r\geq 1\tag 12.14
$$
\endproclaim

We will now further decompose the solution of (12.4) into two parts
$W_0=W_{00}+W_{01}$, where
$$
L_1 W_{01}=-AW_1=F_{01},\qquad\quad W_{01}\big|_{t=0}=\dot{W}_{01}\big|_{t=0}=0\tag 12.15
$$
and
$$
L_1 W_{00}=PF+B_{11} \dot{W}_1+B_{01}W_1=F_{00} ,\qquad
W_{00}\big|_{t=0}=\dot{W}_{00}\big|_{t=0}=0\tag 12.16
$$
For (12.15) we use the estimate in Proposition {11.}7 and for (12.16) we use Proposition {11.}4.
This together with Proposition {12.}3 gives Corollary {12.}4 below.
Our solution to (12.1) is now obtained as a sum of $W=W_1+W_{01}+W_{00}$
so it will satisfy the worst of the estimates in Corollary {12.}4
and this proves Theorem {12.}1.

\proclaim{Corollary {12.}4} Let $0<T\leq c_0\leq 1$ and $c_1<\infty$
be such that (2.7) and (2.8) hold and $x$ is smooth
for $0\leq t\leq T$. Let $\underline{{n}}_s=\sup_{0\leq t\leq T} {n}_s(t)$,
where ${n}_s$ is as in Definition {6.}2.
 Let $W_1$ be the solution of (12.5)-(12.6), let $W_{01}$ be the solution of (12.15)
and let $W_{00}$ be the solutions of (12.16).
 Then there is $K_4$, as in Definition {6.}1,
such that, for $0\leq t\leq T$ and $r\geq 1$,
$$\aligned
\|\dot{W}_1\|_r+\|W_1\|_r &\leq K_4\sum_{s=1}^r \underline{{n}}_{r-1-s}
\!\!\int_0^t\!\!\|F\|_{s}\, d\tau,\\
\|\ddot{W}_1\|_r &\leq K_4\!\sum_{s=1}^r \underline{{n}}_{r-1-s}
\!\big(\!\int_0^t\!\!\! \|F\|_{s}\, d\tau+\|F\|_s\big)
\endaligned \tag 12.17
$$
and
$$
\|\ddot{W}_{01}\|_{r-1}+\|\dot{W}_{01}\|_{r-1}+\|W_{01}\|_{r}
\leq K_4\sum_{s=1}^r \underline{{n}}_{r-s}
\int_0^t\|F\|_{s}\, d\tau \tag 12.18
$$
and
$$\aligned
\|\dot{W}_{00}\|_{r}+\|W_{00}\|_{r}
&\leq K_4\sum_{s=0}^r \underline{{n}}_{r-1-s}
\int_0^t\!\! \|F\|_{s}\, d\tau,\\
\|\ddot{W}_{00}\|_{r-1}&\leq K_4\sum_{s=0}^r \underline{{n}}_{r-1-s}
\int_0^t \!\! \|F\|_{s}\, d\tau
+K_3 \sum_{s=0}^{r-1} \underline{{n}}_{r-1-s}\|F\|_{s}
\endaligned \tag 12.19
$$
\endproclaim
\demo{Proof} (12.17) follows from Proposition {12.}3. By Proposition {11.}7 we have
for $r\geq 1$
$$
\|\ddot{W}_{01}\|_{r-1}+\|\dot{W}_{01}\|_{r-1}+c_0\|W_{01}\|_{r}+\|W_{01}\|_{r-1}
\leq K_4 \sum_{s=0}^{r-1} \underline{{n}}_{r-1-s}
\int_0^t\!\big(\|\dot{F}_{01}\|_{s}+\|F_{01}\|_{s}
+\|\curl{F}_{01}\|_s\big)\, d\tau\tag 12.20
$$
We remark that it follows that also $\ddot{W}_{01}\big|_{t=0}=0$ since
$AW_1\big|_{t=0}=0$.
Here the curl of $F_{01}=AW_1$ vanishes and $\hat{D}_t A W_1=A\dot{W}_1
+\dot{A}W_1-\dot{G}AW_1$ so
$$
\|\dot{F}_{01}\|_{r-1}+\|F_{01}\|_{r-1}
\leq K_4\sum_{s=0}^r \underline{n}_{r-1-s}\big(\|\dot{W}_1\|_s
+\|W_1\|_s\big)\tag 12.21
$$
Using (12.17), (12.20) and (12.21) we obtain (12.18).
Note that the constant $c_0$ in (12.20) can be replaced by $1$ since we have
two consecutive integrals and we assumed that $0\leq t\leq T\leq c_0$.
Finally from Proposition {11.}4 we get
$$\aligned
\|\dot{W}_{00}\|_{r}+\|W_{00}\|_{r}
&\leq K_3 \sum_{s=0}^{r} \underline{{n}}_{r-1-s}\!
\int_0^t\!\!\|{F}_{00}\|_{s}\, d\tau,\\
\|\ddot{W}_{00}\|_{r-1}&\leq K_3\sum_{s=0}^{r} \underline{{n}}_{r-1-s}\!
\int_0^t\!\!\|{F}_{00}\|_{s}\, d\tau+K_3\|F_{00}\|_{r-1}
\endaligned\tag 12.22
$$
Now, the operators $B_{01}$ and $B_{11}$ in (12.16)
 are bounded operator given by (3.25), of the same form
that we have already studied in section 9 and $PF$, the projection is bounded by
the estimates in Proposition {6.}38, so it follows that
$$
\|F_{00}\|_r\leq K_4\sum_{s=0}^r {n}_{r-s}(\|\dot{W}_1\|_s+\|W_1\|_s)
+ K_1\sum_{s=0}^r {n}_{r-s-2}\|F\|_s
\tag 12.23
$$
Combining these inequalities, using interpolation as usual, gives also (12.19).
\qed\enddemo

\head{13. Existence and $L^2$ estimates for the inverse of the
linearized operator.}\endhead
In this section we finally prove existence and estimates for the the
inverse of the linearized operator
$$
L_0 W=F\qquad \quad W\big|_{t=0}=\dot{W}\big|_{t=0}=0,\tag 13.1
$$
where $L_0=\Phi^\prime(x)$ is given by (2.14).
(13.1) can be written
$$
L_1 W=B_3 W+F\qquad \quad W\big|_{t=0}=\dot{W}\big|_{t=0}=0\tag 13.2
$$
where the modified linearized operator $L_1$ is given by (2.49) and $B_3$ is given by (2.57).
In the previous section we proved existence and estimates for the modified linearized
operator $L_1$:
$$
L_1 W=F\qquad \quad W\big|_{t=0}=\dot{W}\big|_{t=0}=0,\tag 13.3
$$
The existence and estimates for (13.3) can now be used to prove existence and
estimate for (13.2), and hence for (13.1), by iteration. We simple define a sequence
by $W_{0}=0$ and for $k\geq 1$:
$$
L_1 W_k=B_3 W_{k-1}+F\qquad \quad W_k\big|_{t=0}=\dot{W}_k\big|_{t=0}=0\tag 13.4
$$
We will use the estimates for (13.3)
to show that $W_k$ converges to a solution of (13.2) and that
the solution of (13.2) satisfies the same as estimates as the solution of (13.3).

\proclaim{Theorem {13.}1}  Let $0<T\leq c_0\leq 1$ and $0<c_1<\infty$
be such that (2.7)-(2.8) hold and $x$ is smooth
for $0\leq t\leq T$. Let $\underline{{n}}_s=\sup_{0\leq t\leq T} {n}_s(t)$,
where ${n}_s$ is as in Definition {6.}2. Then the equation (13.1),
 with $F$ smooth,
 has a smooth solution $W$, for $0\leq t\leq T$. Furthermore, there is $K_4$
as in Definition {6.}1, such that, for $0\leq t\leq T$,
$$
\|\dot{W}\|_{r-1}+\|W\|_{r}
\leq K_4\sum_{s=1}^r \underline{{n}}_{r-s}
\int_0^t\|F\|_{s}\, d\tau, \qquad r\geq 1\tag 13.5
$$
and
$$
\|\ddot{W}\|_{r-1}\leq  K_4\sum_{s=1}^r \underline{{n}}_{r-s}
\int_0^t\|F\|_{s}\, d\tau + K_4\sum_{s=0}^{r-1} \underline{{n}}_{r-1-s}\|F\|_s,
\qquad r\geq 1\tag 13.6
$$
\endproclaim
\demo{Proof} The existence and the
estimates in the theorem for (13.3) were given in Theorem {12.}1.
The estimate for (13.1) follows from the estimate for (13.3) by writing
(13.1) in the form (13.2). If $W$ satisfies
$$
L_1 W=B_3\tilde{W}+F,\qquad W\big|_{t=0}=\dot{W}\big|_{t=0}=0\tag 13.7
$$
where $B_3$ is given by (2.57), then by (13.5) for (13.3)
$$
\|\dot{W}\|_{r-1}+\|W\|_{r}
\leq K_4\sum_{s=1}^r \underline{{n}}_{r-s}
\int_0^t\big(\|F\|_{s}+\|\tilde{W}\|_s\big)\, d\tau,\qquad r\geq 1.\tag 13.8
$$
We claim that (13.5) for (13.3) follows from this with $\tilde{W}=W$,
by induction, for some other $K_4$.
In fact, assume that (13.5) is true for $r\leq k-1$,
then it follows from (13.8) and interpolation that
$$\multline
\|\dot{W}\|_{r-1}+\|W\|_{r}
\leq K_4\sum_{s=1}^r \underline{{n}}_{r-s}
\int_0^t\|F\|_{s}\, d\tau +K_4\int_0^t \|W\|_r\, d\tau\\
+K_4 \sum_{s=1}^{r-1}\sum_{k=1}^{s}\underline{n}_{r-s}
\underline{n}_{s-k}\int_0^t \int_0^{\tau}\|F\|_k\, dz\, \, d\tau
\leq K_4\sum_{s=1}^r \underline{{n}}_{r-s}
\int_0^t\|F\|_{s}\, d\tau +K_4\int_0^t \|W\|_r\, d\tau
\endmultline
\tag 13.9
$$
for some other $K_4$.
By a standard Gronwall argument we can get rid of the $\|W\|_r$,
replacing $K_4$ by some other $K_4$.
Let $g(t)=\int_0^t \|W\|_r\, d\tau $ and $f(t)=\sum_{s=0}^r \underline{{n}}_{r-s}
\int_0^t\|F\|_{s}\, d\tau$.  Then $g^\prime(t)\leq K_4  g+K_4 f$ so
$(g e^{-K_4  t}\big)^\prime\leq K_4 f$ and integrating this up
gives $g\leq K_4 \int_0^t f\, d\tau $ for some other $K_4$ and
for $t\leq T$.

Similarly it follows from (13.5) that the solution of (13.7) satisfies
$$
\|\ddot{W}\|_{r-1}\leq  K_4\sum_{s=1}^r \underline{{n}}_{r-s}
\int_0^t\big(\|F\|_{s}+\|\tilde{W}\|_s\big)\, d\tau
+ K_4\sum_{s=0}^{r-1} \underline{{n}}_{r-1-s}\big(\|F\|_s+\|\tilde{W}\|_s\big),
\qquad r\geq 1\tag 13.10
$$
(13.6) for (13.1) follows from this with $\tilde{W}=W$ using the
estimate (13.5) that we used obtained.

It remains to prove existence for (12.2). We
put up an iteration $W_0=0$ and
$L_1 W_{k}=F+(L_1-L_0)W_{k-1}$, for $k\geq 1$.
Then $L_1 W_1=F$ so $W_1$ satisfies the desired estimate and is smooth.
Let $\overline{W}_k=W_k-W_{k-1}$, for $k\geq 1$. Then $\overline{W}_1=W_1$ and
$L_1\overline{W}_k=(L_1-L_0)\overline{W}_{k-1}$, for $k\geq 2$.
In conclusion
$$
L_1 \overline{W}_1=F,\qquad L_1 \overline{W}_k =B_3 \overline{W}_{k-1},
\qquad k\geq 2,\qquad \overline{W}_k\big|_{t=0}=
\dot{\overline{W}}_k\big|_{t=0}=0\tag 13.11
$$
where $B_3$ is a bounded operator given by (2.57).
Using the estimate (13.8) for each $k$
$$
\sum_{k=1}^{N}\sup_{0\leq \tau\leq t} \big(
\|\dot{\overline{W}}_k(\tau,\cdot)\|_{r-1}+\|\overline{W}_k(\tau,\cdot)\|_{r}\big)
\leq K_4\sum_{s=1}^r \underline{{n}}_{r-s}
\int_0^t\big(\|F\|_{s}+\sum_{k=1}^{N-1}\|\overline{W}_k\|_s\big)\, d\tau,\qquad r\geq 1.
\tag 13.12
$$
Note that the supremum is inside the sum since we use (13.8) for each $\overline{W}_k$
and since in the left of (13.8) we may take the supremum of $\tau \leq t$.
The same argument that lead to the proof of the estimate (13.5) for (13.1)
from (13.8) now gives the uniform estimate
$$
\sum_{k=1}^{N}\sup_{0\leq \tau\leq t} \big(\|\dot{\overline{W}}_k(\tau,\cdot)\|_{r-1}
+\|\overline{W}_k(\tau,\cdot)\|_{r}\big)
\leq K_4\sum_{s=1}^r \underline{{n}}_{r-s}
\int_0^t\|F\|_{s}\, d\tau, \qquad r\geq 1,\tag 13.13
$$
where $K_4$ is independent of $N$.
(One replaces the sum in the right of (13.12)
by the larger sum in the left of (13.12).)
Similarly the uniform estimates corresponding to (13.5) also hold
as is seen by using (13.10) for each $k$ and replacing the
sum in the right by the larger sum in the left and using (13.13)
$$
\sum_{k=1}^{N}\sup_{0\leq \tau\leq t}
\|\ddot{\overline{W}}_k(\tau,\cdot)\|_{r-1}\leq  K_4\sum_{s=1}^r \underline{{n}}_{r-s}
\int_0^t\|F\|_{s}\, d\tau + K_4\sum_{s=0}^{r-1} \underline{{n}}_{r-1-s}
\sup_{0\leq \tau\leq t}\|F(\tau,\cdot)\|_s,
\qquad r\geq 1\tag 13.14
$$
It follows that $W_N=\sum_{k=1}^N\overline{W}_k$
is a Cauchy sequence in $C^2([0,T],H^{r-1}(\overline{\Omega}))$,
for any $T$,
and hence there is a limit $W\in C^2([0,T],H^{r-1}(\overline{\Omega}))$,
for any $T$.
Additional regularity in time follows from differentiating this equation.
We have already proved that $\hat{D}_t^2 W=AW+B_{0}W+B_1\dot{W}+B_3 W\in
C^1([0,T],H^{r-2}(\overline{\Omega}))$, i.e. $\hat{D}_t^2 W$ is continuously
differentiable with respect to time so
$ W\!\in\!  C^3\big([0,T],H^{r-2}(\overline{\Omega})\big)$ and so on.
Since this argument is true for any $r$ it follows that $W$ is smooth.
\qed\enddemo

\comment
\proclaim{Theorem {12.}4}  Let $0<T\leq c_0\leq 1$ be such that (1.6) hold
for $0\leq t\leq T$. Let $\underline{{n}}_s=\sup_{0\leq t\leq T} {n}_s(t)$,
where ${n}_s$ is as in Definition {6.}2. Then the equation (12.2), with $F$ smooth,
 has a smooth solution $W$, for $0\leq t\leq T$. Furthermore, there is $K_4$
as in Definition {6.}1, such that, for $0\leq t\leq T$,
$$
\|\dot{W}\|_{r-1}+\|W\|_{r}
\leq K_4\sum_{s=0}^r \underline{{n}}_{r-s}
\int_0^t\|F\|_{s}\, d\tau ,\tag 13.8
$$
and
$$
\|\ddot{W}\|_{r-1}\leq  K_4\sum_{s=0}^r \underline{{n}}_{r-s}
\int_0^t\|F\|_{s}\, d\tau + K_4\sum_{s=0}^{r-1} \underline{{n}}_{r-1-s}\|F\|_s\tag 13.13
$$
\endproclaim
Summarizing, we have proven that
\proclaim{Theorem {12.}5} Let $0<T\leq c_0\leq 1$ and $c_1$ be constants
and $x$ a smooth function on $\overline{\Omega}\times[0,T]$, such that (2.7) and (2.8) hold
for $0\leq t\leq T$, where $p$ is given by solving (2.4).
 Then the linearized equation (12.2)
has a smooth solution on $\overline{\Omega}\times[0,T]$, such that
$$
\multline
\|\ddot{W}\|_{r-1,T}+\|\dot{W}\|_{r-1,T}+\|W\|_{r,T}\\
\leq K_4\Big(
\|F\|_{r,T} +\|F\|_{0,T}
\big(\|x\|_{r+4,\infty,T}+\|\dot{x}\|_{r+3,\infty,T}
+\|\ddot{x}\|_{r+2,\infty,T}\big)\Big)
\endmultline \tag 12.32
$$
where
$\|W\|_{r,T}=\sup_{0\leq t\leq T}
 \|W(t,\cdot)\|_{H^r(\Omega)}$, $\|x\|_{r,\infty, T}=\sup_{0\leq t\leq T}
 \|x(t,\cdot)\|_{C^r(\Omega)}$ and
$K_4=K_4(c_4)$ is a continuous function of $c_4$, where $c_4$ is
a constant such that
$$
\|x\|_{4,\infty,T}+\|\dot{x}\|_{3,\infty,T}
+\|\ddot{x}\|_{2,\infty,T}+c_1\leq c_4 .\tag 12.33
$$
\endproclaim
\endcomment

\head 14. Estimates for the physical and coordinate conditions.\endhead

We assume that the physical condition and the coordinate condition hold
initially at time $0$ for some constants $c_0>0$ and $c_1<\infty$
and we need to show that this implies that they will hold with $
c_0$ replaced by $c_0/2$ and $c_1$ replaced by $2c_1$,  for
 $0\leq t\leq T$, if $T$ is sufficiently small.

 Let us introduce the space time norms:
$$
\||u\||_{r}=\sup_{0\leq t\leq T}\|u(t,\cdot)\|_{r,\infty},
\qquad\quad
\||u\||_{r,k}=\||u\||_{r}+...+\||D_t^k u\||_{r}\tag 14.1
$$
We have:
 \proclaim{Lemma {14.}1} Let $M(t)=\sup_{y\in\Omega}\sqrt{|\pa x/\pa y|^2+|\pa y/\pa
 x|^2}$. Then
 $$
 M(t)\leq 2M(0),\qquad \text{for}\qquad t\leq T,\qquad
 \text{if}\qquad T \||\dot{x}\||_{1}M(0)\leq 1/8\tag 14.2
 $$
 Let $N(t)=\sup_{y\in\pa\Omega} |\na_N p|^{-1}$. Then assuming
 that $T$ is so small that (14.2) hold we have
 $$
 N(t)\leq 2N(0)\qquad \text{for}\qquad t\leq T,\qquad
 \text{if}\qquad T \||\dot{p}\||_{1}M(0)N(0)\leq 1/8\tag 14.3
 $$
 \endproclaim
 \demo{Proof} We have $|D_t\, \pa x/\pa y|\leq \||\dot{x}\||_{1} $
 and $|D_t\, \pa y/\pa x|\leq |\pa y/\pa x|^2 |D_t \pa x/\pa y|$ so
 $ M^\prime(t)\leq (1+M^2) \||\dot{x}\||_{1}\leq 2M^2
 \||\dot{x}\||_{1}$, since also $M(t)\geq 1$. Hence
 $$
 M(t)\leq M(0) \big(1-2\||\dot{x}\||_{1} M(0) t\big)^{-1}, \qquad
 \text{when}\qquad 2\||\dot{x}\||_{1} M(0) t< 1.\tag 14.4
 $$
 Now, $\na_N p=N^a \pa_a p$, where $N$ is the unit normal,
 so $D_t \na_N p=\na_N D_t p+ (D_t N^a) \pa_a p=\na_N D_t p+(D_t
 N^a) g_{ab} N^b \na_N p$, since $p\big|_{\pa\Omega}=0$.
 Furthermore $0=D_t (g_{ab} N^a N^b)= 2g_{ab} (D_t N^a) N^b+(D_t
 g_{ab}) N^a N^b$ and $N^a=(\pa y^a/\pa x^i) N^i$, where
 $\delta_{ij} N^i N^j=1$. Hence
 $|D_t\na_N p|\leq M \big(|\pa D_t p|+ |\pa D_t x| |\na_N p|\big)$
 Therefore if $N(t)=\sup_{y\in\pa\Omega}|\na_N p|^{-1}$, we have
 $N^\prime\leq M \||\dot{p}\||_{1}N^2 + M \||\dot{x}\||_{1}N/2$
 and if we use (14.2) and multiply with the
 integrating factor, $\tilde{N}(t)=N(t)e^{-tM(0) \||\dot{x}\||_{1}}$
 we get $\tilde{N}^\prime\leq 2e^{1/8}M(0) \||\dot{p}\||_{1}
 \tilde{N}^2$. Hence
 $$
 {N}(t)\leq {N}(0) e^{1/8} \big(1-N(0)2e^{1/8} M(0) \||\dot{p}\||_{1}
 t\big)^{-1}\!\!\!\!,\quad \text{when}\quad
 N(0)2 e^{1/8}M(0) \||\dot{p}\||_{1}t<1\tag 14.5
 $$
 This proves the lemma. \qed
 \enddemo

It now follows from Lemma {14.}2:
\proclaim{Lemma {14.}2} Let $x_0$ be the approximate solution satisfying (2.12)
and suppose that (2.7) and (2.8) holds when $t=0$.
Then there is a $T_0>0$,  depending only and an upper bound for
$\||x_0\||_{4,2}$, $c_1$ and $c_0^{-1}$ such (2.7) and (2.8) hold for $0\leq t\leq T$
with $c_0$ replaced by $c_0/2$ and $c_1$ replaced by $2c_1$ provided that
$$
0<T\leq T_0,\qquad \text{and}\qquad \||x-x_0\||_{4,2}\leq 1,\qquad\text{and}
\qquad  (x-x_0)\big|_{t=0}=D_t (x-x_0)\big|_{t=0}=0 \tag 14.6
$$
\endproclaim
\demo{Proof} We need to satisfy the conditions (14.2) and (14.3) in Lemma {14.}1.
Since $\||\dot{x}\||_1\leq \||x_0\||_{4,2}+1$ (14.2) hold if
$T\leq (8c_1 (\||x_0\||_{4,2}+1))^{-1}$. To satisfy (14.3) we use the estimate
in Lemma {6.}4, where $K_3$ is as in Definition {6.}1, to
obtain $\|\dot{p}\|_{1,\infty}\leq F\big(\|x\|_{3,\infty}+\|\dot{x}\|_{2,\infty}
+\|\ddot{x}\|_{1,\infty}\big)$ for some increasing function $F$.
Hence (14.3) hold if $T\leq c_0(8c_1+F(\||x_0\||_{4,2}+1))^{-1}$.
\qed\enddemo

\head 15. Tame $L^\infty$ estimates for the inverse of the linearized operator.\endhead
We are now going to modify the estimate for the inverse of the linearized operator
in Theorem {13.}1 so it can be used with the Nash-Moser
inverse function theorem in section 18. We want tame estimates for the inverse of the
linearized operator
$$
\Phi^\prime(x)\delta x= \delta \Phi,\qquad 0\leq t\leq T,\qquad\delta x\big|_{t=0}=D_t\,
\delta x\big|_{t=0}=0\tag 15.1
$$
but the norms in Theorem {13.}1 are in terms of
$W^a=\delta x^i\pa y^a/\pa x^i$ and $F^a=\delta \Phi^i\pa y^a/\pa x^i$
and we like to see our operator as
an operator on $\delta x$. Using interpolation and Theorem {13.}1 we get
$$\multline
\|\delta \ddot{x}\|_{r}+\|\delta\dot{ x}\|_{r}+\|\delta x\|_{r}
\leq K_2(\| \ddot{W}\|_{r}+\| \dot{W}\|_{r}+\| W\|_{r})
+K_2 (\|\ddot{x}\|_{r+1}+\|\dot{ x}\|_{r+1}+\|x\|_{r+1})
(\| \ddot{W}\|+\| \dot{W}\|+\| W\|)\\
\leq K_4\sup_{0\leq \tau\leq t} \|F(\tau,\cdot)\|_{r+1}+
 K_4\sup_{0\leq \tau\leq t} (\|\ddot{x}(\tau,\cdot)\|_{r+4,\infty}
+\|\dot{x}(\tau,\cdot)\|_{r+4,\infty}+\|x(\tau,\cdot)\|_{r+4,\infty})
\sup_{0\leq \tau\leq t} \|F(\tau,\cdot)\|_{1}\\
\leq K_4\sup_{0\leq \tau\leq t} \|\delta\Phi(\tau,\cdot)\|_{r+1}+
 K_4\sup_{0\leq \tau\leq t} (\|\ddot{x}(\tau,\cdot)\|_{r+4,\infty}
+\|\dot{x}(\tau,\cdot)\|_{r+4,\infty}+\|x(\tau,\cdot)\|_{r+4,\infty})
\sup_{0\leq \tau\leq t} \|\delta\Phi(\tau,\cdot)\|_{1}
\endmultline \tag 15.2
$$

Another issue is that we have $L^2$ estimates
of $\delta x$ but we need $L^\infty$ estimates for $x$.
The $L^2$ norm is bounded by the $L^\infty$ norm and the $L^\infty$
is by Sobolev's lemma bounded by the $L^2$ norm of an additional $n/2$ derivatives so one
can obviously turn one into the other with an additional loss:
$$
\|u(t,\cdot)\|_{r}\leq c_r\|u(t,\cdot)\|_{r,\infty}\leq C_r \|u(t,\cdot)\|_{r+r_0},\qquad
 r_0=[n/2]+1 \tag 15.3
$$
Furthermore, the Nash-Moser theorems that we will follow are in terms of
 H\"older spaces,
but one can obviously also turn H\"older norms into $L^\infty$ norms with a loss of an additional derivate:
$$
C_k^{-1}\|u(t,\cdot)\|_{k,\infty}\leq \|u(t,\cdot)\|_{a,\infty}\leq C_k
\|u(t,\cdot)\|_{k+1,\infty}, \qquad k\leq  a\leq k+1\tag 15.4
$$
where $\|u(t,\cdot)\|_{a,\infty}$ denotes the H\"older norms in section 17.
Let us now introduce the norms
$$
\||u\||_{a,k}
=\||u\||_{a}+...+\||D_t^k u\||_{a},
,\qquad \text{where}\quad
\||u\||_a=\sup_{0\leq t\leq T} \|u(t,\cdot)\|_{a,\infty}\tag 15.5
$$
It follows that if (2.7) and (2.8) hold then (15.1) has a solution that satisfies
$$
\||\delta x\||_{a,2}
\leq K_4\big( \||\delta \Phi\||_{a+r_0+2}
+\||\delta \Phi\||_{1}\, \||x\||_{a+r_0+6,2} \big),\qquad a\geq 0\tag 15.6
$$

We in fact want to solve for $u$ in (2.13):
 $$
 \tilde{\Phi}(u)=\Phi(u+x_0)-\Phi(x_0)=f_\delta \tag 15.7
 $$
Then $\tilde{\Phi}^\prime(u)=\Phi^\prime(u+x_0)$
and the norm of $x$ in (15.6) may be replaced by the norm of $u=x-x_0$
since
 $$
\||x\||_{a,2}
\leq \||x-x_0\||_{a,2}+\||x_0\||_{a,2}\leq \||x-x_0\||_{a,2}+C_a\tag 15.8
$$
for some constant $C_a$ depending on $x_0$. Hence we have proven:
 \proclaim{Proposition {15.}1}
Suppose that $x$ is smooth for $0\leq t\leq T$ and that the conditions
in Lemma {14.}2 hold. Then if $\delta \Phi$ is smooth for $0\leq t\leq T$
(15.1) has a smooth solution $\delta x$.
Furthermore there are constants $C_a$, depending on the approximate solution
 $x_0$, on $(c_0,c_1)$ in (2.7)-(2,8) and on $a$, such that
$$
\||\delta x\||_{a,2}\leq C_a
 \big( \||\delta\Phi\||_{a+r_0+2}+ \|| x-x_0\||_{a+r_0+6,2}
 \|\delta\Phi\||_{1}\big),\qquad a\geq 0\tag 15.9
 $$
provided that
$$
\||x-x_0\||_{4,2}\leq 1\tag 15.10
$$
\endproclaim

\comment
\proclaim{Proposition {15.}1} Suppose that $0<T< c_0$ and $c_1$ are constants such that
$x\in C^2([0,T],C^\infty(\overline{\Omega}))$ and
(2.7)-(2.8) holds for $0\leq t\leq T$, where $p$ is given by (2.4).
 Then, the linearized equation (15.1) has a solution
 $\delta x\in C^2([0,T],C^\infty(\overline{\Omega}))$
if $\delta \Phi\in  C([0,T],C^\infty(\overline{\Omega}))$.
The solution satisfies the estimates
$$
\||\delta x\||_{a,2}
\leq
C_a\big( \||\delta \Phi\||_{a+r_1}
+\||\delta \Phi\||_{1}\, \||x\||_{a+r_2,2} \big),\qquad a\geq 0\tag 15.7
$$
where $r_1=[n/2]+3$, $r_2=[n/2]+7$. Here $K_4$, depending on $a$,
 is as continuous function of $c_4$, where
$c_4$ is a constant such that
$$
\|x\|_{4,2}+c_1\leq c_4\tag 15.8
$$
\endproclaim

\demo{Proof} Only, the fact that we replaced $C^\infty([0,T]\times \overline{\Omega})$
by $ C^2([0,T],C^\infty(\overline{\Omega}))$ and
 $C([0,T],C^\infty(\overline{\Omega}))$ remains to be proved.
First that we have existence for $\delta \Phi\in
C([0,T],C^\infty(\overline{\Omega}))$ and that then $\delta x\in
C^2([0,T],C^\infty(\overline{\Omega}))$ follows directly from
applying the estimate above to differences, since
$C^\infty\big([0,T],C^\infty(\overline{\Omega})\big)$ is dense in
$ C([0,T],C^\infty(\overline{\Omega}))$, i.e. we can find a
sequence which converges uniformly in
$C\big([0,T],C^k(\overline{\Omega})\big)$, for any $k$. Also, the
fact that we have existence for $x\in
C^2([0,T],C^\infty(\overline{\Omega}))$ follows from this estimate
and the estimate for the second variational derivative in the next
section. In fact, let
 $x_n\in C^\infty\big([0,T],C^\infty(\overline{\Omega}) \big)$
be a sequence which tends to $x\in
C^2([0,T],C^\infty(\overline{\Omega}))$, i.e. $x_n\to x$ uniformly
in $ C^2([0,T],C^k(\overline{\Omega}))$, for any $k$, and let
$\delta x_n$ be defined by
$$
\Phi^\prime(x_n)\delta x_n =\delta \Phi\tag 15.8
$$
then the above estimate gives uniform bounds for $\delta x_n $ and
$$
\Phi^\prime(x_n)(\delta x_n-\delta x_m)=
(\Phi^\prime(x_m)-\Phi^\prime(x_n))\delta x_m\tag 15.9
$$
Here
$$
(\Phi^\prime(x_m)-\Phi^\prime(x_n))\delta x_m=
\int_0^1 (1-s) \Phi^{\prime\prime}(x_n+s(x_m-x_n))(x_m-x_n,\delta x_m)\, ds
\tag 15.10
$$
can be estimated by the estimate in Proposition {16.}1. Hence it
follows that the right hand side of (15.9) tends to $0$ in
$C([0,T],C^\infty(\overline{\Omega}))$, when $m,n\to\infty$, i.e.,
it is a Cauchy sequence in $C([0,T],C^k(\overline{\Omega}))$ for
any $k$. Using the estimate (15.7) applied to (15.9) we get that
$\delta x_n-\delta x_m$ tends to $0$ in
$C^2([0,T],C^\infty(\overline{\Omega}))$, when $m,n\to\infty$.
\qed\enddemo
\endcomment

\head 16. Regularity properties of the Euler map
and tame estimates for the second variational derivative.\endhead
Recall that the Euler map is given by
$$
\Phi(x)_i=D_t^2 x_i+\pa_i p,\qquad \text{in}\qquad [0,T]\times\Omega,
\qquad\text{where}\qquad \pa_i=\frac{\pa y^a}{\pa x^i}\pa_a \tag 16.1
$$
where $p=\Psi(x)$ is given by solving
$$
\triangle p=-(\pa_i V^k)\pa_k V^i,\qquad V^i=D_t\, x^i,\qquad p\big|_{\pa\Omega}=0
\tag 16.2
$$

We will now discuss the regularity properties of $\Phi$ needed
and the definition of derivatives of $\Phi$:
 Let
${\Cal F}= C^\infty\big([0,T]\times \overline{\Omega}\big)$,
${\Cal F}_{M}=\{x\in{\Cal F};\, |\pa x/\pa y|+|\pa y/\pa x|<M\}$ and let
$I_k=I\times\cdot\cdot\cdot\times I$ be $k$ copies of
$I=[-\varepsilon,\varepsilon]$, $\varepsilon>0$.
Suppose that $\overline{x}\in C^m(I_k,{\Cal F}_M)$, $m\geq k$
then we claim that $\Phi(\overline{x})\in  C^m(I_k,{\Cal F})$.
In fact, by the proof of Lemma {7.}3 $\overline{p}=\Psi(\overline{x})\in C^m(I_k,{\Cal F})$,
since there $t\in R$ was just any parameter and we can replace
it by $t\in\bold{R}^k$ and replace the
derivatives with respect to $t$ by partial derivatives.

\demo{Definition {16.}1} Suppose that
 $x\in {\Cal F}= C^\infty\big([0,T]\times \overline{\Omega}\big)$
and $w_j\in {\Cal F}$, for $j\leq k$.
Set $\overline{x}=x+r_1 w_1+...+r_k w_k$ and suppose that
$\Phi(\overline{x})$ is a $C^k$ function of $(r_1,...,r_k)$ close to $(0,...,0)$
with values in ${\Cal F}$.
We define the $k$:th (directional) derivative of $\Phi$ at
the point $x$ in the
directions $w_i$, $i=1,..,k$ by
$$
\Phi^{(k)}(x)(w_1,...,w_k)=\frac{\pa}{\pa r_1}\cdot\cdot\cdot
\frac{\pa}{\pa r_k} \Phi(\overline{x})\big|_{r_1=...=r_k=0},
\qquad \overline{x}=x+r_1 w_1+...+r_k w_k\tag 16.3
$$
We say that $\Phi(x)$ is $k$ times differentiable at $x$
 if $\Phi(\overline{x})$ is a $C^k$ function
of $(r_1,...,r_k)$ close to $(0,...,0)$ with values in ${\Cal F}$,
and if $\Phi^{(j)}(x)(w_1,...,w_j)$
is linear in each of the arguments $w_1,...,w_j$, for $j\leq k$.
\enddemo

It is clear that (16.3) is independent of the order of differentiation,
debut to conclude that it is multi linear in $w_1,...,w_k$
one also needs to assume that it is continuous as a functional of
$x,w_1,...w_k$, see \cite{Ha}. We instead take (16.3) as the definition
of the derivative and once we calculated it the linearity follows
by inspection in our case. We will assume that $\Phi$ is twice differentiable
in which case it follows from the above definition that Taylor's formula
with integral reminder of order two hold:
$$\align
\big(\Phi^\prime(v)-\Phi^\prime(u)\big)w
&=\int_0^1 \Phi^{\prime\prime}\big(u+s(v-u)\big)(v-u,w)\, ds\tag 16.4\\
\Phi(v)-\Phi(u)-\Phi^\prime(u)(v-u)
&=\int_0^1 (1-s)\Phi^{\prime\prime}\big(u+s(v-u)\big)(v-u,v-u)\, ds\tag 16.5
\endalign
$$
The Nash-Moser technique uses these reminder formulas
together with tame estimates for the second variational derivative that
we now will derive:
\proclaim{Proposition {16.}1} Suppose that $x$ is smooth for $0\leq t\leq T$ and that the
conditions in Lemma {14.}3 hold. Then $\Phi$ is twice differentiable and
the second derivative satisfies the estimates
$$\multline
\||{\Phi}^{\prime\prime}(\delta x,\epsilon x)\||_a
\leq C_a \Big( \||\delta x\||_{a+4,1}
\||\epsilon{x}\||_{1,1}+
 \||\delta x\||_{1,1}
\||\epsilon{x}\||_{a+4,1}\Big)\\
+C_a \||x-x_0\||_{a+5,1} \||\delta x\||_{1,1}
\||\epsilon {x}\||_{1,1}\Big)
\endmultline \tag 16.6
$$
provided that
$$
\||x-x_0\||_{4,2}\leq 1 \tag 16.7
$$
Here the norms are as in (15.5).
\endproclaim

Let us now calculate the second derivative of $\Phi$
and afterwards prove the tame estimates for it.
Let us first recall the commutator identities:
\proclaim{Lemma {16.}2} We have
$$\align
[\delta,\pa_i]&=-(\pa_i\delta x^k)\pa_k\tag 16.8\\
[\delta,\pa_i\pa_j]&=-(\pa_i\delta x^k)\pa_j\pa_k-(\pa_j\delta x^k)\pa_i\pa_k
-(\pa_i\pa_j \delta x^k)\pa_k\tag 16.9
\endalign
$$
Furthermore
$$
[\delta,\triangle]=-(\triangle \delta x^k)\pa_k -2(\pa^i\delta x^j)\pa_i\pa_j\tag 16.10
$$
and if $\varepsilon$ is another variation then
$$\multline
\big[\delta,[\epsilon,\triangle]\big]
=\big((\triangle \delta x^l)\pa_l\epsilon x^k+(\pa_l\pa_m\delta x^k)\pa_l \epsilon x^m
+(\triangle \epsilon x^l)\pa_l\delta x^k+(\pa_l\pa_m\epsilon x^k)\pa_l \delta x^m\big)
\pa_k \\
2\big((\pa^k \delta x^m)\pa_m \epsilon x^l +(\pa^k \epsilon x^m)\pa_m \delta x^l
+(\pa^m \delta x^k)\pa_m \epsilon x^l\big)\pa_k\pa_l
\endmultline \tag 16.11
$$
\endproclaim
\demo{Proof} (16.8) was proven in Lemma {2.}2 and (16.9) follows from this
since $[\delta,\pa_i\pa_j]=[\delta,\pa_i]\pa_j+\pa_i [\delta,\pa_j]$.
(16.10) follows from contracting (16.9). (16.11) follows from using (16.9)
and (16.10) applied to $\delta$ as well as $\epsilon$ in place of $\delta$.
\qed\enddemo
Let $\overline{x}(t,y,r)=x(t,y)+r\, \delta x(t,y)$. The
first variational derivative $\Phi^\prime(x)$ of the Euler map
$$
\Phi^\prime(x)\delta x_i=\delta\Phi(x)_i=\frac{\pa \Phi(\overline{x})_i}{\pa
r}\big|_{r=0}\tag 16.12
$$
is given by
\proclaim{Lemma {16.}3}
$$
\Phi^\prime(x)\delta x_i=D_t^2 \delta x_i -\pa_k p\,\, \pa_i\delta x^k
+\pa_i\, p^{\,\prime}(\delta x), \tag 16.13
$$
Here $\delta p=p^{\,\prime}(\delta x)=\Psi^\prime(x)\delta x$
satisfies
$$\align
\triangle \delta p&=\delta\triangle p + \pa_k p\, \triangle \delta
x^k\! +2(\pa_i\pa_k p)\pa^i\delta x^k,\qquad\text{where}\tag 16.14\\
\delta\triangle p&=2\pa_k V^i\, \pa_i\delta x^l\,\pa_l V^k \!-2\pa_k V^i\,\pa_i
\delta v^k\tag 16.15
\endalign
$$
where $\delta v=D_t\delta x$ and $\delta p\big|_{\pa\Omega}=0$.
\endproclaim
\demo{Proof} This follows from a calculation using that
$\delta -\delta x^k\pa_k$ commutes with $\pa_i$ and hence with $\triangle$
or using (16.9).
\qed\enddemo

Let  $\overline{x}(t,y,r,s)=x(t,y)+r\delta x(t,y)+s\,\epsilon
x(t,y)$. Then the second variational derivative is given by
$$
 \Phi^{\prime\prime}(x)(\delta x,\epsilon x)_i=
 \epsilon \delta \Phi(x)_i= \frac{\pa^2\Phi_i(\overline{x})}{\pa r\pa s}
\Big|_{r=s=0},
 .\tag 16.16
$$
is given by:
\proclaim{Lemma {16.}4} Let $\delta v=D_t\delta x$ and $\epsilon v=D_t\epsilon x$. Then
$$
\Phi^{\prime\prime}(\delta x,\epsilon x)_i= \pa_k p\,\, \big(\pa_i
\epsilon x^l\,\pa_l \delta x^k +\pa_i \delta x^l\,\pa_l \epsilon
x^k\big) -\pa_k p^{\,\prime}(\epsilon x)\,\, \pa_i \delta x^k
-\pa_k p^{\,\prime}(\delta x)\,\, \pa_i \epsilon x^k +\pa_i\,
p^{\, \prime\prime}(\delta x,\epsilon x)\tag 16.17
$$
where $\delta p=p^{\,\prime}(\delta x)=\Psi^\prime(x)\delta x$ and
$\delta\epsilon p=p^{\,\prime\prime}(\delta x,\epsilon
x)=\Psi^{\prime\prime}(x)(\delta x,\epsilon x)$ satisfies
$$\triangle (\delta\epsilon p)=[\triangle,\delta\epsilon]p +\delta\epsilon\triangle p,
\qquad [\triangle,\delta\epsilon]p= f_1+2f_2-f_3-2f_4,
\quad \delta\epsilon\triangle p=-2f_5+2f_6-2f_7\tag 16.18$$
where:
$$
f_1=(\triangle \delta x^i)(\partial_i\epsilon p)+(\triangle\epsilon x^i)(\partial_i\delta p)$$
$$f_2=(\partial_i\partial_j\delta p)(\partial_j\epsilon x^i)+(\partial_i\partial_j\epsilon p)(\partial_j\delta x^i)$$
$$f_3=\partial_j p\{(\partial_i\delta x^j)(\triangle\epsilon x^i)+(\partial_i\epsilon x^j)(\triangle\delta x^i)
+2(\partial_k\delta x^i)(\partial_k\partial_i\epsilon x^j)+2(\partial_k\epsilon x^i)(\partial_k\partial_i\delta x^j)\}$$
$$f_4=\partial_i\partial_j p\{(\partial_k\delta x^j)(\partial_k\epsilon x^i)+(\partial_k\delta x^i)(\partial_j\epsilon x^k)
+(\partial_k\epsilon x^i)(\partial_j\delta x^k)\}$$
$$f_5=(\partial_k v^l)(\partial_l v^j)\{(\partial_i\delta x^k)(\partial_j\epsilon x^i)+(\partial_i\epsilon x^k)(\partial_j\delta x^i)\}
+(\partial_i v^k)(\partial_j v^l)(\partial_k\delta x^j)(\partial_l\epsilon x^i)$$
$$f_6=(\partial_k v^j)\{(\partial_i\delta v^k)(\partial_j\epsilon x^i)+(\partial_i\epsilon v^k)(\partial_j\delta x^i)
+(\partial_j\delta v^i)(\partial_i\epsilon x^k)+(\partial_j\epsilon v^i)(\partial_i\delta x^k)\}$$
$$f_7=(\partial_i\delta v^j)(\partial_j\epsilon v^i)$$
and $\delta\epsilon p\big|_{\pa\Omega}=0$.
\endproclaim
\demo{Proof}
A calculation using that $[\delta,\pa_i]=-(\pa_i\delta x^k)\pa_k$ and $\epsilon\delta x=0$
gives (16.17). (16.18) follows from using Lemma {16.}2 and
$$
\triangle \delta\epsilon p=[\delta,\triangle]\epsilon p+[\epsilon,\triangle]\delta p
+\big[\delta,[\epsilon,\triangle]\big]+\delta\epsilon \triangle p\qed\tag 16.19
$$
\enddemo

The estimates for the first and second derivative of $p=\Psi(x)$
are given in the following lemma:
\proclaim{Lemma {16.}5} Let $p=\Psi(x)$ be the solution of
$\triangle p=-(\pa_i V^j)\pa_j V^i$, $p\big|_{\pa\Omega}=0$, where
$V=D_t x$. Let $\delta p=p^{\,\prime}(\delta
x)=\Psi^\prime(x)\delta x$ be the variational derivative. We have
with $D_t\delta x=\delta v$,
 $D_t\epsilon x=\epsilon v$:
$$
\|\delta p\|_{r,\infty}\leq K_3\Big(\|\delta {v}\|_{r,\infty}+\|\delta{x}\|_{r+1,\infty}
+\big(\|{x}\|_{r+2,\infty}+\|{v}\|_{r+1,\infty}\big)
\big(\|\delta x\|_{1,\infty}+\|\delta v\|_{1,\infty}\big)\Big)\tag 16.20
$$
and with $p^{\prime\prime}(\delta x,\epsilon x)=\Psi^{\prime\prime}(x)(\delta x,\varepsilon x)$
the second variational derivative, we have
$$\multline
\|p^{\prime\prime}(\delta x,\epsilon x)\|_{r,\infty}
\leq K_3(\|\delta {v}\|_{r+1,\infty}+\|\delta{x}\|_{r+2,\infty})
(\|\epsilon x\|_{1,\infty}+\|\epsilon v\|_{1,\infty})\\
K_3(\|\epsilon {v}\|_{r+1,\infty}+\|\epsilon {x}\|_{r+2,\infty})
(\|\delta x\|_{1,\infty}+\|\delta v\|_{1,\infty})\\
+K_3(\|{v}\|_{r+2,\infty}+\|{x}\|_{r+3,\infty})
(\|\epsilon x\|_{1,\infty}+\|\epsilon v\|_{1,\infty})
(\|\delta x\|_{1,\infty}+\|\delta v\|_{1,\infty})
\endmultline \tag 16.21
$$
\endproclaim
\demo{Proof} The proof of (6.20) is similar to the estimate
of a time derivative in the proof of Lemma {6.}4.
By Lemma {16.}3, Lemma {6.}3 and Lemma {6.}4
$$\multline
\|\triangle \delta p-\delta\triangle p\|_{r-1,\infty}\leq K_1\|\delta x\|_{r+1,\infty}
\|p\|_{1,\infty}\!+\! K_1\|p\|_{r+1,\infty}\|\delta x\|_{1,\infty}
\!+\!K_1 \|x\|_{r+1,\infty} \| p\|_{1,\infty}\|\delta x\|_{1,\infty}\\
\leq  K_3\|\delta x\|_{r+1,\infty}+\big(\|v\|_{r+1,\infty}
+\|x\|_{r+2,\infty}\big)\|\delta x\|_{1,\infty}
\endmultline \tag 16.22
$$
and
$$
\|\delta\triangle p\|_{r-1,\infty}\leq
K_3\big(\|\delta x\|_{r,\infty}+\|\delta v\|_{r,\infty}\big)
+K_3\big(\|v\|_{r,\infty}+\|x\|_{r,\infty}\big)
\big(\|\delta x\|_{1,\infty}+\|\delta v\|_{1,\infty}\big)
\tag 16.23
$$
which proves (16.20). Similarly by Lemma {16.}4, Lemma {6.}3, Lemma {6.}4
and (6.20)
$$\multline
\|[\triangle,\delta\epsilon]p\|_{r-1,\infty}
\leq K_1\|\delta x\|_{r+1,\infty}
\|\epsilon p\|_{1,\infty}\!+\! K_1\| \epsilon p\|_{r+1,\infty}\|\delta x\|_{1,\infty}
\!+\!K_1 \|x\|_{r+1,\infty} \| \epsilon p\|_{1,\infty}\|\delta x\|_{1,\infty}\\
+K_1\|\delta p\|_{r+1,\infty}
\|\epsilon x\|_{1,\infty}\!+\! K_1\| \epsilon x\|_{r+1,\infty}\|\delta p\|_{1,\infty}
\!+\!K_1 \|x\|_{r+1,\infty} \| \epsilon x\|_{1,\infty}\|\delta p\|_{1,\infty}\\
+K_1\|\delta x\|_{r+1,\infty}\|\epsilon x\|_{1,\infty}\!
+\! K_1\| \epsilon x\|_{r+1,\infty}\|\delta x\|_{1,\infty}
\!+\!K_1( \|x\|_{r+1,\infty} +\|p\|_{r+1,\infty})\| \epsilon x\|_{1,\infty}\|\delta
x\|_{1,\infty}\\
\leq  K_1\big(\| \epsilon x\|_{r+2,\infty}+\| \epsilon v\|_{r+1,\infty}\big)
\|\delta x\|_{1,\infty}
+ K_1\big(\| \delta x\|_{r+2,\infty}+\| \delta v\|_{r+1,\infty}\big)
\|\epsilon x\|_{1,\infty}\\
+K_1( \|x\|_{r+3,\infty} +\|v\|_{r+2,\infty})\| \epsilon x\|_{1,\infty}\|\delta
x\|_{1,\infty}
\endmultline \tag 16.24
$$
and
$$\multline
\|\delta\epsilon\triangle p\|_{r-1,\infty}\leq
K_3\big(\|\delta x\|_{r,\infty}+\|\delta v\|_{r,\infty}\big)
\big(\|\epsilon x\|_{1,\infty}+\|\epsilon v\|_{1,\infty}\big)
\!+K_3\big(\|\epsilon x\|_{r,\infty}+\|\epsilon v\|_{r,\infty}\big)
\big(\|\delta x\|_{1,\infty}+\|\delta v\|_{1,\infty}\big)\\
\!+K_3 \big(\|v\|_{r,\infty}\!+\|x\|_{r,\infty}\big)
\big(\|\delta x\|_{1,\infty}+\|\delta v\|_{1,\infty}\big)
\big(\|\epsilon x\|_{1,\infty}+\|\epsilon x\|_{1,\infty}\big)
\endmultline
\tag 16.25
$$
which proves (16.21).
\qed\enddemo

It now follows from Lemma {16.}1, Lemma {16.}6, the fact that
$\pa_i=(\pa y^a/\pa x^i)\pa/\pa y^a$ and interpolation:
\proclaim{Lemma {16.}6}
$$\multline
\|\Phi^{\prime\prime}(\epsilon x,\delta x)_i\|_{r,\infty}
\leq K_3(\|\delta {v}\|_{r+2,\infty}+\|\delta{x}\|_{r+3,\infty})
(\|\epsilon x\|_{1,\infty}+\|\epsilon v\|_{1,\infty})\\
K_3(\|\epsilon {v}\|_{r+2,\infty}+\|\epsilon {x}\|_{r+3,\infty})
(\|\delta x\|_{1,\infty}+\|\delta v\|_{1,\infty})\\
+K_3(\|{v}\|_{r+3,\infty}+\|{x}\|_{r+4,\infty})
(\|\epsilon x\|_{1,\infty}+\|\epsilon v\|_{1,\infty})
(\|\delta x\|_{1,\infty}+\|\delta v\|_{1,\infty})
\endmultline \tag 6.26
$$
\endproclaim

Finally, also using (15.8) we get Proposition {16.}1.

\head{17. The smoothing operators.}\endhead We will work in
H\"older spaces since the standard proofs of the Nash-Moser
theorem uses H\"older spaces. The H\"older norms for functions
defined on a compact convex set $B$ are given by, if $k<a\leq
k+1$., where $k\geq 0$ is an integer,
$$
\|u\|_{a,\infty}=\| u\|_{H^a}=\sup_{x,y\in B}\sum_{|\alpha|=k} \frac{
|\pa^\alpha u(x)-\pa^\alpha u(y)|}{|x-y|^{a-k}}+ \sup_{x\in
B}|u(x)|\tag 17.1
$$
and $\|u\|_{H^0}=\sup_{x\in B} |u(x)|$.  Since we use the same
notation for the $C^k$ norms, $\|u\|_{k,\infty}=\|u\|_{C^k}$ we will differ
these by simply using letters $a,b,c,d,e,f$ etc for the H\"older
norms and $i,j,k,l, ..$ for the $C^k$ norms. However, since a Lipschitz
continuous function is differentiable almost everywhere and the
norm of the derivative at these points is bounded by the Lipschitz
constant, we conclude that for integer values this is the same as
the $L^\infty$ norm of $\pa^\alpha u$ for $|\alpha|\leq k$, and
furthermore, since all our functions are smooth it is the same as
the supremum norm. Our tame estimates for the inverse of the
linearized operator and the second variational derivative are only
for $C^k$ norms with integer exponents, with $B=\overline{\Omega}$. However, since
$\|u\|_{k,\infty}\leq C\|u\|_{a,\infty}\leq C\|u\|_{{k+1},\infty}$, if $k\leq a\leq k+1$,
see (17.2),  they also hold for non integer values with a loss of
one more derivative.

The H\"older norms satisfy
$$
\|u\|_{a,\infty}\leq C\|u\|_{b,\infty} ,\qquad a\leq b\tag 17.2
$$
and they also satisfy the interpolation inequality
$$
\| u\|_{c,\infty}
\leq C\|u\|_{a,\infty}^\lambda \|u\|_{b,\infty}^{1-\lambda}\tag 17.3
$$
where $a\leq c\leq b$, $0\leq \lambda\leq 1$ and $\lambda a+(1-\lambda)b=c$.

We will use norms which consist of H\"older norms in space
and supremum $C^k$ norms only in time
$$
\|| u\||_{a,k}=\sup_{0\leq t\leq T}\big( \| u(t,\cdot)\|_{a,\infty}+\| D_t
u(t,\cdot)\|_{a,\infty}+...+ \| D_t^k u(t,\cdot)\|_{a,\infty}\big).\tag 17.4
$$

 For the Nash-Moser technique, apart from tame estimates
one also needs smoothing operator $S_\theta$
 that satisfy the properties below with respect to the H\"older norms,
and in fact also with respect to the norms above since the smoothing
operators will be invariant under time translation. We have:
\proclaim{Proposition {17.}1}
$$\align
\|S_\theta u\|_{a,\infty}&\leq C\| u\|_{b,\infty},\qquad a\leq b\tag 17.5\\
\|S_\theta u\|_{a,\infty}&\leq C \theta^{a-b} \|u\|_{b,\infty},\qquad a\geq b \tag 17.6\\
\|(I-S_\theta) u\|_{a,\infty} &\leq C \theta^{a-b}\|u\|_{b,\infty},\quad a\leq b\tag 17.7\\
\| (S_{2\theta}-S_\theta) u\|_{a,\infty} &
\leq C\theta^{a-b} \|u\|_{b,\infty}\, \qquad a,b\geq 0.   \tag 17.8
\endalign
$$
where the constants $C$ only depend on the dimension and an upper bound for
$a$ and $b$.

Moreover, these estimates hold with the norms replaced by the norms (17.4)
for fixed $k$.
\endproclaim
First we note that (17.8) follows from (17.6), when $a\geq b$  and (17.7), when $a\leq b$.
(This alternatively follows from an additional property
$\| d/d\theta S_\theta u\|_{a,\infty}\leq C\theta^{a-b-1}\|u\|_{b,\infty}$, $a\geq 0$.
that also hold.)

For compactly supported functions on $\bold{R}^n$ there are
standard smoothing operators, see \cite{H1}, that satisfy the
above properties (17.5)-(17.8), with respect to the norms defined
in (17.1). However we have functions defined on the compact set
$\overline{\Omega}$ that do not have compact support in $\Omega$.
Therefore we need to extend these functions to have compact
support in some larger set, without increasing the H\"older norms
more than with a multiplicative constant. There is a standard
extension operator in \cite{S} that turns out to have
 these properties, see Lemma {17.}2 below.
If $\tilde{S}_\theta$ is the standard smoothing operator mentioned above, that
satisfies (17.5)-(17.8), the we define our smoothing operator
by
$$
S_\theta u= \tilde{S}_\theta \tilde{u}\big|_{\Omega}, \qquad
\text{where}\qquad  \tilde{u}={\Cal Ext}(u)\tag 17.9
$$
Since $\tilde{S}_\theta$ satisfies (17.5)-(17.8) and since
$\|\tilde{u}\|_{b,\infty}\leq C\|u\|_{b,\infty}$, by Lemma {17.}2, it follows that
$S_\theta $ satisfies (17.5)-(17.8).
\proclaim{Lemma {17.}2} There
is a linear extension operator ${\Cal Ext }$ such that ${\Cal
Ext}(f)=f$ in $\{y;\, |y|\leq 1\}$, supp ${\Cal Ext}
(f)\subset\{y;\, |y|\leq 2\}$ and
$$
\|{\Cal Ext }( f)\|_{a,\infty}\leq C\|f\|_{a,\infty}\tag 17.10
$$
where the norms in the left are H\"older norms in $\{y;\, |y|\leq 2\}$
and the norms in the right are H\"older norms in $\{y;\, |y|\leq 1\}$,
and $C$ is bounded when $a$ is bounded.
\endproclaim
\demo{Proof}
We will introduce polar coordinates and for fixed
angular variables $\omega$ extend a function defined for the radial variable
 $r\leq 1$ to $r\geq 1$. Away from the origin, the change of variables
given by polar coordinates is a diffeomorphism and H\"older continuity is
preserved under composition with a diffeomorphism $\kappa$:
$$
\|f\circ\kappa\|_{a,\infty}\leq C_a \|f\|_{a,\infty}\tag 17.11
$$
 Therefore, let us first remove the
origin by a partition of unity. Let $\chi_0\in C_0^\infty(\bold{R})$ satisfy
$\chi_0(|y|)=1$, when $|y|\leq 1/2$ and $\chi_0(|y|)=0$, when $|y|\geq
3/4$, and let $\chi_1=1-\chi_0$. Furthermore, we multiply with another cutoff
function so that the extension has compact support in $|y|\leq 2$. Let
 $\chi_2\in C_0^\infty(\bold{R})$ satisfy $\chi_2(|y|)=1$, when
$|y|\leq 5/4$ and $\chi_2(y)=0$, when $|y|\geq 3/2$. If ${\Cal
Ext}_1( f)$ is the extension operator in the radial variable,
defined in (17.14) below, we now define the extension ${\Cal Ext}(f)$ of $f$ to
be
$$
{\Cal Ext}(f)=\chi_2 {\Cal Ext}_1(\chi_1 f) +\chi_0 f\tag 17.12
$$
H\"older continuity in $(r,\omega)$ follows from
H\"older continuity of ${\Cal Ext}_1(f)$ in the radial variable
and the linearity and invariance under rotations of ${\Cal Ext }_1(f)$,
using the triangle inequality. In fact if $f_\omega(r)=f(r,\omega)$ then
$\pa_\omega^\alpha {\Cal Ext}_1(f_\omega)={\Cal Ext}_1 (f_\omega^\alpha)$,
where $f_{\omega}^\alpha=\pa_\omega^\alpha f_\omega$ and if $j+|\alpha|=k<a\leq k+1$
then by (17.18) and (17.17)
$$\multline
|\pa_r^j{\Cal Ext}_1(f_{\omega}^\alpha)(r)
-\pa_r^j{\Cal Ext}_1(f_\sigma^\alpha)(\rho)|
\leq |\pa_r^j{\Cal Ext}_1(f_\omega^\alpha)(r)
-\pa_r^ j{\Cal Ext}_1(f_\omega^\alpha)(\rho)|
+ |\pa_r^ j{\Cal Ext}_1(f_\omega^\alpha-f_\sigma^\alpha)(\rho)|\\
\leq
\sup_{r^\prime,\rho^\prime}
\frac{|\pa_r^j\pa_\omega^\alpha f(r^\prime,\omega)-\pa_r^j\pa_\omega^\alpha f(\rho^\prime,\omega)|}
{|r^\prime-\rho^\prime|^{a-k}}\,\,|r-\rho|^{a-k}
 +\sup_{\rho}\sup_{\omega^\prime,\sigma^\prime}
\frac{|\pa_r^j\pa_\omega^\alpha f(\rho,\omega^\prime)-\pa_r^j\pa_\omega^\alpha f(\rho,\sigma^\prime)|}
{|\omega^\prime-\sigma^\prime|^{a-k}}\,\, |\omega-\sigma|^{a-k}
\endmultline
\tag 17.13
$$
It therefore remains to prove the estimates (17.17) and (17.18)
for the extension in the radial variable only given by (17.14).

Suppose that $f(r)$ is a function defined for $r\leq 1$, then we
define the extension $f$ by  ${\Cal Ext}_1 (f)(r)=f(r)$, when $r\leq 1$,
and
$$
{\Cal Ext}_1(f)(r)=\int_1^\infty f(r-2\lambda (r-1))\,
\psi_1(\lambda)\, d\lambda, \qquad r\geq 1\tag 17.14
$$
where $\psi_1$ is a continuous function on $[1,\infty)$, such that
$$
\int_1^\infty \psi_1(\lambda)\, d\lambda =1, \qquad \int_1^\infty
\lambda^k \psi_1(\lambda)\, d\lambda =0,\quad k>0, \qquad
|\psi_1(\lambda)|\leq C_N (1+\lambda)^{-N}, \quad N\geq 0\tag
17.15
$$
The existence of such a function was proved in \cite{S} where the
extension operator was also introduced.
In \cite{S} it was proven that this operator is continuous on the Sobolev spaces
but it was not proven there that it is continuous on the H\"older spaces so we must
prove this.  As pointed out above, we only need to prove that it is H\"older
continuous with respect to the radial variable.

First we note that if $f\in C^k$ then the extension is in $C^k$. In fact
$$
\pa_r^j\, {\Cal Ext}_1(f)(r)=\int_1^\infty f^{(j)} (r-2\lambda
(r-1))(1-2\lambda)^j \, \psi_1(\lambda)\, d\lambda, \qquad r\geq
1\tag 17.16
$$
From the continuity of $\pa_r^j f$ and (17.14)-(17.15) it follows
that $\lim_{r\to 1}\pa_r^j\, {\Cal Ext }_1(f)(r)=\pa_r^j\, f(1)$,
that ${\Cal Ext}_1(f)$ is in $C^k$, and that for $k$ integer
$$
\sup_r |\pa_r^k {\Cal Ext }_1(f)(r)|
\leq C_k\sup_r  |f^{(k)}(r)|\tag 17.17
$$

 Suppose now that
$k<a\leq k+1$ where $k$ is an integer. We will prove that
$$
\sup_{r,\rho} \frac{|\pa_r^k\, {\Cal Ext}_1(f)(r)-\pa_r^k\, {\Cal
Ext}_1(f)(\rho)|} {|r-\rho|^{a-k}}\leq
C_a\sup_{r,\rho} \frac{|f^{(k)}(r)-f^{(k)}(\rho)|}{|r-\rho|^{a-k}} \tag 17.18
$$
If $r\leq 1$ and $\rho\leq 1$ there
is nothing to prove. Also if $r<1<\rho$ or $\rho<1<r$, then
$|r-\rho|\geq |1-\rho| $ and $|r-\rho|\geq |1-r|$ so in this case,
we can reduce it to two estimates with either $r=1$ or $\rho=1$.
Also it is symmetric in $r$ and $\rho$ so it only remains to prove
the assertion when $r>\rho\geq  1$. Then we have
$$\multline
\Big|\int_1^\infty\big( f^{(k)} (r-2\lambda (r-1))- f^{(k)}
(\rho-2\lambda (\rho-1))\big) (1-2\lambda)^k \, \psi_1(\lambda)\,
d\lambda \Big|\\
\leq \sup_{r^\prime\!,\rho^\prime}
\frac{|f^{(k)}(r^\prime)-f^{(k)}(\rho^\prime)|}{|r^\prime-\rho^\prime|^{a-k}}\,\,
|r-\rho|^{a-k}
\int_1^\infty |(1-2\lambda)^{a} \, \psi_1(\lambda)|\,
d\lambda\endmultline \tag
17.19
$$
and using the last estimate in (17.15), (17.18) follows. \qed\enddemo

\head 18. The Nash Moser Iteration.\endhead
At this point, given the results stated in sections 11-14, the problem is now reduced to a
completely standard application of the Nash-Moser technique.
One can just follow the steps of the proof of \cite{AG,H1,H2,K1}
replacing their norms with our norms. The main difference is that we have a boundary,
but we have constructed smoothing operators that satisfy the required properties
for the case with a boundary.
Furthermore, we avoid doing smoothing in the time direction, a similar approach was followed in \cite{K2}.
Alternatively, one could follow
the approach of \cite{Ha}, where it is proven that
$C^\infty$ of a compact manifold with a boundary is also a tame space,
just one small detail is missing which is that the the set $[0,T]\times\overline{\Omega}$
is not smooth at $\{0\}\times \pa\Omega$, and again we get back to the situation were it
is preferable just to do smoothing in the space directions only.

We will follow the formulation from \cite{AG} which however is
similar to \cite{H1,H2}. The theorem there is stated in terms of
H\"older norms, with a slightly different definition of the
H\"older norms for integer values. However, the only properties
that are used of the norms are the smoothing properties,
(17.5)-(17.8) and the interpolation property (17.3) which we
proved with the usual definition, i.e. the one used in \cite{H1}.
\comment
Furthermore, there is a critical lemma, Lemma 4.2.1 in \cite{AG},
for non integer values. This lemma comes from Theorem A.11 in
\cite{H1} where it is stated for $C^\infty$ functions on a compact
set $K$ and with the usual definitions of the H\"older norms, so
it holds in our situation.
\endcomment

Let us also change notation and call $\tilde{\Phi}(u)$ in (2.13) $\Phi(u)$.
Let
$$
\|| u\||_{a,k}=\sup_{0\leq t\leq T} \| u(t,\cdot)\|_{a,\infty}+...+
\| D_t^k u(t,\cdot)\|_{a,\infty},\qquad\quad \||u\||_a=\||u\||_{a,0},\tag 18.1
$$
where $\|u(t,\cdot)\|_a$ are the  H\"older norms, see (17.1).
Proposition {15.}1 and Proposition {16.}1 now says that the conditions
$({\Cal H}_1)$ and $({\Cal H}_2)$ below hold:

(${\Cal H}_1$):  $\Phi$,
is twice differentiable and satisfies
$$\multline
\||\Phi^{\prime\prime}(u)(v_1,v_2)\||_a \leq C_a\Big( \||
v_1\||_{a+\mu,2} \|| v_2\||_{\mu,2}+
 \|| v_1\||_{\mu,2}
\||v_2\||_{a+\mu,2}\Big)\\
 +C_a \|| u\||_{a+\mu,2} \||v_1\||_{\mu,2} \||v_2\||_{\mu,2},
\endmultline  \tag 18.2
$$
where $\mu=5$,
for $u,v_1,v_2\in C^\infty([0,T],C^\infty(\overline{\Omega}))$,
if
$$
\|| u\||_{\mu,2}\leq 1,\qquad \mu=5\tag 18.3
$$

(${\Cal H}_2$):
If  $u\in C^\infty([0,T],C^\infty(\overline{\Omega}))$
satisfies (18.3) then there is a linear map
 $\psi(u)$ from  $C^\infty([0,T],C^\infty(\overline{\Omega}))$
to $ C^\infty([0,T],C^\infty(\overline{\Omega}))$ such that
$\Phi^\prime(u)\psi(u) =Id$ and
$$
\||\psi(u) g\||_{a,2}\leq C_a\big( \|| g\||_{a+\lambda} +\||g
\||_\lambda\, \||u\||_{a+d,2} \big), \tag 18.4
$$
where $\lambda =[n/2]+3$ and $d=[n/2]+7$.
\proclaim{Proposition {18.}1} Suppose that $\Phi$ satisfies (${\Cal H}_1$),
 (${\Cal H}_2$) and $\Phi(0)=0$. Let $\alpha > \mu$,
$\alpha> d$, $\alpha>\lambda+2\mu$,
$\alpha\notin \Bbb{N}$.
Then

i) There is neighborhood $W_\delta=\{ f\in C^\infty([0,T],C^\infty(\overline{\Omega}))
; \, \|| f\||_{\alpha+\lambda}\leq \delta^2\}$, $\delta>0$,
such that, for $f\in W_\delta $,
the equation
$$
\Phi(u)=f\tag 18.5
$$
has a solution $u=u(f)\in C^2([0,T],C^\infty(\overline{\Omega}))$.
Furthermore,
$$
\|| u(f)\||_{a,2}\leq C\||f\||_{\alpha+\lambda},\qquad a<\alpha\tag 18.2
$$
\endproclaim

In the proof, we construct a sequence $u_j\in C^\infty([0,T],C^\infty(\overline{\Omega}))$
converging to $u$, that satisfy
 $\|| u_j\||_{\mu,2}\leq 1$ and $\||S_i u_i\||_{\mu,2}\leq 1$, for all $j$,
 where $S_i$ is the smoothing operator in (18.7).
The estimates (18.2) and (18.4) will only be used for convex
combinations of these and hence within the domain (18.3) for which
these estimates hold.

Following \cite{H1,H2,AG,K1,K2} we set
$$
u_{i+1}=u_i+\delta u_i,\qquad \delta u_i=\psi(S_i u_i) g_i,  \quad u_0=0,
\qquad
S_i=S_{\theta_i}, \quad \theta_i=\theta_0 2^i,\quad \theta_0\geq 1\tag 18.7
$$
and $g_i$ are to be defined so that $u_i$ formally converges to a solution.
We have
$$\multline
\Phi(u_{i+1})-\Phi(u_i)=\Phi^\prime(u_i)(u_{i+1}-u_i)+e^{\prime\prime}_i
=\Phi^\prime(u_i)\psi(S_i u_i) g_i +e^{\prime\prime}_i\\
=(\Phi^\prime(u_i)-\Phi^\prime(S_i u_i))\psi(S_i u_i)
g_i+g_i+e^{\prime\prime}_i =e^\prime_i+e^{\prime\prime}_i+g_i
\endmultline \tag 18.8
$$
where
$$\align
e_i^\prime&=(\Phi^\prime(u_i)-\Phi^\prime(S_i u_i))\delta u_i\tag 18.9 \\
e_i^{\prime\prime}&=\Phi(u_{i+1})-\Phi(u_i)-\Phi^\prime(u_i)\delta u_i \tag 18.10\\
e_i&=e^\prime_i+e^{\prime\prime}_i\tag 18.11
\endalign
$$
Therefore
$$
\Phi(u_{i+1})-\Phi(u_i)=e_i+g_i\tag 18.12
$$
and adding, we get
$$
\Phi(u_{i+1})=\sum_{j=0}^i g_j+S_i E_i +e_i +(I-S_i)E_i,
\qquad E_i=\sum_{j=0}^{i-1} e_j\tag 18.13
$$
To ensure that $\Phi(u_i)\to f$ we must have
$$
\sum_{j=0}^i g_j+S_i E_i=S_i f\tag 18.14
$$
Thus
$$
g_0=S_0 f, \qquad g_i=(S_{i}-S_{i-1})(f-E_{i-1})-S_i e_{i-1}\tag 18.15
$$
and
$$
\Phi(u_i)=S_i f+e_i +(I-S_i) E_i\tag 18.16
$$
Given $u_0,u_1,...,u_i$ these determine $\delta u_0,\delta u_1,...,\delta u_i$
which by (18.9)-(18.10) determine $e_1,...,e_{i-1}$, which by (18.15)
determine $g_i$. The new term $u_{i+1}$ is the determined by (18.7).

\proclaim{Lemma {18.}2} Assume that $\|| u_i\||_{\mu,2}\leq 1$,
$\|| u_{i+1}\||_{\mu,2}\leq 1$ and $\|| S_i u_i\||_{\mu,2}\leq 1$.
Then
$$\multline
\|| e_i^\prime\||_{a}\leq C_a
\Big( \||(I-S_i) u_i\||_{a+\mu,2}\||\delta u_i\||_{\mu,2}
+\||(I-S_i) u_i\||_{\mu,2}\||\delta u_i\||_{a+\mu,2}\Big)\\
+C_r\||S_i u_i\||_{a+\mu,2}
 \||(I-S_i) u_i\||_{\mu,2} \||\delta u_i\||_{\mu,2}
\endmultline\tag 18.17
$$
and
$$
\|| e_i^{\prime\prime}\||_a\leq C_r \Big( \||\delta
u_i\||_{a+\mu,2}\||\delta u_i\||_{\mu,2}
+\||u_i\||_{a+\mu,2}\||\delta u_i\||_{\mu,2}^2\Big) \tag 18.18
$$
\endproclaim
\demo{Proof} The proof of (18.17) makes use of
$$
(\Phi^\prime(u_{i})-\Phi^\prime(S_i u_i))\delta u_i =
\int_0^1  \Phi^{\prime\prime}(S_i u_i+s(I-S_i)u_i)(u_i-S_i u_i ,\delta u_i)\, ds
\tag 18.19
$$
together with (18.2). Note that from the third term in (18.2) we get a term that
is not present in (18.17) since it can be bounded by the others using the
assumptions. In fact, since
$\|| u_i\||_{\mu,2}+\|| S_i u_i\||_{\mu,2}\leq 2$,
$\||(I-S_i )u_i\||_{a+\mu,2} \||(I-S_i) u_i\||_{\mu,2} \||\delta u_i\||_{\mu,2}
\leq 2 \||(I-S_i )u_i\||_{a+\mu,2} \||\delta u_i\||_{\mu,2}$.
 (18.18) makes use of
$$
\Phi(u_{i+1})-\Phi(u_i)-\Phi^\prime(u_i)\delta u_i =
\int_0^1 (1-s) \Phi^{\prime\prime}(u_i+s\delta u_i)(\delta u_i,\delta u_i)\, ds
\tag 18.20
$$
together with (18.2). Here we used that
$\||\delta u_i\||_{a+\mu,2}
\||\delta u_i\||_{\mu,2}^2\leq 2\||\delta u_i\||_{a+\mu,2}\||\delta u_i\||_{\mu,2}
$
\qed\enddemo

Let $\tilde{\alpha}>\alpha$ and $\tilde{\alpha}-\mu>2(\alpha-\mu)$.
 Throughout the proof $C_a$ will stand for constants that
depend on $a$ but is independent of $n$ in (18.21).

Our inductive assumption $(H_n)$ is,
$$
\|| \delta u_i\||_{a,2}\leq \delta \theta_i^{a-\alpha} ,
\qquad 0\leq a\leq \tilde{\alpha},\qquad i\leq n\tag 18.21
$$
If $n=0$ then if $a\leq \tilde{\alpha}$, we have
 $\|| \delta u_0\||_{a,2}\leq C_{\tilde{\alpha}}\|| f\||_{\alpha+\lambda}\leq C_{\tilde{\alpha}} \delta^2$,
so it follows that (18.21) hold for $n=0$ if we choose $\delta$ so small that
$C_{\tilde{a}} \delta \leq \theta_0^{\tilde{\alpha}-\alpha}$.
We must now prove that $(H_n)$ implies $(H_{n+1})$ if $C_{\tilde{\alpha}}^\prime\delta \leq 1$,
where $C_{\tilde{\alpha}}^\prime$ is some
 constant that only depends on $\tilde{\alpha}$ but is
independent of $n$.

\proclaim{Lemma {18.}3} If (18.21) hold then for $i\leq n$
 $$
 \sum_{j=0}^i\||\delta u_j\||_{a,2}\leq
C_a \delta
\big(\min(i,1/|\alpha-a|)+1)(\theta_i^{a-\alpha}+1\big), \qquad
0\leq a\leq \tilde{\alpha}\tag 18.22
$$
\endproclaim
\demo{Proof} Using (18.21) we get $\sum_{j=0}^i\||\delta
u_j\||_{a,2}\leq C_a\delta \sum_{j=0}^i 2^{j(a-\alpha)}$ and
noting that $\sum_{j=0}^i 2^{-sj}\leq C(\min{(1+1/s,i)}+1)$, if
$s>0$, (18.22) follows. \qed\enddemo

\proclaim{Lemma {18.}4} If $(H_n)$, i.e. (18.21),  hold and
$\tilde{\alpha}>\alpha$, then for $i\leq n+1$ we have
$$\align
\|| u_i\||_{a,2}&\leq C_a \delta (\min(i,1/|\alpha-a|)+1)(\theta_i^{a-\alpha}+1),
\qquad 0\leq a\leq \tilde{\alpha}
\tag 18.23\\
\|| S_i u_i\||_{a,2}&\leq  C_a\delta
 (\min(i,1/|\alpha-a|)+1)(\theta_i^{a-\alpha}+1) ,
\qquad a\geq 0 \tag 18.24\\
\|| (I-S_i) u_i\||_{a,2}&\leq  C_a\delta \theta_i^{a-\alpha},
\qquad 0\leq a\leq \tilde{\alpha} \tag 18.25
\endalign
$$
\endproclaim
\demo{Proof} The proof of (18.23) is just summing up the series
$u_{i+1}=\sum_{j=0}^i \delta u_j$, using Lemma {18.}3. (18.24)
follows from (18.22) using (17.5) for $a\leq \tilde{\alpha}$ and
(17.6) with $b=\tilde{\alpha}$ for $a\geq \tilde{\alpha}$. (18.25)
follows from (17.7) with $b=\tilde{\alpha}$ and (18.23) with
$a=\tilde{\alpha}$. \qed\enddemo

Since we have assumed that $\alpha>\mu$, we note
  that in particular, it follows that
 $$
 \||u_i\||_{\mu,2}\leq 1\quad \text{and}\quad
 \||S_i u_i\||_{\mu,2}\leq 1,\qquad
 \text{for}\quad i\leq n+1\quad\text{if}\quad
  C_{\mu} \delta\leq 1.
 \tag 18.26
 $$

As a consequence of Lemma {18.}4 and Lemma {18.}2 we get
\proclaim{Lemma {18.}5} If $(H_n)$ is satisfied and $\alpha>\mu$,
then for $i\leq n$,
$$
\align
\|| e_i^\prime\||_a\leq C_a\delta^2 \theta_i^{a-2(\alpha-\mu)},
\qquad 0\leq a\leq \tilde{\alpha}-\mu \tag 18.27\\
\|| e_i^{\prime\prime}\||_a\leq C_a\delta^2
\theta_i^{a-2(\alpha-\mu)}, \qquad 0\leq a\leq \tilde{\alpha}
-\mu\tag 18.28
\endalign
$$
\endproclaim

As a consequence of Lemma {18.}5 and (17.8) we get \proclaim{Lemma
{18.}6} If $(H_n)$ is satisfied, then for $i\leq n+1$,
$$\align
\||\ S_i e_{i-1}\||_a &\leq  C_a\delta^2 \theta_i^{a-2(\alpha-\mu)},
\qquad a\geq 0\tag 18.29\\
\|| (S_{i}-S_{i-1}) f\||_a &\leq C_a\theta_i^{a-\beta} \|| f\||_\beta,
\qquad a\geq 0 \tag 18.30\\
\||(I-S_i) f\||_a &\leq C_a\theta_i^{a-\beta} \|| f\||_\beta,
\qquad 0\leq a\leq \beta  \tag 18.31
\endalign
$$
Furthermore, if $\tilde{\alpha}-\mu>2(\alpha-\mu)$:
$$\align
\|| (S_{i}-S_{i-1}) E_{i-1}\||_a &\leq  C_a\delta^2
\theta_i^{a-2(\alpha-\mu)},
\qquad a\geq 0\tag 18.32\\
\||(I-S_i) E_{i}\||_a &\leq  C_a\delta^2
\theta_i^{a-2(\alpha-\mu)}, \qquad  0\leq a\leq \tilde{\alpha}-\mu
\tag 18.33
\endalign
$$
\endproclaim
\demo{Proof} (18.29) follows from (18.27); For $a\leq
\tilde{\alpha}-\mu$ we use (17.5) with $b=a$ and for $a\geq
\tilde{\alpha}-\mu$, we use (17.6) with $b=\tilde{\alpha}-\mu$.
(18.30) follows from (17.8) and (18.31) follows from (17.7). Now,
$E_i=\sum_{j=0}^{i-1} e_j$ so by Lemma {18.}5
 $\|| E_i\||_{\tilde{\alpha}-\mu} \leq C_a\delta^2 \sum_{j=0}^{i-1}
 \theta_j^{\tilde{\alpha}-\mu-2(\alpha-\mu)}\leq
 C_a^\prime \delta^2
 \theta_i^{\tilde{\alpha}-\mu-2(\alpha-\mu)}$,
 since we assumed that the exponent is positive.
(18.32) follows from this and (17.8) with $b=\tilde{\alpha}-\mu$
and similarly (18.33) follows from (17.7) with $b=\tilde{\alpha}-\mu$.
 \qed\enddemo

It follows that: \proclaim{Lemma {18.}7} If $(H_n)$ is satisfied,
$\tilde{\alpha}-\mu>2(\alpha-\mu)$, and $\alpha>\mu$ then for
$i\leq n+1$,
$$
\|| g_i\||_a\leq C_a\delta^2\theta_i^{a-2(\alpha-\mu)} +C_a
\theta_i^{a-\beta}\||f\||_\beta ,\qquad a\geq 0. \tag 18.34
$$
\endproclaim

 Using this lemma and (18.4) we get
 \proclaim{Lemma {18.}8} If
$(H_n)$ holds, $\tilde{\alpha}-\mu>2(\alpha-\mu)$, $\alpha>\mu$,
$\alpha>d$ then, for  $i\leq n+1$, we have
$$
\|| \delta u_i\||_{a,2} \leq C_a\delta^2
\theta_i^{a+\lambda-2(\alpha-\mu)} +C_a \|| f\||_\beta
\theta_i^{a+\lambda-\beta} , \qquad a\geq 0. \tag 18.35
$$
\endproclaim
\demo{Proof} Using (18.7), (18.4), (18.34) and (18.24) we get
$$\multline
\|| \delta u_i\||_{a,2}\leq
C_a\big(\delta^2 \theta_i^{a+\lambda-2(\alpha-\mu)} +\|| f\||_\beta
\theta_i^{a+\lambda-\beta}\big)\\
+ C_a\big(\delta^2 \theta_i^{\lambda-2(\alpha-\mu)}
+\|| f\||_\beta \theta_i^{\lambda-\beta}\big)
 \delta
 (\min(i,1/|\alpha-a-d|)+1)(\theta_i^{a+d-\alpha}+1)
\endmultline \tag 18.36
$$
The lemma follows from using that
$\min{(i,1/|\alpha-a-d|)}+1\leq C\theta_i^a/(\theta_i^{a+d-\alpha}+1)$,
where $C$ is a constant depending on  $\alpha-d>0$
but independent of $i$. \qed\enddemo

If, we now pick $\beta=\alpha+\lambda$, and use the assumptions  that
$\lambda+\alpha<2(\alpha-\mu)$, and $\|| f\||_{\alpha+\lambda}\leq \delta^2$,
 we get that for $i\leq n+1$,
$$
\|| \delta u_i\||_{a,2}\leq C_a \delta^2 \theta_i^{a-\alpha},\qquad
a\geq 0, \tag 18.37
$$
If we pick $\delta>0$ so small that
$$
C_{\tilde{\alpha}}\delta \leq 1,\tag 18.38
$$
the assumption $(H_{n+1})$ is proven.

The convergence of the $u_i$ is an immediate consequence of Lemma
{18.}2:
$$
\sum_{i=0}^\infty \||u_{i+1}-u_i\||_{a,2}\leq C_{a}\delta ,\qquad
a<\alpha \tag 18.39
$$
It follows from Lemma {18.}6 that
$$
\|| \Phi(u_i)-f\||_a\leq C_a \delta^2 \theta_i^{a-\alpha-\lambda}
\tag 18.40
$$
which tends to $0$, as $i\to\infty$, if $a<\alpha+\lambda$.

It remains to prove $u\in C^2([0,T],C^\infty(\overline{\Omega}))$.
Note that in Lemma {18.}8 we proved a better estimate than
$(H_n)$. In fact if we let
$\gamma=2(\alpha-\mu)-(\alpha+\lambda)>0$ and
$\alpha^\prime=\alpha+\gamma$, then $\||
f\||_{\alpha^\prime+\lambda}\leq C$ implies that
$$
\|| \delta u_i\||_{a,2}\leq C_a \theta_i^{a-\alpha^\prime},\qquad
a\geq 0\tag 18.41
$$
Using this new estimate, in place of $(H_n)$, we can go back to
Lemma {18.}3-Lemma {18.}8 and replace $\alpha$ by $\alpha^\prime$
and $\delta $ by $1$. Then it follows from Lemma {18.}8 that
$$
\|| \delta u_i\||_{a,2}\leq
C_a\theta_i^{a+\lambda-2(\alpha^\prime-\mu)} +C_a
\theta_i^{a+\lambda-\beta}\|| f\||_\beta\tag 18.42
$$
and if we now pick $\gamma^\prime=2(\alpha^\prime-\mu)-(\lambda-\alpha^\prime)=2\gamma$
and $\alpha^{\prime\prime}=\alpha^\prime+\gamma^\prime=\alpha+2\gamma$, and use that
$\|| f\||_{\alpha^\prime+\gamma^\prime}\leq C$ we see that
$$
\|| \delta u_i\||_{a,2}\leq C_a
\theta_i^{a-\alpha^{\prime\prime}},\qquad a\geq 0\tag 18.43
$$
Since the gain $\gamma>0$ is constant, repeating this process
yields that (18.41) holds for any $\alpha^\prime$ and hence that
(18.39)-(18.40) hold for any $\alpha\geq 0$, (with $\delta$
replaced by $1$). It follows that $u_j$ is a Cauchy sequence in
$C^2\big([0,T],C^k(\overline{\Omega})\big)$, for any $k$, and
hence that $u_j\to u\in
C^2\big([0,T],C^\infty(\overline{\Omega})\big)$ and $\Phi(u_j)\to
f\in C\big([0,T],C^\infty(\overline{\Omega})\big)$. (18.6) follows
from (18.37) with $\delta^2=\||f\||_{\alpha+\lambda}$.  This
concludes the proof of Proposition {18.}1.

\subheading{Acknowledgments} I would like to thank Demetrios
Christodoulou, Richard Hamilton and Kate Okikiolu for many long and helpful discussions.

\enddocument